\documentclass[11pt,leqno]{article}
\usepackage{amsthm,amsfonts,amssymb,epsfig,graphics,amsmath,eufrak}
 \usepackage[latin1]{inputenc}\relax

\newcommand{\bDelta}{{\bar{\Delta}}}
\newcommand{\bU}{{\bar{u}}}
\newcommand{\sgn}{{\text{\rm sgn}}}
\newcommand{\pv}{{\text{\rm P.V.}}}

\newcommand{\ft}{{\cal{F}}}
\newcommand{\lt}{{\cal{L}}}
\newcommand{\ess}{{\text{\rm ess }}}
\newcommand{\calE}{{\cal{E}}}
\newcommand{\CalE}{{\cal{E}}}
\newcommand{\CalV}{{\cal{V}}}
\newcommand{\CalD}{{\cal{D}}}

\newcommand{\CalA}{{\cal{A}}}
\newcommand{\CalB}{{\cal{B}}}
\newcommand{\CalO}{{\cal{O}}}

\newcommand{\CalT}{{\cal{T}}}
\newcommand{\bfO}{{\bf O}}
\newcommand{\btheta}{{\bar \theta}}
\newcommand{\Txi}{{\tilde \xi}}
\newcommand{\Tmu}{{\tilde \mu}}
\newcommand{\hU}{{\hat U}}
\newcommand{\R}{\Re}
\newcommand{\const}{\text{\rm constant}}
\newcommand{\osc}{\text{\rm oscillation }}
\newcommand{\bdiag}{\text{\rm block-diag }}

\newcommand{\V}{\bold{v}}

\newcommand{\RR}{{\mathbb R}}
\newcommand{\WW}{{\mathbb W}}

\newcommand{\ZZ}{{\mathbb Z}}

\newcommand{\CC}{{\mathbb C}}
\newcommand{\mA}{{\mathbb A}}

\newcommand{\II}{{\mathbb I}}

\newcommand{\FF}{{\mathbb F}}
\newcommand{\GG}{{\mathbb G}}

\newcommand\adots{\mathinner{\mkern2mu\raise1pt\hbox{.}
\mkern3mu\raise4pt\hbox{.}\mkern1mu\raise7pt\hbox{.}}}

\newcommand{\rank}{{\rm rank }}
\newcommand{\range}{{\rm range }}

\renewcommand{\div}{{\rm div}}

\newcommand{\Span}{\text{\rm Span\ }}

\newtheorem{theo}{Theorem}[section]
\newtheorem{prop}[theo]{Proposition}
\newtheorem{cor}[theo]{Corollary}
\newtheorem{lem}[theo]{Lemma}
\newtheorem{defi}[theo]{Definition}
\newtheorem{ass}[theo]{Assumption}
\newtheorem{assums}[theo]{Assumptions}
\newtheorem{obs}[theo]{Observation}

\newtheorem{exam}[theo]{Example}
\newtheorem{rem}[theo]{Remark}
\newtheorem{rems}[theo]{Remarks}
\newtheorem{exams}[theo]{Examples}
\newtheorem{ex}[theo]{Exercise}
\newtheorem{exs}[theo]{Exercises}

\numberwithin{equation}{section}

 \title{     
\bf Planar stability criteria for 
viscous shock waves of systems with real viscosity
}

\author{\sc \small 
Kevin Zumbrun\thanks{Indiana University, Bloomington, IN 47405;
kzumbrun@indiana.edu:
K.Z. thanks CIME and especially organizers C.M. Dafermos and
P. Marcati for the opportunity to participate
in the summer school at which this material was originally presented.
Each section corresponds to a single $90$-minute lecture.
Thanks to G. M\'etivier, and M. Williams for their interest in
the work and for many helpful conversations,
and to O. Gues,  G. M\'etivier, and M. Williams for their indirect contribution 
through our concurrent joint investigations of the closely
related small-viscosity problem for real viscosity systems [GMWZ.4].
We note in particular that the simplified approach of 
obtaining high-frequency resolvent bounds entirely through Kawashima-type 
energy estimates was suggested to us by some small-viscosity investigations of
M. Williams; indeed, the argument given here is a large-amplitude
version of an argument developed by him in an earlier, since discarded
version of [GMWZ.4] (see also related constant-coefficient 
and one-dimensional analyses in [KSh] and [HuZ], respectively).
The novelty of the current presentation lies rather in the 
development of a nonlinear iteration scheme depending only on 
such bounds (standard for the small-viscosity problem,
new for long-time stability).
Thanks also to B. Texier for his careful reading and many helpful suggestions.
Research of K.Z. was partially supported
under NSF grants number DMS-0070765 and DMS-0300487.
 }}

\begin{document}

\maketitle

\begin{abstract}
We present a streamlined account of recent developments
in the stability theory for planar viscous shock waves,
with an emphasis on applications
to physical models with ``real,'' or partial viscosity.
The main result is the establishment of necessary, or ``weak'',
and sufficient, or ``strong'', conditions for nonlinear stability 
analogous to those established by Majda [Ma.1--3] in the inviscid case
but (generically) 
separated by a codimension-one set in parameter space rather than
an open set as in the inviscid case.
The importance of codimension one is that transition between
nonlinear stability and instability is thereby determined,
lying on the boundary set between the open regions of
strong stability and strong instability (the latter defined
as failure of weak stability).
Strong stability holds always for small-amplitude shocks 
of classical ``Lax'' type [PZ.1--2, FreS];
for large-amplitude shocks, however, strong instability may occur [ZS, Z.3].
\end{abstract}

\clearpage
\tableofcontents

\clearpage
\bigbreak

\section{Introduction: structure of physical equations}

Many equations of physics take the form of 
{\it hyperbolic conservation laws}
\begin{equation}
\label{hyperbolic}
U_t + \sum_j F^j(U)_{x_j}=0,
\end{equation}
with associated {\it viscous conservation laws}
\begin{equation}
\label{viscous}
U_t + \sum_j F^j(U)_{x_j}=\nu\sum_{j,k} (B^{jk}(U)U_{x_k})_{x_j}
\end{equation}
incorporating neglected transport effects of viscosity, heat conduction, etc.
(more generally, hyperbolic and viscous
{\it balance laws}\footnote
{ Outside the scope of these lectures, but accessible to the same
techniques; see, e.g., [God, Z.3, MaZ.1, MaZ.5, Ly, LyZ.1--2, JLy.1--2].}
including also zero-order
derivative terms $C(U)$, as especially in relaxation and combustion
equations).
Here, $U$, $F^j \in \RR^n$, $B^{jk}\in \RR^{n\times n}$, $x\in \RR^d$, and $t\in \RR$.

\begin{exams}
\textup{
Euler and Navier--Stokes equations, respectively, of 
gas- or magnetohydrodynamics (MHD);
see \eqref{nseq} below.
}
\end{exams}

A fundamental feature of \eqref{hyperbolic} is the appearance
of {\it shock waves}
\begin{equation}
\label{shock}
U(x,t)={\bar {\bar U}}(x-st)=
\begin{cases}
U_- & x_1<st,\\
U_+ & x_1\ge st,\\
\end{cases}
\end{equation}
discontinuous weak, or distributional solutions of \eqref{hyperbolic}
determined by the Rankine--Hugoniot conditions
\begin{equation}
\tag{RH}
s[U]=[F(U)];
\end{equation}
see, e.g., [La.1--2, Sm].
Here and elsewhere, $[h(U)]:= h(U_+)-h(U_-)$.
Solution \eqref{shock} may be uniquely identified
by the ``shock triple'' $(U_-,U_+,s)$.

Such waves in fact occur in applications, and do
well-approximate experimentally observed behavior.
However, in general they occur in only one direction,
despite the apparent symmetry $(x,t,s) \to (-x,-t,s)$
in equations \eqref{hyperbolic}.
That is, only one of the shock triples $(U_-,U_+,s)$ 
$(U_+,U_-,s)$ is typically observed, though they are
indistinguishable from the point of view of (RH).
The question of when and why a particular shock triple
is physically realizable, known as the {\it shock
admissibility problem}, is one of the oldest and most
central problems in the theory of shock waves.
For an interesting discussion of this issue from a general
and surprisingly modern point of view, see the 1944
roundtable discussion of [vN].

Two basic approaches to admissibility are:

1. Hyperbolic stability in the Hadamard sense, i.e., short-time 
bounded stability, or well-posedness of \eqref{shock}
as a solution of \eqref{hyperbolic}, 
also known as {\it dynamical stability} [BE]:
that is, internal consistency of the hyperbolic model.
Here, there exists a well-developed theory; see, e.g.,
[Ma.1--3, M\'e.1--4, FM\'e] and references therein. 

2. Consistency with viscous or other regularization, in this case
the viscous conservation law \eqref{viscous}. 
(a) A simple version is the ``viscous profile condition'',
requiring existence of an associated family of nearby traveling-wave solutions
\begin{equation}
\label{profile}
U(x,t)=\bar U\Big(\frac{x-st}{\nu}\Big),
\qquad \lim_{z\to \pm \infty}\bar U(z)=U_\pm
\end{equation}
of the viscous conservation law \eqref{viscous};
see, e.g., [Ra, Ge, CF, Be, Gi, MP, Pe, MeP], and references therein.
This is the planar version of the prepared-data ``vanishing-''
or ``small-viscosity'' problem (SV) treated for general, curved shocks
in the parallel article of Mark Williams in this volume [W].
The viscous profile condition, augmented with the requirement that
$\bar U$ be a transverse connection with respect to the associated
traveling-wave ODE,
is sometimes known as {\it structural stability} [BE, ZS, Z.3--4].
(b) A more stringent version is the ``stable viscous profile condition'',
requiring stability under perturbation of individual profiles \eqref{profile}
with viscosity coefficient $\nu$ held fixed.  
We denote by (LT) the associated problem of 
determining long-time viscous stability.

\begin{defi}
\label{LongT}
Long-time viscous stability is defined 
as the property that, 
for some appropriately chosen norms $|\cdot|_X$ and $|\cdot|_Y$,
for initial data $U_0$ suffiently close to profile $\bar U$
in $|\cdot|_X$, the viscous problem \eqref{viscous},
$\nu=1$, has a (unique) global solution $U(\cdot, t)$ that converges to
$\bar U$ as $t\to \infty$ in $|\cdot|_Y$.
We refer to the latter property as asymptotic $|\cdot|_X\to |\cdot|_Y$
viscous stability.
\end{defi}

\begin{rem}
\label{SVrem}
\textup{
One may also consider the question whether solutions
of \eqref{viscous} converge on a bounded time interval
for fixed initial data as $\nu \to 0$ to a solution of \eqref{hyperbolic},
that is, the unprepared-data (SV) problem.
This was considered for small-amplitude shock waves in one dimension
by Yu [Yu], and, more recently, for
general small-variation solutions in one dimension in the
fundamental work of Bianchini and Bressan [BB]; in multiple dimensions
the problem remains completely open.
This is a more stringent requirement than either of 2(a) or 2(b);
indeed, it appears to be overly restrictive as an admissibility
condition.  
In particular, as discussed in [Fre.1--2, FreL, L.4, Z.3--4],
there arise in (MHD) certain nonclassical ``overcompressive'' shocks
that are both stable for fixed $\nu$ and play an important role in solution
structure, yet which do not persist as $\nu \to 0$.
The (LT) and (SV) problems are related by the scaling
\begin{equation}
\label{scaling}
(x,t,\nu)\to (x/T, t/T,\nu/T),
\end{equation}
with $T\to \infty$, $0\le t\le T$,
the difference lying in the prescription of initial data.
}
\end{rem}
\medbreak
Conditions 1 and 2(a) may be formally derived by
matched asymptotic expansion using the rescaling \eqref{scaling}
and taking the zero-viscosity limit,
as described, e.g., in [Z.3], Section 1.3.
They have the advantage of simplicity, and for this reason
have received the bulk of the attention in
the classical mathematical physics literature; see, e.g., the
excellent surveys [BE] and [MeP].
However, rigor (and also rectitude; see Section 1.4, [Z.3]) of the theory
demands the study of the more complicated,
but physically correct condition 2(b), 
motivating the study of the long-time viscous stability problem (LT).
It is this problem that we shall consider here.

In contrast to the hyperbolic stability theory, progress in
the multidimensional viscous stability theory has come
only quite recently.
The purpose of this article is to present an account of these
recent developments, with an emphasis on
(i) connections with, and refinement of the hyperbolic stability theory,
and (ii) applications to situations of physical interest, i.e.,
{\it real viscosity}, {\it large amplitude}, and 
{\it real} (e.g., van der Waals-type) {\it gas equation of state}.
Our modest goal is to present sharp and (at least numerically)
computable planar stability criteria analogous to the 
Lopatinski condition obtained by Majda [K, Ma.1--3] in the hyperbolic case.

This is only the first step toward a complete theory;
in particular, evaluation of the stability 
criteria/classification of stability remain important open problems.
Preliminary results in this direction include
stability of general small-amplitude shock profiles
[Go.1--2, HuZ, PZ.1--2, FreS];
geometric conditions for stability, yielding instability of certain 
large-amplitude shock profiles [GZ, BSZ, FreZ, God, Z.2--4, 
Ly, LyZ.1--2];
and the development of efficient algorithms for numerical testing
of stability [Br.1--2, BrZ, BDG, KL].
On the other hand, the techniques we use here
are completely general, applying
also to relaxation, combustion, etc.; see, e.g., [Z.3] and references therein.

We begin in this section
with some background discussion of a mainly historical nature
concerning the common structural properties relevant to our investigations
of various equations arising in mathematical physics,
at the same time introducing some basic energy estimates 
of which we shall later make important use.

\medbreak
{\bf 1.1. Symmetry and normal forms.}
We may write \eqref{hyperbolic} and \eqref{viscous} in quasilinear form as
\begin{equation}
\label{quasihyperbolic}
U_t + \sum_j A^j U_{x_j}= 0
\end{equation}
and
\begin{equation}
\label{quasiviscous}
U_t + \sum_j A^j U_{x_j}=
\nu\sum_{j,k} (B^{jk} U_{x_k})_{x_j},
\end{equation}
where $A^j:=dF^j(U)$ and $B^{jk}:=B^{jk}(U)$.
We assume the further structure
\begin{equation}
\label{hypparab}
U=\begin{pmatrix}
u^I \\ u^{II}
\end{pmatrix},
\qquad
A^j=\begin{pmatrix}
A^j_{11} & A^j_{12} \\ 
A^j_{21} & A^j_{22} \\ 
\end{pmatrix},
\qquad
B^{jk}=\begin{pmatrix}
0&0\\
b^{jk}_{I} & b^{jk}_{II} \\ 
\end{pmatrix},
\end{equation}
$u^I\in \RR^{n-r}$, $u^{II}\in \RR^r$
typical in physical applications,
identifying a distinguished, ``inviscid'' variable $u^I$,
with
\begin{equation}
\label{goodb}
\R \sigma \sum \xi_j \xi_k b^{jk}_{II} \ge \theta |\xi|^2,
\end{equation}
$\theta>0$, for all $\xi\in \RR^d$.
Here and below, $\sigma M$ denotes spectrum of a matrix or linear
operator $M$. 
In the case of an unbounded operator, we use
the simplest definition of spectrum 
as the complement
of the resolvent set $\rho(M)$, defined as the set of $\lambda\in \CC$
for which $\lambda-M$ possesses a bounded inverse with respect to
a specified norm and function space; see, e.g., [Kat], or Appendix A.
\medbreak

{\bf 1.1.1. Inviscid equations.}
Local stability of constant solutions (well-posedness)
requires hyperbolicity of \eqref{quasihyperbolic},
defined as the property that $A(\xi):= \sum_j A^j \xi_j$
have real, semisimple eigenvalues for all $\xi \in R^d$.
As pointed out by Godunov and Friedrichs [G,Fr],
this may be guaranteed by {\it symmetrizability}, defined as existence of
a ``symmetrizer'' $\tilde A^0$ such that $\tilde A^0$ is symmetric
positive definite and $\tilde A^j:= \tilde A^0A^j$ are symmetric.
Left-multiplication by $\tilde A^0$ converts \eqref{quasihyperbolic} to
(quasilinear) {\it symmetric hyperbolic form}
\begin{equation}
\label{symmhyp}
\tilde A^0 U_t + \sum_j \tilde A^j U_{x_j}=0.
\end{equation}

A nonlinear version of this procedure is an invertible coordinate
change $U\to W$ such that \eqref{hyperbolic} considered as an
equation in $W$ takes the form
\begin{equation}
\label{nonlinearsymmhyperbolic}
U(W)_t + \sum_j F^j(U(W))_{x_j}=0,
\end{equation}
where $\tilde A^0:=\partial U/\partial W$ is symmetric positive definite
and $\tilde A^j:= dF^j (\partial U/\partial W)$ are symmetric.
This is more restrictive, but has the advantage of preserving divergence form;
see Remarks \ref{hentropyrems} 1-2 below.

Symmetric form yields hyperbolic properties directly
through elementary energy estimates/integration by parts,
using the {\it Friedrichs symmetrizer relation}
\begin{equation}
\label{symmrel}
\R \langle U, S U_{x_j} \rangle=
-\frac{1}{2} \langle U, S^j_{x_j}U \rangle
\end{equation}
for self-adjoint operators $S\in \CC^{n\times n}$, $U\in \CC^n$ (exercise).
Here and below, 
$\langle \cdot, \cdot \rangle$ denotes the standard, (complex) $L^2$
inner product with respect to variable $x$.
%
We shall require also the following elementary bounds.

\begin{lem}[Strong Sobolev embedding principle]\label{soblem}
For $s>d/2$, 
\begin{equation}\label{sobolev}
|f|_{L^\infty(\RR^d)}\le |\hat f|_{L^1(\RR^d)}\le C|f|_{H^s(\RR^d)},
\end{equation}
where $\hat f$ denotes Fourier transform of $f$.
\end{lem}

\begin{proof}
The first inequality follows by Hausdorff--Young's inequality,
the second by 
\begin{equation}
\begin{aligned}
|\hat f(\xi)|_{L^1}&=
|\hat f(\xi)(1+|\xi|^s)(1+|\xi|^s)^{-1}|_{L^1} \\
&\le
|(1+|\xi|^s)^{-1}|_{L^2} |\hat f(\xi)(1+|\xi|^s)|_{L^2} 
\le C|f|_{H^s}.
\end{aligned}
\end{equation}
\end{proof}

\begin{lem}[Weak Moser inequality]\label{moslem}
For $s=\sum |\alpha_j|$ and $k\ge d/2$,
\begin{equation}
\label{moser}
\begin{aligned}
|(\partial^{\alpha_1}v_1)\cdots (\partial^{\alpha_r}v_r)|_{L^2(\RR^d)}
&\leq \sum^r_{i=1}|v_i|_{H^s(\RR^d)}
\Big(\prod_{j\neq i}|\hat v_j|_{L^1(\RR^d)}\Big)\\
&\leq C\sum^r_{i=1}|v_i|_{H^s(\RR^d)}
\Big(\prod_{j\neq i}|v_j|_{H^k(\RR^d)}\Big).\\
\end{aligned}
\end{equation}
\end{lem}

\begin{proof} By repeated application of the Hausdorff--Young
inequality $|f*g|_{L^p}\le |f|_{L^2}|g|_{L^p}$, where $*$
denotes convolution, we obtain the first inequality,
\begin{equation}
\begin{aligned}
|(\partial^{\alpha_1}v_1)&\cdots (\partial^{\alpha_r}v_r)|_{L^2}
=
|(\xi^{\alpha_1}\hat v_1)*\cdots* (\xi^{\alpha_r}\hat v_r)|_{L^2}\\
&\leq 
\Big| |\xi^s \hat v_1|* \cdots* |\hat v_r|\Big|_{L^2}
+\cdots+
\Big| |\hat v_1|* \cdots* |\xi^s \hat v_r|\Big|_{L^2}\\
&\leq \sum^r_{i=1}|v_i|_{H^s}
\Big(\prod_{j\neq i}|\hat v_j|_{L^1}\Big).\\
\end{aligned}
\end{equation}
The second follows by \eqref{sobolev}.
\end{proof}

\begin{prop} [{[Fr, G]}]
\label{hypwellposed}
Symmetric form \eqref{symmhyp} 
implies local well-posedness in $H^s$,
and bounded local stability $|U(t)|_{H^s}\le C|U_0|_{H^s}$,
provided $\tilde A^j(\cdot)\in C^s$ and 
$s\ge [d/2]+2$. 
\end{prop}

\begin{proof}
Integration by parts together with \eqref{symmhyp} yields
the basic $L^2$ estimate
\begin{equation}
\label{0friedrichs}
\begin{aligned}
\frac{1}{2}\langle U, \tilde A^0 U\rangle_t 
&=
\frac{1}{2}\langle U, \tilde A^0_t U\rangle+
\langle U, \tilde A^0 U_t \rangle\\
&=
\frac{1}{2}\langle U, \tilde A^0_t U\rangle-
\langle U, \sum_j\tilde A^j U_{x_j} \rangle\\
&=
\frac{1}{2}\langle U, (\tilde A^0_t+\sum_j \tilde A^j_{x_j}) U \rangle\\
&\le
C|U|_{L^2}^2|U|_{W^{1,\infty}}
\le
C|U|_{L^2}^2|U|_{H^s}.\\
\end{aligned}
\end{equation}
Here, we have used \eqref{symmrel} in equating
$\langle U, \sum_j\tilde A^j U_{x_j} \rangle
=
-\frac{1}{2}\langle U, \sum_j \tilde A^j_{x_j} U \rangle$,
original equation $U_t= -\sum_j A^j U_{x_j}$ in estimating
$\tilde A^0_t \le C|U_t|\le C_2 |U_x|$ in the second-to-last inequality,
and \eqref{sobolev} in the final inequality.

A similar, higher-derivative calculation yields the $H^s$ estimate
\begin{equation}
\label{sfriedrichs}
\begin{aligned}
\frac{1}{2}\big(
\sum_{r=0}^s
\langle \partial_x^r, \tilde A^0 \partial_x^r U\rangle
\big)_t 
&=
\frac{1}{2}
\sum_{r=0}^s
\langle \partial_x^r U, \tilde A^0_t \partial_x^r U\rangle+
\sum_{r=0}^s
\langle \partial_x^r U, \tilde A^0 \partial_x^r U_t \rangle\\
&=
\sum_{r=0}^s
\frac{1}{2}\langle \partial_x^r U, \tilde A^0_t \partial_x^r U\rangle-
\sum_{r=0}^s
\langle \partial_x^r U, \sum_j\tilde A^0 \partial_x^r A^j U_{x_j} \rangle\\
&=
\frac{1}{2}\sum_{r=0}^s
\langle \partial_x^r U, (\tilde A^0_t+\sum_j \tilde A^j_{x_j}) 
\partial_x^r U \rangle \\
&\quad + 
\sum_{r=0}^s \sum_{\ell=1}^r 
\langle \partial_x^r U, \partial_x^\ell \tilde A^0 \partial_x^{r-\ell} A^j U_{x_j} \rangle \\
&\le
C|U|_{H^s}^2
\big(|U|_{H^s}+
\big(|U|_{H^s}^s\big),
\end{aligned}
\end{equation}
where the final inequality follows by \eqref{moser}; 
see Exercise \ref{diff} below.

Defining
\begin{equation}
\zeta(t):=
\frac{1}{2}\big(
\sum_{r=0}^\ell
\langle \partial_x^r, \tilde A^0 \partial_x^r U\rangle
\big),
\end{equation}
we have, therefore, the Ricatti-type inequality
\begin{equation}
\zeta_t\le C\big( \zeta+ \zeta^s\big),
\end{equation}
yielding $\zeta(t)\le C\zeta(0)$ for small $t$, provided
$\zeta(0)$ is sufficiently small.  Observing that
$\zeta^{1/2}$ is a norm equivalent to $|U|_{H^\ell}$, we 
obtain bounded local stability, provided a solution exists.

Essentially the same a priori estimate can be used to show
existence and uniqueness of solutions. Define the standard
nonlinear iteration scheme (see, e.g., [Fr, Ma.3])
\begin{equation}
\label{it}
\tilde A^0(U^n)U^{n+1}_t + \sum_j \tilde A^j(U^n)U^{n+1}_{x_j}=0,
\qquad U(0)=U_0.
\end{equation}
For $U^n\in H^s$, an $H^s$ solution $\CalT U^{n+1}$ of \eqref{it}
may be obtained by 
linear theory; see Remark \ref{linearexistence}, Section 3.2.
By the estimate already obtained, we find, 
that $\CalT$ takes the ball 
$\CalB:=\{U: \, |U|_{L^\infty([0,\tau]; H^s(x))}\le 2|U_0|_{H^s(x)}\}$ 
to itself,
for $\tau>0$ sufficiently small. 
A similar energy estimate on the variation $e:=\CalT(U_1)- \CalT(U_2)$,
for $U_j\in \CalB$ yields that $\CalT$ is contractive in 
$L^\infty([0,\tau];L^2(x))$ on the invariant
set $\CalB$, and stable in $L^\infty[0,\tau]; H^s(x)$,
yielding existence of a fixed-point solution $U\in L^\infty([0,T]; H^s(x))$
(Exercise \ref{gen}).
Likewise, uniqueness of solutions may be obtained by a stability estimate
on the nonlinear variation $e:= U_1-U_2$, where $U_1$ and $U_2$
denote solutions of \eqref{hyperbolic}.
\end{proof}

\begin{rems}
\label{blowup}
\textup{
1. Clearly, we do not obtain global well-posedness by this argument,
since solutions of Ricatti-type equations in general blow up in finite time.
Indeed, it is well-known that shock-type discontinuities
may form in finite time even for arbitrarily smooth initial data, 
corresponding to blow-up in $H^1$; see, e.g., [La.1--2, J, KlM, Si].
}
\smallbreak

\textup{ 2.  
Using the {\it strong Moser inequality} 
\begin{equation}
\label{smoser}
|(\partial^{\alpha_1}v_1)\cdots (\partial^{\alpha_r}v_r)|_{L^2(\RR^d)}
\leq C\sum^r_{i=1}|v_i|_{H^s(\RR^d)}
\Big(\prod_{j\neq i}|v_j|_{L^\infty(\RR^d)}\Big)
\end{equation}
for $s=\sum |\alpha_j|$
(proved using Gagliardo--Nirenberg inequalities [T]),
the same argument may be used to show that smooth continuation
of the solution is possible so long as $|U|_{W^{1,\infty}}$ remains
bounded; see [Ma.3], Chapter 2.
}
\end{rems}

\begin{ex}\label{diff}
Verify the final inequality in \eqref{sfriedrichs}
by showing that
\begin{equation}
|\partial_x^{r-s} A(U)U_{x}|
\le C\sum_{\sum |\alpha_j| \le r+1-s} \Pi_{1\le j\le r+1-s} |\partial_x^{\alpha_j}|.
\end{equation}
\end{ex}

\begin{ex}\label{gen}
If $u_n\in H^s(x)$ are uniformly bounded in $H^s$,
and convergent in $L^2(x)$, show that $\lim_{n\to \infty }u_n\in H^s$,
using the definition of $H^s$ as the set of $v\in L^2$
such that, for all $1\le r\le s$,
$\langle v, \partial_x^r \phi\rangle\le C|\phi|_{L^2}$
for all test functions $\phi\in C^\infty_0$,
together with the fact that limits and distributional derivatives commute.
%
\end{ex}

Symmetrizability is at first sight a rather restrictive requirement 
in more than one spatial dimension.  However, it turns out to be
satisfied in many physically interesting situations, in particular
for gas dynamics and MHD.
Indeed, a fundamental observation of Godunov
is that symmetrizability is closely related with existence of
an associated convex entropy.

\begin{defi}\label{hypentropy}
A hyperbolic entropy, entropy flux ensemble is a set of scalar functions
$(\eta, q^j)$ such that 
\begin{equation}
\label{etaq}
d\eta dF^j=dq^j,
\end{equation}
or equivalently
\begin{equation}
\label{entropyeqn}
\eta(U)_t + \sum_j q^j(U)_{x_j}=0
\end{equation}
for any smooth solution $U$ of \eqref{hyperbolic}.
\end{defi}

\begin{prop}
[{[God, Mo, B, KSh]}]
\label{godunov}
For $U$ lying in a convex set $\cal{U}$, 
existence of a convex entropy $\eta$ is equivalent to
symmetrizability of \eqref{hyperbolic}
by an invertible coordinate change $U\to W:=d\eta$
(known as an ``entropy variable''),
i.e., writing
\begin{equation}
\label{Weq}
U(W)_t + \sum_j F^j(W)_{x_j}=0,
\end{equation}
we have $\tilde A^0:= (\partial U/\partial W)=(d^2\eta)^{-1}$ symmetric
positive definite and 
$\tilde A^j:= (\partial F^j/\partial W)=A^j\tilde A^0$ symmetric.
\end{prop}

\begin{proof}
(Exercise)
($\Rightarrow$) Differentiate \eqref{etaq} and use symmetry of $d^2\eta$. 
($\Leftarrow$) Reverse the calculation to obtain \eqref{etaq}
with $d\eta:=W$, then note that $d\eta$ is exact, due to
symmetry of $dW:=\tilde A^0$. 
\end{proof}

\begin{rems}
\label{hentropyrems}
\textup{
1. Symmetrizability by coordinate change implies symmetrizability
in the usual quasilinear sense (exercise), but not the converse.
In particular, (nonlinear) symmetrization by
coordinate change preserves divergence form, whereas
(quasilinear) symmetrization by a left-multiplier
$\tilde A^0$ does not,  cf. \eqref{Weq} and \eqref{hypparab}.
\smallbreak
2. Symmetrizability by coordinate change implies also \eqref{entropyeqn},
which yields the additional information that $\int \eta dx$, without
loss of generality equivalent to the $L^2$ norm of $U$, is conserved
for smooth solutions.  Thus, we find in the discussion of Remark
\ref{blowup} that blowup occurs in a derivative of $U$
and not in $U$ itself.
\smallbreak
3. For gas dynamics and MHD, there exists a convex entropy in
the neighborhood of any thermodynamically stable state,
namely the negative of the thermodynamical entropy $s$; see, e.g.
[Kaw, MaZ.4, Z.4].
In particular,
for an ideal gas, there exists a global convex entropy.
\smallbreak
4. For the hyperbolic shock stability problem, 
hypotheses of hyperbolicity, symmetrizability, etc.,
are relevant only in neighborhoods of $U_\pm$ and
not between.
Thus, a shock may be stable even if $U_\pm$ are entirely
separated by unstable constant states, as, e.g., for
phase-transitional shocks in van der Waals gas dynamics; see [Fre.3, B--G.2--3].
}
\end{rems}

\medbreak
{\bf 1.1.2. Viscous equations.}
Analogous to \eqref{symmhyp} in the setting of the viscous equations
\eqref{viscous} is the {\it symmetric hyperbolic--parabolic form}
\begin{equation}
\label{symmhypparab}
\tilde A^0 W_t + \sum_j \tilde A^j W_{x_j}= \sum_{j,k} 
(\tilde B^{jk} W_{x_k})_{x_j} +
\begin{pmatrix} 0\\ \tilde g \end{pmatrix},
\end{equation}
$\tilde A^0$ symmetric positive definite, $\tilde A^j_{11}$ symmetric, 
$\tilde B^{jk}=\bdiag \{0, \tilde b^{jk}\}$ with 
\begin{equation}
\label{goodbtilde}
\sum \xi_j \xi_k \tilde b^{jk} \ge \theta |\xi|^2,
\quad \theta>0,
\end{equation}
for all $\xi\in \RR^d$, and
\begin{equation}
\label{Gtilde}
\tilde G=
\begin{pmatrix} 0\\ \tilde g(\partial_x W) \end{pmatrix}
\end{equation}
with $\tilde g={\cal{O}}(|W_x|^2)$, 
to be achieved by an invertible coordinate change $U\to W$
combined with left-multiplication by an invertible lower block-triangular
matrix $S(U)$;
see [Kaw, KSh] and references therein, or Appendix A1, [Z.4].

\begin{rem}
\label{varchoice}
\textup{
As pointed out in the references, $A^0$ may be taken without loss of 
generality to be block-diagonal, thus identifying ``hyperbolic''
and ``parabolic'' variables $w^I$ and $w^{II}$.
Indeed, $w^I$ may be taken without loss of generality as $u^I$
and $w^{II}$ as any variable satisfying the (clearly necessary) 
integrability condition 
$B^{jk}U_{x_k}= \beta^{jk}(U)w^{II}_{x_k}$ for all $j$, $k$ [GMWZ.4, Z.4].
}
\end{rem}

Similarly as in the hyperbolic case, we may deduce local well-posedness
of \eqref{symmhypparab} directly from the structure of the equations,
using energy estimates/integration by parts.
Here, we shall require also 
a standard but essential tool for multidimensional parabolic systems,
the {\it G\"arding inequality}
\begin{equation}
\label{garding}
\sum_{j,k}\langle \partial_{x_j}f, \tilde b^{jk} \partial_{x_k}f\rangle
\ge \tilde \theta |\partial_x f|_{L^2}^2 - C|f|_{L^2}^2
\end{equation}
for Lipshitz $\tilde b^{jk}\in \RR^{n\times n}$ satisfying
uniform ellipticity condition \eqref{goodbtilde}, and $0<\tilde \theta<\theta$,
with $C=C(\tilde \theta, |\partial_x \tilde b^{jk}|_{L^\infty})
\le C_2 |\partial_x \tilde b^{jk}|_{L^\infty}/|\theta-\tilde \theta|$
for some uniform $C_2>0$.

\begin{ex}
\label{gardingex}
(i) Prove \eqref{garding} with $C=0$ in the case $\tilde b^{jk}\equiv
\const$, using the Fourier transform and Parseval's identity.
(ii) Prove \eqref{garding} with $C=0$ in the case that 
$\tilde b^{jk}$ varies by less that $|\theta-\tilde \theta|/d^2$
from some constant value, i.e, 
$\osc \tilde b^{jk}\le 2|\theta-\tilde \theta|/d^2$.
(iii) Prove the general case using a partition of unity 
$\{\chi_r\}$ such that $\osc \tilde b^{jk}\le 2|\theta-\tilde \theta|/d^2$
on the support of each $\chi_r$, and the estimates
\begin{equation}
\sum_{j,k}\langle \partial_{x_j}f, \tilde b^{jk} \partial_{x_k}f\rangle
=\sum_r\sum_{j,k}\langle \partial_{x_j}\chi_r f, 
\tilde b^{jk} \partial_{x}\chi_r f\rangle
+{\cal{O}}(|\partial_x f|_{L^2}|(\partial_x \chi_r)f|_{L^2})
\end{equation}
and
\begin{equation}
|\partial_x f|_{L^2}=
\sum_r |\chi_r \partial_x f|_{L^2}
=
\sum_r \Big(
 |\partial_x \chi_r f|_{L^2}+ |(\partial_x \chi_r)f|_{L^2})
\Big).
\end{equation}
\end{ex}

\begin{rem}
\label{pseudodiff}
\textup{
The G\"arding inequality \eqref{garding} 
is an elementary example of a pseudodifferential estimate.
Pseudodifferential techniques 
play a fundamental role 
in the analysis of the curved shock problem;
see [Ma.1--3, M\'e.4],
[GMWZ.1, GMWZ.3--4, W],
and references therein.
}
\end{rem}

\begin{prop} [{[Kaw]}]
\label{parabwellposed}
Symmetric form \eqref{symmhypparab} 
implies local well-posedness in $H^s$,
and bounded local stability $|U(t)|_{H^s}\le C|U_0|_{H^s}$,
provided $\tilde A^j(\cdot)$, $\tilde B^{jk}(\cdot)\in C^s$
and 
$s\ge [d/2]+2$.
\end{prop}

\begin{proof}
For simplicity, take $\tilde g\equiv 0$; the general case is similar.
Similarly as in \ref{0friedrichs}, we have
\begin{equation}
\label{v0friedrichs}
\begin{aligned}
\frac{1}{2}\langle W, \tilde A^0 W\rangle_t 
&=
\frac{1}{2}\langle W, \tilde A^0_t W\rangle+
 \langle W, \tilde A^0 W_t \rangle\\
&=
\frac{1}{2}\langle W, \tilde A^0_t W\rangle-
 \langle W, \sum_j\tilde A^j W_{x_j} - 
\sum_{j,k}(\tilde B^{jk}W_{x_k})_{x_j} \rangle\\
&=
\frac{1}{2}\langle W, \tilde A^0_t W\rangle+ 
\langle w^I, \sum_j \tilde A^j_{11,x_j} w^I\rangle 
+ {\cal{O}}(|\partial_x w^{II}||W|) \\
&\quad
- \sum_{j,k} \langle w^{II}_{x_j}, \tilde B^{jk}w^{II}_{x_k}\rangle\\
&\le
C|W|_{L^2}^2|W|_{W^{1,\infty}}\\
&\le 
C|W|_{L^2}^2|W|_{H^{[d/2]+2}},\\
\end{aligned}
\end{equation}
where in the second-to-last inequality we have used \eqref{garding}
with $\tilde \theta=\theta/2$, together with Young's inequality
$|\partial_x w^{II}||W|\le (1/2C) |\partial_x w^{II}|^2 + (C/2)|W|^2$
with $C>0$ sufficiently large, and used the original equation to bound
$|W_t|\le C(|W_x|+|\partial_x^2 w^{II}|)$, and in the final inequality
we have used the Sobolev inequality \eqref{sobolev}. 

Likewise, we obtain by a similar calculation
\begin{equation}
\label{vsfriedrichs}
\begin{aligned}
\frac{1}{2}\big(
\sum_{r=0}^s
\langle \partial_x^r, \tilde A^0 \partial_x^r W\rangle
\big)_t 
&\le
C|W|_{H^s}^2
\big(|W|_{H^{[d/2]+2}}+
|W|_{H^{[d/2]+2}}^s\big)\\
&\le
C\big(|W|_{H^s}^3+
|W|_{H^s}^{s+1}\big),\\
\end{aligned}
\end{equation}
by the Moser inequality \eqref{moser} and $s\ge [d/2]+2$,
from which the result follows by the
same argument as in the proof of Proposition \ref{hypwellposed}.
\end{proof}

The following beautiful results of Kawashima et al
(see [Kaw, KSh] and references therein) 
generalize to the viscous case the observations of Godunov et al
regarding entropy and the structure of the inviscid equations.

\begin{defi}[{[Kaw, Sm]}]
\label{ventropy}
A viscosity-compatible convex entropy, entropy flux ensemble 
$(\eta, \, q^j)$ for \eqref{quasiviscous} is a convex hyperbolic
entropy, entropy flux ensemble such that 
\begin{equation}
\label{Bcompatible}
d^2\eta \sum_{j,k} \tilde B^{jk} \xi_j \xi_k \ge 0
\end{equation}
and 
\begin{equation}
\label{Brank}
\rank \, \R d^2\eta \sum_{j,k} \tilde B^{jk} \xi_j \xi_k 
= \rank \sum_{j,k} \tilde B^{jk} \xi_j \xi_k \equiv r
\end{equation}
for all nonzero $\xi\in \RR^d$.
\end{defi}

\begin{ex}
\label{suff}
Using \eqref{hypparab}--\eqref{goodb}, and the fact that
\begin{equation}
(d^2\eta)^{1/2}B^{jk}(d^2\eta)^{-1/2} 
=
(d^2\eta)^{-1/2}(d^2\eta B^{jk})(d^2\eta)^{-1/2} 
\end{equation}
is similar to $B^{jk}$, show that $d^2\eta B^{jk}$ symmetric
is sufficient for \eqref{Bcompatible}--\eqref{Brank}.
\end{ex}

\begin{prop}[{[KSh, Yo]\footnote{
Established for $d^2\eta B^{jk}$ symmetric in [KSh].
The observation that \eqref{Bcompatible}--\eqref{Brank} suffice
seems to be due to Yong [Yo].
}
}]
\label{entropyprop}
A necessary and sufficient condition that a system
\eqref{viscous}, \eqref{goodb} can be put in symmetric
hyperbolic--parabolic form \eqref{symmhypparab} 
with $\tilde G\equiv 0$ and $\tilde A^j$ symmetric
($\tilde A^0$ not necessarily block-diagonal)
by a change of coordinates $U\to W(U)$
for $U$ in a convex set $\cal{U}$, 
is existence of a convex viscosity-compatible 
entropy, entropy flux ensemble $\eta$, $q^j$ defined on $\cal{U}$,
with $W=d\eta(U)$.
\end{prop}

\begin{proof}
Identical with that of Proposition \ref{godunov} as concerns $\tilde A^0$
and $\tilde A^j$.
Regarding
$\tilde B^{jk} =B^{jk}(d^2\eta)^{-1}=
(d^2\eta)^{-1}(d^2\eta B^{jk})(d^2\eta)^{-1}$, we have that conditions
$\sum_{j,k}\tilde B^{jk}\xi_j\xi_k\ge 0$ and 
$\rank \, \R \tilde B^{jk}\xi_j\xi_k \equiv r$ 
for all nonzero $\xi\in \RR^d$ 
are equivalent both to \eqref{Bcompatible}--\eqref{Brank}
and to \eqref{goodbtilde} and
$\tilde B^{jk}=\bdiag\{0,\tilde b^{jk}\}$ (exercise). 
The latter assertion depends on the observation that
$\sum_{j,k}\tilde B^{jk}\xi_j\xi_k\ge 0$ 
together with structure \eqref{hypparab} implies that
$\R\sum_{j,k}\tilde B^{jk}\xi_j \xi_k$ is block-diagonal and vanishing
in the first diagonal block.
This completes the argument.
\end{proof}

\begin{defi}[{[Kaw]}]
\label{gc}
Augmenting the local conditions \eqref{symmhypparab}--\eqref{goodbtilde},
we identify the time-asymptotic stability conditions\footnote{
This refers to stability of constant solutions;
see Proposition \ref{conststab} below.}
of ``first-order symmetry,''
$\tilde A^j$ symmetric,
and ``genuine coupling:''
\begin{equation}
\tag{GC}
\text{\rm No eigenvector of $\sum_j \xi_j dF^j$ lies in 
$\ker \sum \xi_j \xi_k B^{jk}$, for $\xi\ne 0\in \RR^d$.}
\end{equation}
\end{defi}

\begin{prop}
[{[Kaw]\footnote{Established in [Kaw] for $d\ge 3$, or $d\ge 1$ in the case
of a convex entropy.}}]
\label{conststab}
Symmetric hyperbolic--parabolic form \eqref{symmhypparab}
together with the time-asymptotic stability conditions
of first-order symmetry, 
$\tilde A^j$ symmetric, and genuine coupling, (GC),
implies $L^1\cap H^s \to H^s$ time-asymptotic stability of 
constant solutions, $\bar U\equiv \const$, provided that
$\tilde A^j$, $\tilde B^{jk}\in C^s$ and 
$s\ge [d/2]+2$,
with rate of decay
\begin{equation}
\label{constrate}
|U-\bar U|_{H^s}(t)\le C(1+t)^{-\frac{d}{4}} |U-\bar U|_{H^s\cap L^1}(0)
\end{equation}
for $|U-\bar U|_{H^s\cap L^1}(0)$ sufficiently small,
equal to that of a $d$-dimensional heat kernel.
\end{prop}

\begin{prop}
[{[Kaw]\footnote{Established in [Kaw] for
``physical'' convex entropies satisfying \eqref{ellipt} below.}}]
\label{globalwellposed}
Existence of a convex viscosity-compatible entropy, 
together with genuine coupling, (GC),
implies global well-posedness of \eqref{quasiviscous}--\eqref{goodb}
in $H^s$, and bounded stability of constant solutions,
\begin{equation}
\label{gwellposed}
|U-\bar U|_{H^s}(t)\le C |U-\bar U|_{H^s}(0),
\end{equation}
provided $\tilde A^j$, $\tilde B^{jk}\in C^s$ and 
$s\ge [d/2]+2$.
\end{prop}

Propositions \ref{conststab} and \ref{globalwellposed} depend
on a circle of ideas associated with the phenomenon 
of ``hyperbolic--parabolic smoothing'', as indicated by the existence of
parabolic-type energy estimates
\begin{equation}
\label{Kenergy}
|W(t)|^2_{H^s}+ \int_0^t \Big(|\partial_x W|^2_{H^{s-1}}+ 
|\partial_x w^{II}|^2_{H^s}\Big)ds\le
C\Big(|W(0)|^2_{H^s} + \int_0^t |W(s)|^2_{L^2}ds\Big)
\end{equation}
for $s$ sufficiently large.
We defer discussion of these important concepts to the detailed treatment
of Sections 3 and 4.

\begin{ex}
\label{gposed}
Result \eqref{gwellposed}
was established for gas dynamics
in the seminal work of Matsumura and Nishida [MNi]. 
The key point is existence of an $L^2$ energy estimate
\begin{equation}
\label{L2}
|W(t)|_{L^2}^2+ \int_0^t |\partial_x w^{II}|^2_{L^2} \le C|W(0)|_{L^2}^2,
\end{equation}
which, coupled with the general machinery used to obtain \eqref{Kenergy},
yields an improved version of \eqref{Kenergy} 
in which the $\int|U(s)|^2_{L^2} ds$ term
on the righthand side does not appear.  This in turn implies
\eqref{gwellposed}.
Prove \eqref{L2} in the general case using the viscous version
\begin{equation}
\label{ventropyeqn}
\eta(U)_t + \sum_j q^j(U)_{x_j}=\sum_{j,k}(d\eta B^{jk}U_{x_k})_{x_j}
-\sum_{j,k} U_{x_j}^t d^2\eta B^{jk} U_{x_k}
\end{equation}
of entropy equation \eqref{entropyeqn}
and Exercise \ref{gardingex}(ii), under the assumption that $|W|_{H^s}$
(hence $|W|_{L^\infty}$) remains sufficiently small, for $s$ sufficiently
large.
\smallbreak
The G\"arding inequality is not required for gas dynamics or MHD, for
which there hold the strengthened ellipticity condition
\begin{equation}
\label{ellipt}
\sum_{j,k}\langle f_{x_j}, d^2\eta B^{jk} f_{x_k}\rangle \ge 
\theta |\partial x f|^2.
\end{equation}
Under this assumption, show using \eqref{ventropyeqn}
that \eqref{gwellposed} holds for $s=0$,
independent of the size of $|W(0)|_{L^2}$.
\end{ex}


\begin{rem}
\label{vrem}
\textup{
Similarly as in the inviscid case, 
the negative of the thermodynamical entropy 
serves as a convex viscosity-compatible entropy
for gas and MHD in the neighborhood of any 
thermodynamically stable state; see [Kaw, MaZ.4, Z.4].  
In particular, for an ideal gas, there exists a global convex 
viscosity-compatible entropy.
}
\end{rem}

\begin{exam}
\label{ns}
\textup{
The Navier--Stokes equations of compressible gas dynamics,
may be written as
\begin{equation}
\label{nseq}
\begin{aligned}
\rho_t + \div (\rho u)=&0,\\
(\rho u)_t + \div(\rho u\otimes u)+ \nabla p&= 
\overbrace{
\mu \Delta u + (\lambda + \mu) \nabla \div u}^{\div \tau},\\
(\rho(e+ \frac{1}{2} u^2))_t + \div(\rho(e+\frac{1}{2}u^2) u + pu)&=
\div(\tau\cdot u)+ \kappa \Delta T,\\
\end{aligned}
\end{equation}
where $\rho>0$ denotes density, $u\in \RR^d$ fluid velocity, $T>0$
temperature, $e=e(\rho,T)$ internal
energy, and $p=p(\rho, T)$ pressure. 
Here, $\tau:= \lambda \div(u)I + 2\mu Du$, 
where $Du_{jk}=\frac{1}{2}(u^j_{x_j}+
u^j_{x_k})$ is the deformation tensor, $\lambda(\rho, T)>0$ and 
$\mu(\rho,T)>0$
are viscosity coeffients, and $\kappa(\rho,T)>0$ is the coefficient
of thermal conductivity.
The thermodynamic entropy $s$ is defined implicitly by
the underlying thermodynamic relations 
$e=\hat e(v,s)$, $T=\hat e_s$, $p=-\hat e_v$, 
where $v=\rho^{-1}$ is specific volume.
(These may or may not be integrable for arbitrary choices of
$p(\cdot, \cdot)$, $e(\cdot, \cdot)$; on the other hand, each choice
of $\hat e$ satisfying $\hat e_s=T>0$ uniquely determines functions $p$ 
and $e$.)
Thermodynamic stability is
the condition that $\hat e$ be a convex function of $(v,s)$,
or equivalently (exercise) $e_T$, $p_\rho>0$.
The Euler equations of compressible gas dynamics
are \eqref{nseq} with the righthand side set to zero.}

\textup{
Evidently, \eqref{nseq}
is of form \eqref{viscous}, \eqref{hypparab}, with ``conservative variables''
$U=\big(\rho, \rho u, \rho(e+ \frac{1}{2}u^2)\big)$.
By Remark \ref{varchoice}, we have that symmetric hyperbolic--parabolic
form \eqref{symmhypparab}, 
if it exists, may be expressed in the ``natural variables''
$W=(\rho, u, T)$.
Indeed, this can be done [KSh], with 
\begin{equation}
\label{1}
\tilde A^0=
\begin{pmatrix}
p_\rho/\rho & 0 & 0\\
0 & \rho I_d & 0\\
0 & 0 & \rho e_T/T\\
\end{pmatrix},
\end{equation}
\begin{equation}
\label{2}
\sum_j \tilde A^j \xi_j=
\begin{pmatrix}
(p_\rho)u\cdot \xi& p_\rho \xi & 0\\
p_\rho \xi^t & \rho (u\cdot \xi)I_d &  p_T\xi^t\\
0 & p_T\xi & (\rho e_T/T)u\cdot \xi\\
\end{pmatrix},
\end{equation}
and
\begin{equation}
\label{3}
\sum_{j,k}\tilde B^{jk}\xi_j\xi_k=
\begin{pmatrix}
0&0&0\\
0 & \mu |\xi|^2 I_d + (\mu+\lambda)\xi^t\xi & 0 \\
0 & 0 & T^{-1}|\xi|^2 \\
\end{pmatrix},
\end{equation}
whenever $p_\rho$, $e_T>0$ (thermodynamic stability), in 
which case we also have the time-asymptotic stability conditions
of symmetry of $\tilde A^j$ and genuine coupling (GC) (by inspection,
equivalent to $p_\rho\ne 0$).
If only $e_T>0$, as for example for a van der Waals-type equation
of state,
we may still achieve symmetric hyperbolic--parabolic
form, but without symmetry of $\tilde A^j$ or genuine coupling,
by dividing through $p_\rho$ from the first equation ($\sim$ first
row of each matrix $\tilde A^0$, $\tilde A^j$, $\tilde B^{jk}$).
That is, we obtain in this case local well-posedness, but not asymptotic
stability of constant states.
Similar considerations hold in the case of MHD; see [Kaw, KSh].
}
\end{exam}

%

\bigbreak
\clearpage
\section{Description of results}
We now turn to the long-time stability
of viscous shock waves.  Fixing the viscosity coefficient 
$\nu$, consider a viscous shock solution
\eqref{profile} of a system of viscous conservation laws
\eqref{viscous} of form \eqref{hypparab}--\eqref{goodb}.
Changing coordinates if necessary to a rest frame moving with
the speed of the shock, we may without loss of generality arrange
that shock speed $s$ vanish, so that \eqref{profile} becomes a stationary,
or standing-wave solution.
Hereafter, we take $\nu\equiv 1$ and $s\equiv 0$, and suppress
the parameters $\nu$ and $s$.

\medbreak
{\bf 2.1. Assumptions.}
The classical Propositions \ref{conststab} and \ref{globalwellposed}
concern long-time stability of a single, thermodynamically stable
equilibrium, whereas 
a viscous shock solution typically consists of  
two different thermodynamically stable
equilibria connected by a profile that in
general may pass through regions of thermodynamical instability;
see Remark \ref{hentropyrems}.4, and examples, [MaZ.4, Z.4].
Accordingly, we make the following structural assumptions,
imposing 
stability at the endstates, but only local well-posedness along the profile. 

\medbreak
\begin{assums} \label{A}{ $ \,$ }

\textup{(A1) }
For $U$ in a neighborhood of profile $\bar U(\cdot)$,
there is an invertible change of coordinates $U\to W$ 
such that \eqref{viscous} may be expressed in $W$ coordinates
in symmetric hyperbolic--parabolic form, i.e, in form
\eqref{symmhypparab} with
$\tilde A^0$ symmetric positive definite 
and (without loss of generality) block-diagonal 
and $\tilde A^j_{11}$ symmetric,
$\tilde B^{jk}=\bdiag\{0, \tilde b^{jk}\}$ with $\tilde b^{jk}$
satisfying uniform ellipticity condition \eqref{goodbtilde}, and
$\tilde G=\begin{pmatrix} 0\\\tilde g\\\end{pmatrix}$ with
$\tilde g={\cal{O}}(|\partial x W|^2)$.
\medbreak

\textup{(A2)}
At endstates $U_\pm$, the coefficients $\tilde A^j$,
$\tilde B^{jk}$ defined in (A1)
satisfy the constant-coefficient stability
conditions of first-order symmetry,
$\tilde A^j_\pm$ symmetric, and genuine coupling, (GC).
\end{assums}
\medbreak

\begin{rem}\label{satisfaction}
\textup{
Conditions (A1)--(A2) are satisfied for gas- and magnetohydrodynamical
shock profiles connecting thermodynamically stable endstates,
under the mild assumption $e_T>0$ on the equation of state
$e=e(\rho, T)$ relating internal energy $e$ to density $\rho$
and temperature $T$:
in particular for both ideal and van der Waals-type equation of state.
(Exercise: verify this assertion for gas dynamics, starting from
form \eqref{1}--\eqref{3} and modifying all coeffients $M$ by
the transformation 
\begin{equation}
M \to \chi M+ (1-\chi)\bdiag \{p_\rho^{-1},I_d,1\}M,
\end{equation}
where $\chi=\chi(W)$ is a smooth cutoff function supported on the region of
thermodynamic stability $p_\rho >0$.)
}
\end{rem}

To Assumptions \ref{A}, we add the following technical hypotheses.

\begin{assums} \label{Hj}{ $ \,$ }

\textup{(H0) }
$F^j$, $B^{jk}$, $W(\cdot)$, $\tilde A^0 \in C^{q+1}$,
$q\ge q(d):= [d/2]+2$.
\medbreak

\textup{(H1)}
$\tilde A^1_{11}(\bar U)$ is uniformly definite, without loss of generality
$\tilde A^1_{11}\ge \theta >0$.
(For necessary conditions, $\det \tilde A^1_{11}\ne 0$ is sufficient.)

\medbreak

 \textup{(H2)}   
$\det \tilde A^1_\pm \ne 0$, or equivalently $\det A^1_\pm \ne 0$.

\medbreak

\textup{(H3)}  
Local to $\bar U$, the solutions of \eqref{profile} form a smooth
manifold $\{\bar U^\delta\}$, $\delta \in {\cal{U}}\subset \RR^\ell$,
with $\bar U^0=\bar U$.
\end{assums}

\medbreak
Condition (H0) gives the regularity needed for our analysis.
Condition (H2) is the standard inviscid requirement 
that the hyperbolic shock triple $(U_-,U_+,0)$ be noncharacteristic.
Condition (H1) states that convection in
the reduced, hyperbolic part of \eqref{symmhypparab} governing coordinate
$w^I$ is either up- or downwind, in particular uniformly noncharacteristic,
everywhere along the profile; 
as discussed in [MaZ.4, Z.4], this is satisfied for all gas-dynamical
shocks and at least generically for magnetohydrodynamical shocks.
Condition (H3) is a weak form of (implied by but not implying)
transversality of the connection $\bar U$ as a solution of the
associated traveling-wave ODE.

These hypotheses suffice for the investigation of necessary conditions for
stability.  In our investigation of sufficient conditions, we
shall require two further hypotheses at the level of the
inviscid stability problem.

\begin{ass} \label{H4}{ $ \,$ }
\medbreak
\textup{(H4) }
The eigenvalues of
$\sum_j A^j_\pm \xi_j$ have constant multiplicity
with respect to $\xi\in \RR^d$, $\xi\ne 0$.
\end{ass}

Denote by $a_j^\pm(\xi)$, $j=1, \dots, n$ the eigenvalues
of $\sum_j A^j_\pm \xi_j$, necessarily real by (A2), 
indexed by increasing order.
These are positive homogeneous degree one and, by (H4), 
locally analytic on $\xi\in \RR^d\setminus\{0\}$.\footnote{
In fact, they are globally analytic, since they maintain fixed order;
however, this is unimportant in the analysis.}
Here, and elsewhere, real homogeneity refers to homogeneity with
respect to positive reals.
Let
\begin{equation}
\label{symbol}
P_\pm(\xi,\tau):= i\tau + \sum_j i \xi_j A^j_\pm
\end{equation}
denote the frozen-coefficient symbols associated with \eqref{quasihyperbolic}
at $U=U_\pm$.
Then, $\det P_\pm(\xi, \tau)=0$ has $n$ locally analytic, positive homogeneous
degree one roots
\begin{equation}
\label{tau}
i \tau=-i a_r^\pm(\xi), \qquad r=1,\dots,n,
\end{equation}
describing the dispersion relations for the frozen-coefficient
initial-value problem (IVP).
The corresponding objects for the initial--boundary-value problem (IBVP)
are relations
\begin{equation}
\label{mu}
i\xi_1=\mu_r(\tilde \xi, \tau), \qquad r=1,\dots,n,
\end{equation}
$\xi=(\xi_1, \dots, \xi_d)=:(\xi_1, \tilde \xi)$,
describing roots of $\det Q_\pm(\tau)= 0$, where 
\begin{equation}
\label{IBVsymbol}
Q_\pm(\xi,\tau):= i\xi_1 +(A^1_\pm)^{-1}(i\tau + \sum_{j=2}^d i \xi_j A^j_\pm)
\end{equation}
denote the frozen-coefficient initial--boundary-value
symbols associated with \eqref{quasihyperbolic} at $U=U_\pm$.
Evidently, graphs \eqref{tau} and \eqref{mu} 
describe the same sets, since
$\det A^1 \det Q_\pm(\tau)= \det P_\pm $ and $\det A^1_\pm \ne 0$;
the roots $i\tau$ describe characteristic rates of temporal
decay, whereas $\mu=i\xi_1$ describe characteristic rates
of spatial decay in the $x_1$ direction.

\begin{defi}\label{glancing}
Setting $\xi=(\xi_1,\dots,\xi_d)=: (\xi_1, \tilde \xi)$,
we define the glancing sets ${\cal{G}}(P_\pm)$ as the
set of all $(\tilde \xi, \tau)$ such that, 
for some real $\xi_1$ and $1\le r\le n$,
$\tau=-a_r^\pm(\xi_1,\tilde \xi)$
and $(\partial a_r^\pm/\partial \xi_1)(\xi_1, \tilde \xi)=0$: 
that is, the projection onto $(\tilde \xi, \tau)$ of the
set of real roots $(\xi, \tau)$ of $\det P_\pm=0$ 
at which \eqref{tau} is not analytically invertible as a function \eqref{mu}.
The roots $(\xi, \tau)$ are called glancing points.
\end{defi}

\begin{ass}{$ $}
\medbreak
\textup{(H5) }
Each glancing set
${\cal{G}}(P_\pm)$ is the (possibly intersecting) 
union of finitely many smooth curves $\tau=\eta_q^\pm(\tilde \xi)$,
on which the root $\xi_1$ of $i\tau+ a_r^\pm(\cdot, \tilde \xi)=0$ 
has constant multiplicity $s_q$ (by definition $\ge 2$).
\end{ass}

Among other useful properties (see, e.g.,  [M\'e.3, Z.3, M\'eZ.2]), 
condition (H4) implies the
standard block structure condition [Ma.1--3, M\'e.2] 
of the inviscid stability theory [M\'e.3].
Condition (H5), introduced in [Z.3], imposes a further,
laminar structure on the glancing set that is convenient
for the viscous stability analysis.
In either context, glancing is the fundamental obstacle in 
obtaining resolvent estimates; see, e.g., [K, Ma.1--3, Z.3, M\'eZ.1, GMWZ.1--4].
The term ``glancing'' derives from the fact that, at glancing
points $(\xi,\tau)$, null bicharacteristics of $\det P_\pm$ 
lie parallel to the shock front $x_1\equiv 0$.

\begin{rems}\label{glancerems}
\textup{
1. Condition (H5) holds automatically in dimensions one 
(vacuous) and two (trivial),
and also for the case that all characteristics $a_r(\xi)$ 
are either linear or convex/concave in $\xi_1$. (Exercise [GMWZ.2]: 
prove the latter assertion using the Implicit Function Theorem.)
In particular, (H5) holds for both gas dynamics and MHD.
}

\textup{
2. Condition (H4) holds always for gas dynamics in any dimension, 
and generically for MHD in one dimension,
but fails always for MHD in dimension greater than or equal to two.
For a refined treatment applying also to multidimensional
MHD, see [M\'eZ.3].
}
\end{rems}

\medbreak
\begin{ass}
\label{lax}
For definiteness, consider a classical, ``pure'' Lax $p$-shock, satisfying
\begin{equation}
\tag{L}
a_p^->s=0>a_p^+,
\end{equation}
where $a_j^\pm:=a_j^\pm(1,0,\dots,0)$ denote the eigenvalues of 
$A^1_\pm:= dF^1(U_\pm)$, and $\ell=1$ in (H3).
See, e.g.,  [ZS, Z.3--4] for discussions of the general case, 
including the interesting situation of nonclassical 
over- or undercompressive shocks.
\end{ass}
\medbreak
{\bf 2.2. Classical stability conditions.}
Recall the classical admissibility conditions described in
the introduction of:

1. {\it Structural stability}, defined 
as existence of a transverse viscous profile.  
For gas dynamics, this is equivalent to the 
(hyperbolic) Liu--Oleinik admissibility condition; see [L.2, Gi, MeP].
In general, it may be a complicated ODE problem involving 
both \eqref{hyperbolic} and the specific form of the viscosity $B^{jk}$.

2. {\it Dynamical stability}, defined as local hyperbolic
well-posedness, or bounded, bounded-time stability of
ideal shock \eqref{shock}.
%
%
Following Majda [Ma.1--3], we define {\it weak dynamical stability}
(Majda's Lopatinski condition) as the absence of unstable spectrum 
$\R \lambda >0$ for the linearized operator 
(appropriately defined) 
about the shock, an evident necessary condition for stability.
(Recall, equations \eqref{hyperbolic} are positive homogeneous degree one,
hence spectrum lies on rays through the origin and 
instabilities, should they occur, occur to all orders.)

In the present context of a Lax $p$-shock, under assumptions
(A1)--(A2) and (H2), the necessary condition 
of weak dynamical stability
may be expressed in terms of the {\it Lopatinski determinant}
\begin{equation}
\label{Delta}
\Delta(\tilde \xi, \lambda):=
\det \Big(r_1^-, \cdots, r_p^-, r_{p+1}^+, 
\cdots, r_n^+, \lambda[U]+ i[F^{\tilde \xi}]\Big)
\end{equation}
as
\begin{equation}
\label{weaklop}
\Delta(\tilde \xi, \lambda)\ne 0 \quad \text{\rm for $\tilde \xi\in \RR^{d-1}$
and $\R \lambda >0$,}
\end{equation}
where $r_j^\pm=r_j^\pm(\tilde \xi, \lambda)$ denote 
bases for the unstable (resp. stable) subspace of
the coefficient matrix 
\begin{equation}
\label{CalAhyp}
\CalA_\pm:=(A^1)^{-1}(\lambda + iA^{\tilde \xi})_\pm
\end{equation}
arising in the IBVP symbol $Q_\pm$ with $\lambda:=i\tau$,
$F^{\tilde \xi}:=\sum_{j\ne 1}F^j \xi_j$, and
$A^{\tilde \xi}:= \sum_{j\ne 1}A^j \xi_j$.
Functions $r_j^\pm$, and thus $\Delta$, are well-defined 
on $\R \lambda>0$ and may be chosen analytically in $(\Txi, \lambda)$,
by a standard lemma of Hersch [H] asserting that
$\CalA_\pm$ has no center subspace on this domain;
see exercise \ref{hersch} below.
Zeroes of $\Delta$ correspond to ``normal modes'', or solutions
$(U,X)=e^{\lambda t}e^{i\tilde \xi \cdot \tilde x}
\big({\hat {\hat U}}(x_1), {\hat{\hat X}}\big)$
of the linearized perturbation equations, where 
$x_1=X(\tilde x,t)$ denotes location of the shock surface.
We define {\it strong dynamical instability} as failure
of \eqref{weaklop}. 

Under assumptions (A1)--(A2), (H2), and (H4), $r_j^\pm$, and
thus $\Delta$, may be extended continuously to the boundary
$\R \lambda =0$; see [M\'e.3, CP].
{\it Strong dynamical stability} (Majda's uniform Lopatinski condition)
is then defined as 
\begin{equation}
\label{stronglop}
\Delta(\tilde \xi, \lambda)\ne 0 \quad \text{\rm for $\tilde \xi\in \RR^{d-1}$,
$\R \lambda \ge 0$, and $(\tilde \xi, \lambda)\ne (0,0)$.}
\end{equation}
Strong dynamical stability under our assumptions
is sufficient for dynamical stability,
by the celebrated result of Majda 
[Ma.1--3] together with the result of M\'etivier [M\'e.3]
that (H4) implies Majda's block structure condition.

\begin{rems}
\label{majda}
\textup{
1. In the one-dimensional case, $\tilde \xi=0$, both
\eqref{weaklop} and \eqref{stronglop} reduce to nonvanishing
of the Liu--Majda determinant
\begin{equation}
\label{delta}
\delta:= \det \Big(r_1^-, \cdots, r_p^-, r_{p+1}^+, 
\cdots, r_n^+, [U] \Big),
\end{equation}
where $r_j^\pm$ denote eigenvectors of $A^1_\pm$ associated
with characteristic modes outgoing from the shock.
As a single condition, this is generically satisfied, whence
we find that shocks of Lax type (L) are generically
stable in one dimension from the hyperbolic point of view.
}

\textup{
2. As pointed out by Majda [Ma.1--3], for gas dynamics
the region of strong dynamical stability is typically
separated from the region of strong dynamical instability
by an open set in parameter space $\alpha:=(U_-, U_+, s)$\footnote{
More precisely, the subset of $\alpha$ corresponding to Lax $p$-shocks
on which $\Delta$ was defined, or an open neighborhood 
with $\Delta$ extended by \eqref{Delta}.
}
of indeterminate stability, 
on which $\Delta(\tilde \xi, \lambda)$ has roots $\lambda$
lying precisely on the imaginary axis, 
but none with $\R \lambda>0$;
{\it thus, the point of transition from stability to instability is not
determined in the hyperbolic stability theory.}
A similar situation holds in MHD [Bl, BT.1--4, BTM.1--2].
The reason for this at first puzzling phenomenon,
as described in [BRSZ, Z.4], is that,
under appropriate normalization,
$\Delta(\tilde \xi,\cdot)$ takes imaginary $\lambda=i\tau$ to imaginary
$\Delta$ for $|\tau|$ sufficiently large relative to $|\Txi|$: more precisely,
for $(\tilde \xi, \tau)$ lying in the ``hyperbolic regions'' ${\cal{H}}_\pm$
bounded by glancing sets ${\cal{G}}(P_\pm)$.
Thus, fixing $\tilde \xi$, and explicitly noting the
dependence on parameters $\alpha$,
we find that imaginary zeroes $\lambda=i\tau$ 
of $\Delta^\alpha(\tilde \xi, \cdot)$ for which
$(\tilde \xi, \tau)\in {\cal{H}}_\pm$,
such as occur for gas dynamics and MHD,
generically persist under perturbation in $\alpha$,
rather than moving into the stable or unstable complex
half-planes $\R \lambda <0$ and $\R \lambda >0$ as would otherwise
be expected; see exercise \ref{supersonic} below.
}
\end{rems}

\begin{ex}
\label{hersch}
Under (A1)--(A2), (H2), show that 
$\det \big(\CalA_\pm(\Txi, \lambda) - i\xi_1 \big)=0$
implies $\det (\sum_j i\xi_j A^j + \lambda)=0$, 
$\CalA_\pm$ defined as in \eqref{CalAhyp}, violating
hyperbolicity of \eqref{hyperbolic} if $\R \lambda \ne 0$,
$\Im \xi_1=0$.  Thus, $\CalA_\pm$ has no center subspace
on the domain $\R \lambda > 0$.  
It follows by standard matrix perturbation theory [Kat] that the
stable and unstable subspaces of $\CalA_\pm$ have constant dimension
on $\R \lambda>0$, and the associated $\CalA_\pm$-invariant
projections are analytic in $(\Txi, \lambda)$.
This implies the existence of analytic bases $r_j^\pm$, by 
a lemma of Kato asserting
that analytic projections induce analytic bases on 
simply connected domains;
see [Kat], pp. 99--102.
\end{ex}

\begin{ex}
\label{supersonic}
Let $f(\alpha,\lambda)$ be continuous in $\alpha$ and $\lambda\in \CC$,
with the additional property
that $f(\alpha, \cdot): \RR \to \RR$.
(a) Show that in the vicinity of any root $(0,\lambda_0)$, $\lambda_0$ real,
such that $f_\lambda(0,\lambda_0)$ exists and is nonzero,
there is a continuous family of roots $(\alpha, \lambda(\alpha))$
with $\lambda(\alpha)$ real (Intermediate Value Theorem).
More generally, odd topological degree of $ f(0, \cdot )$ at $\lambda=\lambda_0$
is sufficient;
for even degrees
$f(\alpha, \lambda)= (\alpha + \lambda^2)^m$,
$\alpha \in \RR$, $\lambda_0=0$
is a counterexample. 
(b) If $f$ is analytic in $\lambda$, show that 
$\bar f(\alpha, \bar \lambda)= f(\alpha, \lambda)$,
where $\bar z$ denotes complex conjugate of $z$,
hence nonreal roots occur in conjugate pairs.
\end{ex}

\medbreak
{\bf 2.3. Viscous stability conditions.}
We'll both augment and refine the classical stability conditions
through a viscous stability analysis, at the same time providing
their rigorous justification.

\medbreak
{\bf 2.3.1. Spectral stability.}  
We begin by adding to the conditions of structural and
dynamical stability a third condition of (viscous) spectral stability.
Let $L$ denote the linearized operator about the wave, 
i.e.,  the generator of the linear equations
$U_t=LU$ obtained by linearizing \eqref{viscous} about the
profile $\bar U$, and $L_\Txi$ the family of operators 
obtained from $L$ by Fourier transform in the directions
$\tilde x$ parallel to the front, indexed by frequency $\Txi$.
Explicit representations of $L$ and $L_\Txi$ are given in Section 3.2.

\begin{defi}\label{spectralstab}
\textup{
Similarly as in the hyperbolic stability theory, we
define {\it weak spectral stability} (clearly necessary
for viscous stability) as the absence of
unstable $L^2$ spectrum $\R \lambda>0$ for the 
linearized operator $L$ about the wave, or equivalently
\begin{equation}
\label{weakspec}
\lambda \not \in \sigma(L_\Txi)
\quad \text{\rm for $\Txi\in \RR^{d-1}$ and $\Re \lambda >0$.}
\end{equation}
We define {\it strong spectral stability} (neither necessary
nor sufficient for viscous stability) as
\begin{equation}
\label{strongspec}
\lambda \not \in \sigma(L_\Txi)
\quad \text{\rm for $\Txi\in \RR^{d-1}$, $\Re \lambda >0$, 
and $(\Txi,\lambda)\ne (0,0)$.}
\end{equation}
We define {\it strong spectral instability} as failure of \eqref{weakspec}.
}
\end{defi}

Conditions \eqref{weakspec} and \eqref{strongspec}
may equivalently be expressed in terms of the {\it Evans function} 
$D(\tilde \xi, \lambda)$ (defined Section 5), 
a spectral determinant analogous to the Lopatinski determinant of
the hyperbolic case, as
\begin{equation}
\label{weakevans}
D(\tilde \xi, \lambda)\ne 0 \quad \text{\rm for $\tilde \xi\in \RR^{d-1}$
and $\R \lambda >0$}
\end{equation}
and
\begin{equation}
\label{strongevans}
D(\tilde \xi, \lambda)\ne 0 \quad \text{\rm for $\tilde \xi\in \RR^{d-1}$,
$\R \lambda \ge 0$, and $(\tilde \xi, \lambda)\ne (0,0)$,}
\end{equation}
respectively, where
zeroes of $D$ correspond to normal modes, or solutions
$U=e^{\lambda t}e^{i\tilde \xi \cdot \tilde x} {\hat{\hat U}}(x_1)$
of the linearized perturbation equations: equivalently,
eigenvalues $\lambda$ of the operator $L_{\tilde \xi}$ 
obtained from $L$ by Fourier transform in the directions
$\tilde x$ parallel to the front.
Under assumptions (A1)--(A2), 
$D$ may be chosen analytically in $(\Txi, \lambda)$ on 
$\{\Txi, \lambda:\, \R \lambda \ge 0\}\setminus \{(0,0)\}$;
see [ZS, Z.3--4] or Section 5 below.

\begin{rems}
\label{spectrem}
\textup{
1. It is readily verified that $L\bar U'=L_0 \bar U' =0$, a consequence of 
translational invariance of the original equations \eqref{viscous}, 
from which we obtain $0\in \sigma_{\ess}(L)$ and in particular $D(0,0)=0$.
The exclusion of the origin in definition \eqref{strongevans} 
is therefore necessary for the application to 
stability of viscous shock profiles.
If there held the stronger condition of uniform spectral stability,
\begin{equation}
\label{hille}
\R \sigma(L)<0,
\end{equation}
we could conclude exponential stability of the linearized solution operator
$e^{LT}$ by the generalized Hille--Yosida theorem 
(Proposition \ref{semigpcriteria},
Appendix A) together with the high-frequency
resolvent bounds we obtain later.
However, in the absence of a spectral gap between $0$ and $\R \sigma (L)$,
even bounded linearized stability is a delicate question.
Moreover, we here discuss nonlinear asymptotic stability, for
which we require a rate of linearized decay sufficient to
carry out a nonlinear iteration.
This is connected with the rate at which the spectrum of
of $L_{\Txi}$ moves into the stable complex half-plane $\R \lambda <0$
as $\Txi$ is varied about the origin, 
information that is encoded in the
conditions of structural and refined dynamical stability 
defined just below.
}
\smallbreak
\textup{
2. In the one-dimensional case $\tilde \xi\equiv 0$, strong spectral
stability is roughly equivalent to linearized stability with respect
to zero-mass initial perturbations; see [ZH, HuZ, MaZ.2--4].
Zero-mass stability 
is exactly the property that 
was used by Goodman and Xin [GoX] to rigorously justify the 
small-viscosity
matched asymptotic expansion about small-amplitude shock waves
in one dimension; indeed, their argument implies validity of
matched asymptotic expansion in one dimension for any
structurally and dynamically stable wave about which
the associated profiles exhibit zero-mass stability.
}
\end{rems}
\medbreak

{\bf 2.3.2. Refined dynamical stability conditions.}
A rigorous connection between viscous stability
and the formal conditions of structural and dynamical
stability is given by the following fundamental result
(proved in the Lax case in Section 5).

\begin{prop}[{[ZS, Z.3, M\'eZ.2]\footnote{
Proved in [ZS] in the strictly hyperbolic case, along rays through
the origin lying outside the glancing set; full version proved in [Z.3]
under the additional hypothesis (H5), and in [M\'eZ.2] for the
general case.
}
 }]
\label{ZS}
Given (A1)--(A2) and (H0)--(H4), 
and under appropriate normalizations of $D$ and $\Delta$,
\begin{equation}
\label{tangent}
D(\tilde \xi, \lambda) = 
\overbrace{\gamma \Delta(\tilde \xi, \lambda)
}^{{\cal{O}}(|\tilde \xi, \lambda|^\ell)} +
{{o}}(|\tilde \xi, \lambda|^{\ell}),
\end{equation}
uniformly on $\tilde \xi\in \RR^{d-1}$, $\R \lambda \ge 0$,
for $\rho:=|\tilde \xi, \lambda|$ sufficiently small,
where $\gamma$ is a constant measuring transversality of 
connection $\bar U$ as a solution
of the associated traveling-wave ODE,
$\Delta$ is an appropriate Lopatinski condition, 
and $\ell$ is as in (H3).
In the present, Lax case, $\Delta$ is as in \eqref{Delta} and $\ell=1$.
\end{prop}

\begin{rem}
\label{rays}
\textup{
In the absence of (H4), \eqref{tangent} still holds
along individual rays through the origin, provided that they 
lie in directions of analyticity for $\Delta$, in particular
for $\R \lambda>0$, by the original argument of [ZS].
This is sufficient for the necessary stability conditions derived
below; however, uniformity is essential for our sufficient stability
conditions.\footnote{
This point was not clearly stated in [Z.3]; indeed, 
\eqref{tangent} is hidden
in the detailed estimates of the technical Section 4 of that reference,
which depend on the additional hypothesis (H5).
For a more satisfactory treatment, see [M\'eZ.2, GMWZ.3--4, Z.4].
}
}
\end{rem}

\medbreak

\begin{defi}
\label{refstab}
\textup{
We define {\it refined weak dynamical stability} 
as \eqref{weaklop} augmented with
the second-order condition $ \R\beta(\tilde \xi,i\tau) \ge 0 $
for any real $\Txi$,  $\tau$ such that $\Delta$ and 
$D^{\Txi, \lambda}(\rho):=D(\rho \Txi,\rho \lambda)$ are 
analytic at $(\Txi, i\tau)$ and $(\Txi,i\tau,0)$, respectively,
with $\Delta(\Txi, i\tau)=0$ and 
$\Delta_\lambda(\Txi,i\tau)\ne 0$, where
\begin{equation}
\label{beta}
\beta(\Txi,i\tau):= 
{ (\partial/\partial \rho )^{\ell+1}D(\rho\Txi,\rho i\tau)|_{\rho=0}
\over
(\partial/\partial \lambda)  \bar\Delta(\Txi,i\tau)}
=
{(\partial/\partial \rho )^{\ell+1}D(\rho\Txi,\rho \lambda)
\over
(\partial/\partial \lambda)
(\partial/\partial \rho)^{\ell}
D(\rho\Txi,\rho \lambda)}_{|_{\rho=0, \lambda=i\tau}}. 
\end{equation}
}

\textup{
Under (A1)--(A2), (H2), (H4), analyticity of 
$\Delta$ and $D^{\Txi, \lambda}(\rho)$ are equivalent,
and hold for all real
$(\tilde \xi, \tau)$ lying outside the glancing sets
${\cal{G}}(P_\pm)$ [ZS, Z.3--4].
In this case, we define {\it strong refined dynamical stability} 
as \eqref{weaklop} augmented with the conditions that:
(i) $\R\beta(\Txi,i\tau) > 0 $,
$\Delta_\lambda(\tilde \xi, i\tau)\ne 0 $, 
and $\{r_1^-, \dots, r_{p-1}^-, r_{p+1}^+, \dots, r_n^+\}(\tilde \xi, i\tau)$
are independent, where $r_j^\pm$ are defined as in \eqref{Delta}
(automatic for ``extreme'' shocks $p=1$ or $n$),
for any real $(\tilde \xi$, $\tau)\not \in {\cal{G}}(P_\pm)$ 
such that $\Delta(\tilde \xi, i\tau)=0$ and $\Delta$ is analytic,
and (ii) the zero-level set in $(\Txi, \tau)$ of $\Delta(\Txi, i\tau)=0$
intersects the glancing sets ${\cal{G}}_\pm$ transversely
in $\RR^{d-1}\times \RR$ (in particular, their intersections are trivial
in dimension $d=2$).
}
\end{defi}

\begin{rem}
\label{correction}
\textup{
The coefficient $\beta$ may be recognized as a viscous correction
measuring departure from homogeneity: note that $\beta$ would
vanish if $D$ were homogeneous degree $\ell$.
It has a physical interpretation as an ``effective viscosity'' 
governing transverse propagation of deformations in the front 
[Go.2, GM, ZS, Z.3--4, HoZ.1--2].
}
\end{rem}

\medbreak
{\bf 2.3.3. Main results.}
With these definitions, we can now state our main results,
to be established throughout the remainder of the article.

\begin{theo}[{[ZS, Z.3]}]
\label{mainnec}
Given (A1)--(A2), (H0)--(H3), structural stability $\gamma \ne 0$,
and one-dimensional inviscid stability $\Delta(0,1)\ne 0$,
weak spectral stability and weak refined dynamical stability
are necessary for bounded $C^\infty_0 \to L^p$, 
viscous stability in all dimensions $d\ge 2$, for any $1\le p\le \infty$,
and for profiles of any type.
\end{theo}

\begin{theo}[{[Z.3--4]}]
\label{mainsuf}
Given (A1)--(A2), (H0)--(H5),
strong spectral stability plus structural and strong refined
dynamical stability are sufficient for asymptotic
$L^1\cap H^s\to H^s\cap L^p$ viscous stability of pure, Lax-type profiles
in dimensions $d\ge 3$, for $q(d)\le s\le q$, $q$ and
$q(d)$ as defined in (H0),
and any $2\le p\le \infty$, 
or linearized viscous stability in dimensions $d\ge 2$, for $1\le s\le q$
and $2\le p\le \infty$,
with rate of decay (in either case)
\begin{equation}
\label{rate}
|U(t)-\bar U|_{H^s} \le C(1+t)^{-(d-1)/4+ \epsilon}|U(0)-\bar U|_{L^1\cap H^s}
\end{equation}
approximately equal to that of a $(d-1)$-dimensional heat kernel,
for any fixed $\epsilon>0$ and 
$|U(0)-\bar U|_{L^1\cap H^s}$ sufficiently small.
Strong spectral plus structural and strong (inviscid, not refined) 
dynamical stability are sufficient for nonlinear
stability in all dimensions $d\ge 2$, with rate of decay
\eqref{rate}, $\epsilon=0$, exactly equal to that of
a $(d-1)$-dimensional heat kernel.
Similar results hold for profiles of nonclassical, 
over- or undercompressive type, as described in [Z.3--4].
\end{theo}

\begin{rems}\label{endtwo}
\textup{
1. Rate \eqref{rate} with $\epsilon=0$ is sharp for weakly inviscid stable
shocks, as shown by the scalar case [Go.3, GM, HoZ.1--2].
(As pointed out by Majda, scalar shocks are always weakly inviscid stable 
[Ma.1--3, ZS, Z.3].)
}

\textup{
2. Given (A1)--(A2), (H0)--(H3), 
$\gamma \Delta(0,1)\ne 0$ is necessary for one-dimensional 
viscous stability [ZH, MaZ.3] with respect to $C^\infty_0$
perturbations. (Note: $C^\infty_0(x_1)\not \subset C^\infty_0(x)$.)
}

\textup{
3. Under (A1)--(A2), (H0)--(H3) plus the mild additional 
assumptions that the principal characteristic speed $a_p(\xi)$
be simple, genuinely nonlinear, and strictly convex (resp. concave) 
with respect to $\tilde \xi$ at $\xi=(\xi_1, \tilde \xi)=(1,0_{d-1})$,
the sufficient stability conditions
of Theorem \ref{mainsuf} are satisfied for sufficiently small-amplitude
shock profiles [FreS, PZ].
On the other hand, the necessary conditions of Theorem
\ref{mainnec} are known to fail for certain large-amplitude shock
profiles under (A1)--(A2), (H0)--(H3) [GZ, FreZ, ZS, Z.3--4].
}

\textup{
4. The boundary between the necessary conditions of Theorem
\ref{mainnec} and the sufficient conditions of Theorem \ref{mainsuf}
is generically of codimension one (exercise; see [Z.4], Remark
1.22.2).  Thus, transition from viscous stability to instability is
in principle determined, in contrast to the situation of the inviscid case
(Remark \ref{majda}).
}

\textup{
5. It is an important open problem which of the conditions of
refined dynamical and structural stability in practice
determines the transition from viscous stability to instability
as shock strength is varied; see [MaZ.4, Z.3--4] for further discussion.
}
\end{rems}


\bigbreak
\clearpage

\section{Analytical preliminaries}

The rest of this article is devoted to the proof of Theorems
\ref{mainnec} and \ref{mainsuf}.
We begin in this section by assembling some needed background results.
\medbreak

{\bf 3.1. Profile facts.}
Consider the standing wave ODE
\begin{equation}
\label{ode}
\begin{aligned}
(F^1)^I(U)'&=0,\\
(F^1)^{II}(U)'&=((B^1)^{II} U')',\\
\end{aligned}
\end{equation}
$U^I$, $(F^1)^I \in \RR^{n-r}$, $U^{II}$, $(F^1)^{II}\in \RR^r$,
where ``$'$'' denotes $d/dx$ and superscripts $I$ and ${II}$ refer
respectively to first- and second-block rows.
Integrating from $-\infty$
to $x$ gives an implicit first-order system
\begin{equation}
\label{II}
(B^1)^{II}U'= (F^1)^{II}(U) -(F^1)^{II}(U_-)
\end{equation}
on the $r$-dimensional level set
\begin{equation}
\label{I}
(F^1)^{I}(U) \equiv (F^1)^{I}(U_-).
\end{equation}

\begin{lem}[{[MaZ.3]}]\label{odelem}
Given (A1), the weak version $\det A^1_{11}\ne 0$ of (H1)
is equivalent to the property that \eqref{I} determines a nondegenerate
$r$-dimensional manifold
on which \eqref{II} determines a nondegenerate first-order ODE.
Moreover, assuming (A1)--(A2) and $\det A^1_{11}\ne 0$,
(H2) is equivalent to hyperbolicity of rest states $U_\pm$.
\end{lem}

\begin{rem}
\label{oderem}
\textup{
This result motivated the introduction of the weak
form $\det A^1_{11}\ne 0$ of (H1) in [Z.3, MaZ.3]; the importance of the
strong form $A^1_{11}>0$ ($<0$) was first pointed out in [MaZ.4].
}
\end{rem}
\medbreak

\begin{cor}[{[MaZ.3]}]\label{odecor}
Given (A1)-(A2), and (H0)--(H2)
\begin{equation}
\label{expdecay}
|\partial_x^k(\bar U-U_\pm)|\le Ce^{-\theta|x|}, 
\qquad x\gtrless 0,
\end{equation}
for $0\le k \le q+2$, where $q$ is as defined in (H0).
\end{cor}

\begin{proof} [Proof of Lemma \ref{odelem}]
This result was proved in somewhat greater generality in [MaZ.3].
Here, we give a simpler proof making use of structure (A1)--(A2).
From (A1), it can be shown (exercise; see [Z.4], appendix A1) that
$\partial U/\partial W$ is lower block-triangular and $(B^1)^{II}
(\partial U/\partial W)=\bdiag\{0, \hat b\}$, $\hat b$
nonsingular, with $\partial (F^1)^I/\partial w^I\ne 0$ 
if and only if $\tilde A^1_{11}\ne 0$,
with $\tilde A^1$ as in (A1).
Considering \eqref{I}--\eqref{II} with respect to coordinate
$W$, we readily obtain the first assertion, with resulting ODE 
\begin{equation}
\label{wode}
\hat b (w^{II})'= (F^1)^{II}(U(w^I(w^{II}),w^{II})- (F^1)(U_-),
\end{equation}
where $w^{I}(w^{II})$ is determined by the Implicit Function Theorem
using relation $(F^1)^I(U(W))\equiv (F^1)^{I}(U_-)$.

Linearizing \eqref{wode} about rest point $W=W_\pm$
and using the Implicit Function Theorem to 
calculate $\partial w^{I}/\partial w^{II}$,
we obtain
\begin{equation}
\label{linwode}
(w^{II})'= M_\pm w^{II}:=
(\hat b)^{-1} (-\alpha_{21}\alpha_{11}^{-1}\alpha_{12}+ \alpha_{22})_\pm
w^{II},
\end{equation}
where $\alpha_\pm:=(\partial F^1/\partial W)_\pm$, and $\hat b_\pm$ are
evaluated at $W_\pm$.
Noting that $\det (-\alpha_{21}\alpha_{11}^{-1}\alpha_{12}+ \alpha_{22})=
\det \alpha/\det\alpha_{11}$, with $\alpha_{11}\ne 0$
as a consequence of $\det \tilde A^1_{11}\ne 0$, 
we find that the coefficient matrix $M_\pm$ of the 
linearized ODE has zero eigenvalues if and only
if $\det \alpha_\pm$, or equivalently $\det A^1_\pm$ vanishes.
On the other hand, existence of nonzero pure imaginary eigenvalues $i\xi_1$
would imply, going back to the original equation \eqref{symmhypparab}
linearized about $W_\pm$, that
$\det (i\xi_1\tilde A^1 -\xi_1^2 \tilde B^{11})_\pm=0$ for some real
$\xi_1\ne 0$.
But, this is precluded by (A1)--(A2) (nontrivial exercise; see
(K3), Lemma \ref{KSh}, below).
Thus, the coefficient matrix $M$ has a center subspace if and only
if (H2) fails, verifying the second assertion.
\end{proof}

\begin{proof}[Proof of Corollary \ref{odecor}]
Standard ODE estimates.
\end{proof}

For general interest, and to show that our assumptions on
the profile are not vacuous, we include without proof the following
generalization of a result of [MP] in the strictly parabolic case,
relating structure of a viscous profile to its hyperbolic type.

\begin{lem}[{[MaZ.3]}]
\label{type}
Let $d_+$ denote the dimension of the stable subspace
of $M_+$ and
$d_-$ the dimension of the unstable subspace of $M_-$,
$M_\pm$ defined as in \eqref{linwode}.
Then, existence of a connecting profile $\bar U$ together with
(A1)--(A2) and (H1)--(H2), implies that
\begin{equation}
\label{indrel}
\hat \ell:= d_++d_--r=i_++i_--n,
\end{equation}
where $i_+$ denotes the dimension of the stable subspace of $A^1_+$
and
$i_-$ denotes the dimension of the unstable subspace of $A^1_-$;
in the Lax case, $\hat \ell= 1$.
\end{lem}

\begin{rem}
\label{ell}
\textup{
$\hat \ell=\ell$ in the case
of a transverse profile, $\ell$ as in (H3),
whereas $i_++i_--n=1$ may be 
recognized as the Lax characteristic condition (L).
Assumption \ref{lax} could be rephrased in this case as $\hat\ell=
1= i_++i_--n$.
}
\end{rem}

\begin{proof}
See [MaZ.3], Appendix A1, or [Z.4], Appendix A2.
\end{proof}

\medbreak

{\bf 3.2. Spectral resolution formulae.}
Consider the linearized equations
\begin{equation}
\label{linearized}
\begin{aligned}
U_t&=LU:= -\sum_j (A^j U)_{x_j}+ \sum_{jk} (B^{jk}U_{x_k})_{x_j},
\qquad
U(0)=U_0,
\end{aligned}
\end{equation}
where
\begin{equation}
\label{coeffs}
B^{jk}:= B^{jk}(\bar U(x_1)),
\qquad
A^jv:= dF^j(\bar U(x_1))v - (dB^{j1}(\bar U(x_1))v)\bar U'(x_1).
\end{equation}

Since the coefficients depend only on $x_1$, it is natural
to Fourier-transform in $\tilde x$, $x=(x_1, \dots, x_d)=:(x_1, \tilde x)$,
to reduce to a family of partial differential equations (PDE) 
\begin{equation}
\label{FTlinearized}
\begin{aligned}
\hat U_t&=L_\Txi \hat U := 
\overbrace{ (B^{11}\hU')'-(A^1\hU)'}^{L_0\hU} 
 -i \sum_{j\not= 1}A^j \xi_j \hU 
+ i\sum_{j\not= 1} B^{j1}\xi_j \hU' \\
&\quad + i\sum_{k\not= 1}(B^{1k}\xi_k \hU)' 
-\sum_{j,k\not= 1} B^{jk}\xi_j \xi_k \hU,
\qquad \hat U(0)=\hat U_0 
\end{aligned}
\end{equation}
in $(x_1,t)$ indexed by frequency $\tilde \xi\in \RR^{d-1}$,
where ``$\, '\, $'' denotes $\partial / \partial x_1$ and
$\hat U=\hat U(x_1, \tilde \xi, t)$ denotes the Fourier transform
of $U=U(x,t)$.

Finally, taking advantage of autonomy of the equations, 
we may take the Laplace transform in $t$ to reduce, formally,
to the resolvent equation
\begin{equation}
\label{LT}
(\lambda - L_\Txi){\hat {\hat U}}= \hat U_0,
\end{equation}
where ${\hat {\hat U}}(x_1, \Txi, \lambda)$ 
denotes the Laplace--Fourier transform 
of $U=U(x,t)$ and $\hat U_0(x_1)$ denotes the Fourier transform 
of initial data $U_0(x)$.
A system of ODE in $x_1$ indexed by $(\Txi, \lambda)$,
\eqref{LT} may be estimated sharply using either explicit
representation formulae/variation of constants for systems
of ODE or Kreiss symmetrizer techniques; see, e.g., [Z.3--4]
and [GMWZ.2], respectively.
The following proposition makes rigorous sense of this procedure.
%

\begin{prop}[{[MaZ.3, Z.4]}]\label{specres}
Given (A1), $L$ generates a $C^0$ semigroup
$|e^{L t}|\le Ce^{\gamma_0 t}$ on $L^2$
with domain $\CalD(L):=\{ U: \, U,\, LU \in L^2\}$
($Lu$ defined in the distributional sense),
satisfying the generalized spectral resolution (inverse
Laplace--Fourier transform) formula
\begin{equation}
\label{ILFT}
\begin{aligned}
e^{L t}f(x)=
\pv \int_{\gamma -i\infty}^{\gamma+i\infty}
\int_{\RR^{d-1}} 
e^{i\tilde \xi \cdot \tilde x + \lambda t}
(\lambda- L_\Txi)^{-1}\hat f(x_1, \Txi) \, d\tilde \xi d\lambda 
\end{aligned}
\end{equation}
for $\gamma>\gamma_0$, $t\ge 0$, and
$f \in \CalD(L)$, where
$\hat f$ denotes Fourier transform of $f$.
Likewise,
\begin{equation}
\label{inhomILFT}
\begin{aligned}
\int_0^t e^{L(t-s)}f(s) \, ds &=
\pv \int_{\gamma -i\infty}^{\gamma+i\infty}
\int_{\RR^{d-1}} 
e^{i\tilde \xi \cdot \tilde x + \lambda t}
(\lambda- L_\Txi)^{-1}
{\widehat {\widehat {f^T}}}(x_1, \Txi, \lambda) \, d\tilde \xi d\lambda 
\end{aligned}
\end{equation}
for $\gamma>\gamma_0$, $0\le t\le T$,
and $f\in L^1([0,T]; \CalD(L))$, 
where $\hat{\hat g}$ denotes Laplace--Fourier transform of $g$ and
\begin{equation}
\label{truncation}
f^T(x,s):=
\begin{cases}
f(x,s) & \text{\rm for $0\le s\le T$}\\
0 & \text{\rm otherwise}.\\
\end{cases}
\end{equation}
\end{prop}


\begin{rem}
\label{inhomILFTrem}
\textup{
Bound \eqref{inhomILFT} concerns the
inhomogeneous, zero-initial-data problem
$U_t-LU=f$, $U(0)=0$.
Provided that $L$ generates a $C^0$ semigroup $e^{LT}$,
the ``mild'', or semigroup
solution of this equation is defined (see, e.g., [Pa], pp. 105--110)
by Duhamel formula
\begin{equation}
\label{specduhamel}
U(t)=\int_0^t e^{L(t-s)}f(s) \, ds.
\end{equation}
This corresponds to a solution in the usual
weak, or distributional sense, with somewhat stronger regularity in $t$.
Formally, ${\hat{\hat{V}}}:=(\lambda-L)^{-1}
{\widehat{\widehat{f^T}}}$ describes the Laplace--Fourier
transform of the solution of the truncated-source
problem $V_t-LV=f^T$, $V(0)=0$.  Evidently,
$U(s)=V(s)$ for $0\le s\le T$.
}
\end{rem}

\begin{proof}
Let $u_n\to u$ and $Lu_n \to f$, $u_n\in \CalD(L)$,   
$u$, $f\in L^2$, ``$\to$'' denoting convergence in $L^2$ norm.
Then, $Lu_n \rightharpoondown Lu$ by $u_n\to u$,
where $Lu$ is defined in the distributional sense, but also
$Lu_n\rightharpoondown f$ by $Lu_n \to f$,
whence $Lu=f$ by uniqueness of weak limits, and therefore $u\in \CalD(L)$.
Moreover, $\CalD(L)$ is dense in $L^2$, since $C^\infty_0\subset \CalD(L)$.
Thus, 
$L$ is a closed, densely defined operator on $L^2$ 
with domain $\CalD(L)$;
see Definition \ref{closed}, Appendix A. 

By standard semigroup theory, 
therefore (Proposition \ref{semigpcriteria}, Appendix A),
$L$ generates a $C^0$-semigroup 
$|e^{LT}|\le Ce^{\gamma_0t}$
if and only if it satisfies resolvent estimate
\begin{equation}
\label{resk}
|(\lambda- L)^{-k}|_{L^2} \le
\frac{C}{|\lambda - \gamma_0|^k}
\end{equation}
for some uniform $C>0$,
for all real $\lambda >\gamma_0$ and all $k\ge 1$,
or equivalently
\begin{equation}
\label{enk}
|U|_{L^2} \le
\frac{C |(\lambda- L)^{k}U|_{L^2}}
{|\lambda - \gamma_0|^k}
\end{equation}
for all $k\ge 1$ and 
\begin{equation}
\label{range}
\range (\lambda-L)=L^2.
\end{equation}

In the present case, the latter condition is superfluous, 
since it is easily verified by results of Henry [He, BSZ] that
$\lambda -L$ is Fredholm for sufficiently large real $\lambda$.
More generally,
\eqref{range} may be verified by a corresponding estimate 
\begin{equation}
\label{adjen}
|U|_{L^2} \le
\frac{C |(\lambda- L^*)U|_{L^2}}{|\lambda - \gamma_0|}
\end{equation}
on the adjoint operator $L^*$, 
equivalent to \eqref{range} given \eqref{enk} 
(or, given \eqref{range}, to \eqref{enk}, $k=1$)
and often available (as here) by the same techniques;
see Corollary \ref{adjver}, Appendix A.

A sufficient condition for \eqref{enk}, and the standard
means by which it is proved, is
\begin{equation}
\label{res}
|W|_{L^2} \le
\frac{ |S(\lambda- L)S^{-1}W|_{L^2}}
{|\lambda - \gamma_0|}
\end{equation}
for all real $\lambda >\gamma_0$, 
for some uniformly invertible transformation $S$: in this case
the change of coordinates 
$S: U \to W$ defined by 
\begin{equation}
\label{newcoords}
W:= (\tilde A^0)^{1/2}(\partial W/\partial U)(\bar U(x_1)) U.
\end{equation}
Once \eqref{enk}--\eqref{range} have been established, we obtain 
by standard properties of $C^0$ semigroups the bounds
$|Se^{L t}S^{-1}|\le e^{\gamma_0 t}$ and 
$|e^{L t}|\le Ce^{\gamma_0 t}$ and also the inverse Laplace transform formulae
\begin{equation}
\label{ILT1}
e^{Lt}f(x)=
\pv \int_{\gamma -i\infty}^{\gamma+i\infty}
e^{\lambda t}
(\lambda- L)^{-1}f \, d\lambda
\end{equation}
for $f \in \CalD(L)$ and
\begin{equation}
\label{inhomILT1}
\int_0^T e^{L(T-t)}f(x)=
\pv \int_{\gamma -i\infty}^{\gamma+i\infty}
e^{\lambda t}
(\lambda- L)^{-1}\widehat{f^T}(\lambda) \, d\lambda
\end{equation}
for $f\in L^1([0,T];\CalD(L))$, for any $\gamma>\gamma_0$,
where $\hat f$ denotes Laplace transform of $f$;
see Propositions \ref{semigpcriteria}, \ref{ILTprop}, and 
\ref{inhomspecres},
Appendix A.
Formulae \eqref{ILFT}--\eqref{inhomILFT} then follow from 
\eqref{ILT1}--\eqref{inhomILT1} by
Fourier transform/distribution theory.
Thus, it remains only to verify \eqref{res} in order to complete the proof.

To establish \eqref{res}, rewrite resolvent equation
$(\lambda-L)U= f$ in $W$-coordinates, as
\begin{equation}
\label{linearizedW}
\lambda \tilde A^0 W + \sum_j \tilde A^j W_{x_j}
+ \sum_{jk} (\tilde B^{jk}W_{x_k})_{x_j}+ \tilde CW 
+ \sum_j \tilde D^j w^{II}_{x_j} = \tilde A^0 \tilde f,
\end{equation}
where $\tilde f:= (\partial W/\partial U)(\bar U(x_1))f$,
$C$ is uniformly bounded, and coefficients $\tilde A^j$, $\tilde B^{jk}$
satisfy (A1) with $\tilde G\equiv 0$.
(Exercise: check that properties (A1) are preserved up to
lower-order terms $\tilde CW+ \sum_j \tilde D^j w^{II}_{x_j}$,
using commutation of coordinate-change and linearization.)
Using block-diagonal form of $\tilde A^0$, we arrange by
coordinate change $W\to (\tilde A^0)^{1/2}W$ if necessary 
that $\tilde A^0=I$ without loss of generality;
note that lower-order commutator terms arising through change
of coordinates are again of the form 
$\tilde CW+ \sum_j \tilde D^j w^{II}_{x_j}$,
hence do not change the structure asserted in \eqref{linearizedW}.

Taking the real part of the complex inner product of $W$ against
\eqref{linearizedW}, and using \eqref{goodbtilde}, \eqref{garding}
together with symmetry of $\tilde A^j_{11}$ terms and Young's
inequality, similarly as in the proof of Proposition \ref{parabwellposed},
we obtain
\begin{equation}
\R \lambda |W|_{L^2} + \theta|w^{II}|_{H^1}^2 \le
\gamma_0 |W|_{L^2}^2 + |\tilde f|_{L^2}|W|_{L^2}
\end{equation}
for $\gamma_0>0$ sufficiently large.
Dropping the favorable $w^{II}$ term and dividing both sides by
$|W|_{L^2}$, we obtain
\begin{equation}
(\R \lambda- \gamma_0) |W|_{L^2} \le |\tilde f|_{L^2}
\end{equation}
Noting that $S(\lambda-L)S^{-1} W= \tilde f$, we are done.
\end{proof}

\begin{rem}
\label{linearexistence}
\textup{
The properties and techniques used to establish Proposition
\ref{specres} are linearized versions of the ones used to establish 
Proposition \ref{parabwellposed}.  For a still more direct link,
see Remark \ref{lumer} and Exercise \ref{hilex}, Appendix A,
concerning the Lumer--Phillips Theorem, which assert that,
in the favorable $W$-coordinates, existence of a $C^0$ contraction
semigroup $|e^{Lt}|\le e^{\gamma_0 t}$ is equivalent to the a priori estimate
\begin{equation}
\label{apriori}
(d/dt)(1/2)|W|_{L^2}^2=\R \langle u, Lu\rangle \le \gamma_0|W|^2,
\end{equation}
the autonomous version of \eqref{0friedrichs}.
That is, Proposition \ref{specres}
essentially concerns local linearized well-posedness, with no global stability
properties either asserted or required so far.
}

\textup{
Besides validating the spectral resolution formulae, the
semigroup framework gives also a convenient means to generate
solutions of the nonautonomous linear equations arising in
the nonlinear iteration schemes described in Section 1,
based only on the already-verified instantaneous version of \eqref{apriori};
specifically, Proposition \ref{evolutionprop}, 
Appendix A, assuming \eqref{apriori},
guarantees a solution satisfying $|u(t)|_{L^2}\le e^{\gamma_0 t}|u(0)|_{L^2}$.
A standard way of constructing the autonomous semigroup is by 
semidiscrete approximation
using the first-order implicit Euler scheme (see Remark \ref{euler}),
the nonautonomous solution operator then being approximated by
a concatenation of frozen-coefficient semigroup solutions on each 
mesh block; see [Pa].
%
It is interesting to compare this approach to Friedrichs' original
construction of solutions by finite difference approximation [Fr].
}
\end{rem}

%
%

\medbreak
{\bf 3.3. Asymptotic ODE theory: the gap and conjugation lemmas.}
Consider a general family of first-order ODE 
\begin{equation}
\label{gfirstorder}
\WW'-{\mathbb A}(x, \Lambda)\WW=\FF
\end{equation}
indexed by a spectral parameter $\Lambda \in \Omega \subset \CC^m$, where
$W\in \CC^N$, $x\in \RR$ and ``$'$'' denotes $d/dx$.

\begin{exams}
\textup{
1. Eigenvalue equation  
$(L_\Txi-\lambda)U=0$, written in phase coordinates
$\WW:=(W, (w^{II})')$, $W:=(\partial W/\partial U)(\bar U(x))U$, 
with $\Lambda:=(\tilde \xi, \lambda)$ and $\FF:=0$.
}

\textup{
2. Resolvent equation  
$(L_\Txi-\lambda)U=f$, written in phase coordinates
$\WW:=(W, (w^{II})')$, $W:=(\partial W/\partial U)(\bar U(x))U$, 
and $|\FF|={\cal{O}}(|f|)$.
}
\end{exams}

\begin{ass}\label{h0}{$\,$}
\medbreak
\textup{(h0) }
Coefficient ${\mathbb A}(\cdot,\Lambda)$, considered
as a function from $\Omega$ into 
$C^0(x)$
is analytic in $\Lambda$.
Moreover, ${\mathbb A}(\cdot, \Lambda)$ approaches
exponentially to limits $\mA_\pm$ as $x\to \pm \infty$, 
with uniform exponential decay estimates
\begin{equation}
\label{expdecay2}
|(\partial/\partial x)^k(\mA- \mA_\pm)| 
\le C_1e^{-\theta|x|/C_2}, \, 
\quad
\text{\rm for } x\gtrless 0, \, 0\le k\le K,
\end{equation}
$C_j$, $\theta>0$, 
on compact subsets of $\Omega $.
\end{ass}
\medbreak

\begin{lem}[{The gap lemma [KS, GZ, ZH]}]
\label{gaplemma}
Consider the homogeneous version $\FF\equiv 0$ of \eqref{gfirstorder}, 
under assumption (h0).
If $V^-(\Lambda)$ is an 
eigenvector of $\mA_-$ with eigenvalue $\mu(\Lambda)$, both 
analytic in $\Lambda$, 
then there exists a solution of \eqref{gfirstorder} of form
\begin{equation}
 \WW(\Lambda, x) = V (x,\Lambda ) e^{\mu(\Lambda) x},
\end{equation}
where $V$ is $C^{1}$ in $x$ and locally analytic in $\Lambda$ and,
for any fixed $\btheta < \theta$, satisfies 
\begin{equation}
\label{3.6g}
V(x,\Lambda )=  V^-(\Lambda ) + \bfO (e^{-\bar \theta|x|}|V^- (\Lambda)|),\quad x < 0.
\end{equation}
\end{lem}

\begin{proof}
Setting $\WW(x)=e^{\mu x}V(x)$, we may rewrite $\WW'=\mA \WW$ as 
\begin{equation}
\label{3.7g}
V' = (\mA_- - \mu I)V+\theta V,
\qquad
\theta
:= (\mA - \mA_-)=\bfO(e^{-\theta|x|}),
\end{equation}
and seek a solution $V(x,\Lambda )\to V^-(x)$ as $x \to \infty$.
Choose $ \btheta< \theta _1 < \theta $
such that there is a spectral gap 
$|\R \big(\sigma \mA_- - (\mu+\theta_1)\big)|>0$
between $\sigma \mA_-$ and $\mu + \theta_1$.
Then, fixing a base point $\Lambda_0$, we can define on some neighborhood
of $\Lambda_0$ to the complementary $\mA_-$-invariant projections
$P(\Lambda)$ and $Q(\Lambda)$ where $P$ projects onto the direct sum
of all eigenspaces of $\mA_-$ with eigenvalues $\Tmu$ 
satisfying 
$ \R(\Tmu) < \R(\mu) + \theta_1, $
and $Q$ projects onto the direct sum of the remaining eigenspaces,
with eigenvalues satisfying 
$ \R(\Tmu)  > \R(\mu) +\theta _1.  $
By basic matrix perturbation theory (eg. [Kat]) it follows that 
$P$ and $Q$ are analytic in a neighborhood of $\Lambda_0$,  
with
\begin{equation}
\label{3.10g}
\left|e^{(\mA_- - \mu I)x} P \right| 
\le C (e^{\theta_1 x}),
\quad x>0, 
\qquad
\left|e^{(\mA_- - \mu I)x} Q \right| 
\le C (e^{\theta_1 x}), 
\quad x<0.
\end{equation}

It follows that, for $M>0$ sufficiently large, the map $\CalT$ defined by 
\begin{equation}
\label{3.11g}
\begin{aligned}
\CalT V(x) 
&= V^- + \int^x_{-\infty} e^{(\mA_- - \mu I)(x-y)} P
\theta (y) V(y) dy \\
&\quad - \int^{-M}_x e^{(\mA_- - \mu I)(x-y)} Q \theta (y) V(y) dy
\end{aligned}
\end{equation}
is a contraction on $L^\infty(-\infty, -M]$.
For, applying \eqref{3.10g}, we have
\begin{equation}
\label{3.12g}
\begin{aligned}
\left|\CalT V_1 - \CalT V_2 \right|_{(x)} 
&\le C |V_1 - V_2|_\infty 
\bigg(\int^x_{-\infty} e^{\theta_1(x-y)} e^{\theta y} dy 
+ \int^{-M}_{x} e^{\theta _1(x-y)}e^{\theta y}dy\bigg)\\
&\le  C_1 |V_1 - V_2|_\infty 
\bigg(e^{\theta_1 x}e^{(\theta-\theta_1)y}|^x_{-\infty}+
e^{\theta_1 x}e^{(\theta-\theta_1)y}|^{-M}_{x}\bigg) \\
&\le C_2  |V_1 - V_2|_\infty e^{-\btheta M} <
\frac{1}{2} |V_1 - V_2|_\infty.
\end{aligned}
\end{equation}

By iteration, we thus obtain a solution $V \in L^\infty (-\infty, -M]$ of $V = 
\CalT V$ with $V\le C_3|V^-|$; since $\CalT$ clearly preserves analyticity
$V(\Lambda, x)$  is 
analytic in $\Lambda$ as the uniform limit of analytic 
iterates (starting with $V_0=0$). 
Differentiation shows that $V$ is a bounded solution of
$V=\CalT V$ if and only if it is a bounded solution of \eqref{3.7g}
(exercise).
Further, taking $V_1=V$, $V_2=0$ in \eqref{3.12g}, we obtain from
the second to last inequality that 
\begin{equation}
\label{3.13}
|V-V^-| = |\CalT(V) - \CalT(0)| \le C_2  e^{\btheta x} |V| 
\le  C_4e^{\btheta x}|V^-|,
\end{equation}
giving \eqref{3.6g}.
Analyticity, and the bounds \eqref{3.6g},  
extend to $x<0$ by standard analytic dependence for the initial value
problem at $x=-M$. 
\end{proof}

\begin{rem}\label{gapexplanation}
\textup{
The title ``gap lemma'' alludes to the fact that we
do not make the usual assumption of
a spectral gap between $\mu(\Lambda)$ 
and the remaining eigenvalues of $\mA_-$, as
in standard results on asymptotic behavior of ODE [Co];
that is, the lemma asserts that exponential
decay of $\mA$ can substitute for a spectral gap.
Note also that we require only analyticity of $\mu$
and not its associated eigenprojection $\Pi_\mu$, allowing crossing
eigenvalues of arbitrary type
(recall, $\Pi_\mu$ is analytic only if $\mu$ is semisimple;
indeed, $\Pi_\mu$ blows up at a nontrivial Jordan block [Kat]).
This is important in the
following application; see Exercise \ref{jordanex} below.
}
\end{rem}
\medbreak

\begin{cor}[{The conjugation lemma [M\'eZ.1]}]
\label{conjugation}
Given (h0), there exist locally to any given $\Lambda_0\in \Omega $
invertible linear transformations $P_+(x,\Lambda)=I+\Theta_+(x,\Lambda)$ and 
$P_-(x,\Lambda) =I+\Theta_-(x,\Lambda)$ defined
on $x\ge 0$ and $x\le 0$, respectively,
$\Phi_\pm$ analytic in $\Lambda$ as functions from $\Omega$
to $C^0 [0,\pm\infty)$, such that: 
\medbreak
(i)
For any fixed $0<\btheta<\theta$ and $0\le k\le K+1$,  $j\ge 0$, 
\begin{equation}
\label{Pdecay} 
|(\partial/\partial \Lambda)^j(\partial/\partial x)^k
\Theta_\pm |\le C(j) C_1 C_2 e^{-\theta |x|/C_2}
\quad
\text{\rm for } x\gtrless 0. 
\end{equation}
%
\smallbreak
(ii)  The change of coordinates $\WW=:P_\pm \ZZ$,
$\FF=: P_\pm \GG$ reduces \eqref{gfirstorder} to 
\begin{equation}
\label{glimit}.
\ZZ'-\mA_\pm \ZZ = \GG
\quad
\text{\rm for } x\gtrless 0.
\end{equation}
\end{cor}

Equivalently, solutions of \eqref{gfirstorder} may be factored as 
\begin{equation}
\label{Wfactor}
\WW=(I+ \Theta_\pm)\ZZ_\pm, 
\end{equation}
where $\ZZ_\pm$ satisfy the limiting, constant-coefficient 
equations \eqref{glimit} and $\Theta_\pm$ satisfy bounds \eqref{Pdecay}. 

\begin{proof}
Substituting $\WW=P_- Z$ into \eqref{gfirstorder}, equating
to \eqref{glimit}, and rearranging, we obtain the defining equation
\begin{equation}
\label{matrixODE}
P_-'= \mA_-P_- - P_- \mA, \qquad P_- \to I \quad \text{\rm as}
\quad
x\to -\infty.
\end{equation}
Viewed as a vector equation, this has the form
$ P_-'=\CalA P_-, $
where $\CalA$ approaches exponentially as $x\to -\infty$ to its
limit $\CalA_-$, defined by
\begin{equation}
\label{CalA}
\CalA_- P:= 
\mA_-P- P \mA_-.
\end{equation}
The limiting operator $\CalA_-$ evidently has 
analytic eigenvalue, eigenvector pair $\mu\equiv 0$, $P_-\equiv I$,
whence the result follows by Lemma \ref{gaplemma} for $j=k=0$.
The $x$-derivative bounds $0<k\le K+1$ then follow from the ODE and its first
$K$ derivatives,
and the $\Lambda$-derivative bounds from 
standard interior estimates for analytic functions.
A symmetric argument gives the result for $P_+$.
\end{proof}

\begin{rem}
\label{frozencoeffs}
\textup{
Equation \eqref{Wfactor} gives an
explicit connection to the inviscid, bi-constant-coefficient case,
for bounded frequencies $|(\tilde \xi, \lambda)|\le R$ 
(low frequency $\sim$ inviscid regime).
For high-frequencies, the proper analogy is rather to the 
frozen-coefficient
case of local existence theory.
}
\end{rem}

\begin{ex}[{[B]}] 
\label{doubleduhamel}
Use Duhamel's formula to show that
\begin{equation}
\label{Peq}
P V(x) = V^- + \int^x_{-\infty} e^{(\mA_- - \mu I)(x-y)} 
P \theta (y) V(y) dy
\end{equation}
and
\begin{equation}
\label{Qeq}
QV(x)= QV(-M) + \int^{x}_{-M} e^{(\mA_- - \mu I)(x-y)} Q \theta (y) V(y) dy,
\end{equation}
hence $V(x)=\CalT V(x)$ for the unique solution $V$ of \eqref{3.7g}
determined by conditions $V(-\infty)=V_-=PV_-$ and $QV(-M)=0$.
\end{ex}

\begin{ex}\label{jordanex}
(i) If $\{r_j\}$ and $\{l_k\}$ are dual bases of (possibly generalized)
right and left eigenvectors of $\mA_-$,
with associated eigenvalues $\mu_j$ and $\mu_k$, show that
$\{r_jl_k^*\}$ is a basis of (possibly generalized) 
right eigenvectors of the operator $\CalA_-$ defined in \eqref{CalA},
with associated eigenvalues $\mu_j-\mu_k$.
(ii) Show that $\mu=0$ is a semisimple eigenvalue of $\CalA_-$ 
if and only if each eigenvalue of $\mA_-$ is semisimple.
\end{ex}

\medbreak
{\bf 3.4. Hyperbolic--parabolic smoothing.}
We next introduce the circle of ideas associated with 
hyperbolic--parabolic smoothing and estimate \eqref{Kenergy}.

\begin{lem} [{[Hu, MaZ.5]}]
\label{skewlem}
Let $A=\bdiag\{a_jI_{m_j}\}$ be diagonal, with real entries $a_j$ appearing 
with prescribed multiplicities $m_j$ in order of increasing size, 
and let $B$ be arbitrary.
Then, there exists a smooth skew-symmetric matrix-valued function 
$K(A,B)$ such that
\begin{equation}
\label{skew}
\text{\rm Re }\left( B- KA  \right)=
\R \bdiag B,
\end{equation}
where $\bdiag B$ 
denotes the block-diagonal part of $B$, with blocks of dimension
$m_j$ equal to the multiplicity of the corresponding eigenvalues of $A$.
\end{lem}
\medbreak

\begin{proof}
It is straightforward to check that the symmetric matrix
$\R KA=(1/2)(KA-A^tK)$
may be prescribed arbitrarily on off-diagonal blocks,
by setting $K_{ij}:= (a_i-a_j)^{-1}M_{ij}$, where $M_{ij}$
is the desired block, $i\ne j$.
Choosing $M=\R B$, we obtain
$\R(B- KA )=\R \bdiag(B)$ as claimed.
\end{proof}
\medbreak

\begin{lem}[{[KSh]}]
\label{KSh} 
Let $\tilde A^0$, $A$, and $B$ denote real-valued matrices such that
$\tilde A^0$ is symmetric positive definite and $\tilde A:=\tilde A^0 A$ and
$\tilde B:=\tilde A^0 B$ are symmetric, $\tilde B\ge 0$.
Then, the following are equivalent:

\medbreak
(K0) (Genuine coupling)
No eigenvector of $A$ lies in $\ker B$ (equivalently,
in $\ker \tilde B$).

\smallbreak
(K1) $\bdiag L BR= R^t \tilde B R >0$, 
where $L:= \tilde O^t (\tilde A^0)^{1/2}$ and
$R:=(\tilde A^0)^{-1/2}\tilde O$
are matrices of left and right eigenvectors of $A$
block-diagonalizing $LAR$, with $O$ orthonormal.
Here, as in Lemma \ref{skewlem}, 
$\bdiag M$ denotes the matrix formed from the diagonal
blocks of $M$, with blocks of dimension equal to the
multiplicity of corresponding eigenvalues of $LAR$.

\medbreak

(K2) (hyperbolic compensation) 
There exists a smooth skew-symmetric matrix-valued function
$K(\tilde A,\tilde B,\tilde A^0)$ such that
\begin{equation}
\label{skew0}
\text{\rm Re }\left(\tilde B- K A \right) > 0.
\end{equation}
\medbreak
(K3) (Strict dissipativity)
\begin{equation}
\label{diss}
\R\sigma(-i\xi A
- |\xi|^2 B)
\le -\theta
|\xi|^2/(1+|\xi|^2), \quad \theta>0.
\end{equation}
\end{lem}

\begin{proof}
(K0) $\Leftrightarrow$ (K1) by the property of symmetric
nonnegative matrices $M$ that $v^tMv=0$ if and only if
$Mv=0$ (exercise), so that $\bdiag\{R^t_j \R \tilde B R_j\}$
has a kernel if and only if $\alpha^t R_j^t \R \tilde B R_j\alpha=0$,
if and only if $\R \tilde B R_j \alpha=0$, where $R_j$ denotes
a block of eigenvectors with common eigenvalue.

(K1) $\Rightarrow$ (K2) follows 
readily from Lemma \ref{KSh},
by first converting to the case of symmetric $A$ and $B\le 0$
by the transformations $A\to (\tilde A^0)^{1/2}A\tilde A^0)^{-1/2}$,
$B\to (\tilde A^0)^{1/2}B\tilde A^0)^{-1/2}$, from which the original result
follows by the fact that $M>0\Leftrightarrow (\tilde A^0)^{1/2}M(\tilde A^0)^{1/2}>0$,
then converting by an orthonormal change of coordinates to the
case that $A$ is diagonal and $B\le 0$. 
Variable multiplicity eigenvalues may be handled by partition of
unity/interpolation, noting that $\R(B-KA)<0$ persists under perturbation.

(K2) $\Rightarrow$ (K3) follows upon rearrangement of energy estimate
\begin{equation}
\label{FTenergy}
\begin{aligned}
0&=\R\langle ((C+|\xi|^2)\tilde A^0 + i\xi K)w, 
(\lambda  + i\xi A + |\xi|^2 B)w\rangle\\
&=
\R \lambda \langle w, 
\Big((C+|\xi|^2) \tilde A^0 +i\xi K \Big)w \rangle
+ |\xi|^2  \langle w, \R(\tilde B- KA)w\rangle\\
&\quad +
\R \langle w, -i|\xi|^3 K(\tilde A^0)^{-1}\tilde Bw\rangle 
+ C|\xi|^4\langle  w, \tilde B w \rangle,\\
\end{aligned}
\end{equation}
which yields
\begin{equation}
\begin{aligned}
\R \lambda \langle w,& 
\Big((1+C|\xi|^2) \tilde A^0 +i\xi K \Big)w \rangle\\
&\le - \Big( \theta|\xi|^2|w|^2 + M|\xi|^3(|w||\tilde B w| + 
(C/M)|\xi|^4|\tilde B w|^2\Big)\\
&\le 
- \theta|\xi|^2|w|^2
\end{aligned}
\end{equation}
and thereby
\begin{equation}
\label{finen}
 \R \lambda (1+ |\xi|^2)|w|^2 \le 
-\theta_1 |\xi|^2|w|^2,
\end{equation}
for $M$, $C>0$ sufficiently large and $\theta$, $\theta_1>0$
sufficiently small, by positivity of
$\tilde A^0$ and $\R(\tilde B-KA)$.
%

Finally, (K3) $\Rightarrow$ (K1) follows by
by first-order Taylor expansion at $\xi=0$ of the spectrum of 
$ LAR -i\xi LBR$ (well-defined, by symmetry of $LAR$), together
with symmetry of $R^t\tilde B R$.
\end{proof}
\medbreak

\begin{cor} Under (A1)--(A2), there holds the uniform dissipativity
condition
\begin{equation}
\label{multidiss}
\R \sigma( \sum_j i\xi_j A^j - \sum_{j,k} \xi_j \xi_k B^{jk})_\pm\le 
-\theta |\xi|^2/(1+ |\xi|^2).
\end{equation}
Moreover, there exist smooth skew-symmetric 
``compensating matrices'' $K_\pm(\xi)$, homogeneous degree
one in $\xi$, such that 
\begin{equation}
\label{Krel}
\R \Big( 
\sum_{j,k} \xi_j \xi_k
\tilde B^{jk} - K(\xi) (\tilde A^0)^{-1}
\sum_k \xi_k \tilde A^k \Big)_\pm 
\ge \theta >0
\end{equation}
for all $\xi \in \RR^d \setminus \{0\}$.
\end{cor}

\begin{proof} 
By the block--diagonal structure of $\tilde B^{jk}$
(GC) holds also for $A^j_\pm$ and 
$\hat B^{jk}:=(\tilde A^0)^{-1}\R \tilde B^{jk}$, since 
\begin{equation}
\ker \sum_{j,k}\xi_j\xi_k\hat B^{jk}=\ker \sum_{j,k}\xi_j\xi_k\R \tilde B^{jk}
=\ker \sum\xi_j\xi_k \tilde B^{jk}
=\ker \sum\xi_j\xi_k B^{jk}.
\end{equation}
Applying Lemma \ref{KSh} to 
\begin{equation}
\tilde A^0:=\tilde A^0_\pm, 
\quad A:= \Big((\tilde A^0)^{-1} \sum_k \xi_k\tilde A^k \Big)_\pm,
\quad 
B:= \Big( (\tilde A^0)^{-1} \sum_{j,k} \xi_j \xi_k \R\tilde B^{jk} \Big)_\pm,
\end{equation}
we thus obtain \eqref{Krel} and
\begin{equation}
\label{intdis}
\R \sigma \Big[(\tilde A^0)^{-1} \big(-\sum_j i\xi_j \tilde A^j 
- \sum_{j,k} \xi_j\xi_k \R \tilde B^{jk}\big)\Big]_\pm
\le -\theta_1
|\xi|^2/(1+|\xi|^2), 
\end{equation}
$\theta_1>0$,
from which we readily obtain
\begin{equation}
\big(-\sum_j i\xi_j \tilde A^j 
- \sum_{j,k} \xi_j\xi_k \tilde B^{jk}\big)_\pm
\le -\theta_2
|\xi|^2/(1+|\xi|^2)
\end{equation}
and thus \eqref{multidiss} (exercise,
using $M>\theta_1 \Leftrightarrow
(\tilde A^0)_\pm^{-1/2}M(\tilde A^0)^{-1/2}_\pm >\theta$
and
$\sigma (\tilde A^0)_\pm^{-1/2}M(\tilde A^0)^{-1/2}_\pm >\theta
\Leftrightarrow \sigma (\tilde A^0)_\pm^{-1}M>\theta$,
together with $S>\theta \Leftrightarrow \sigma S>\theta$
for $S$ symmetric).
Because all terms other than $K$ in the lefthand side of \eqref{Krel} 
are homogeneous, it is 
evident that we may choose $K(\cdot)$ homogeneous as well
(restrict to the unit sphere, then take homogeneous extension).
\end{proof}

{\bf 3.4.1. Basic estimate.}
Energy estimate \eqref{FTenergy}--\eqref{finen} may be recognized
as the Laplace--Fourier
transformed version 
of a corresponding time-evolutionary estimate
\begin{equation}
\label{kaw}
\begin{aligned}
(d/dt)&(1/2)\Big(\langle \tilde A^0 W,W\rangle + 
\langle \tilde A^0 W_x,W_x\rangle + 
\langle K\partial_x W,W\rangle \Big)\\
&= -\langle W_x, \R(\tilde B- KA)W_x\rangle 
- \langle W_x, K(\tilde A^0)^{-1}\tilde BW_x\rangle 
- C\langle  W_{xx}, \tilde B W_{xx} \rangle,\\
&\le -\theta(|W_x|_{L^2}^2 + |\tilde BW_{x}|_{H^1}^2)
\end{aligned}
\end{equation}
for the one-dimensional, linear constant-coefficient equation
\begin{equation}
\tilde A^0 W_t + \tilde AW_x=\tilde BW_{xx},
\end{equation}
which in the block-diagonal case $\tilde B=\bdiag\{0,\tilde b\}$,
$\tilde b>0$, yields \eqref{Kenergy} for $s=1$
by the observation that 
\begin{equation}
{\cal{E}}(W):=(1/2)\Big(\langle \tilde A^0 W,W\rangle 
+\langle K\partial_x W,W\rangle \Big)
+C\langle \tilde A^0 W_x,W_x\rangle
\end{equation}
for $C>0$ sufficiently large determines a norm 
${\cal{E}}^{1/2}(\cdot)$ equivalent to $|\cdot|_{H^1}$.
This readily generalizes to 
\begin{equation}
\label{Eexp}
(d/dt){\cal{E}}(t)\le -\theta ( |\partial_x W|_{H^{s-1}}^2
+ |\tilde B \partial_x W|_{H^s}^2), 
\end{equation}
\begin{equation}
\begin{aligned}
\label{genE}
{\cal{E}}(W)&:=(1/2)\Big(\langle \tilde A^0 W,W\rangle 
+\langle K \partial_x W,W\rangle \Big)\\
&\quad
+C\Big(\langle \tilde A^0 W_x,W_x\rangle 
+ \langle K\partial_x W_{x},W_{x}\rangle \Big) + \cdots + \\
&\quad +C^{s-1}\Big(
\langle \tilde A^0 d_x^{s-1}W,d_x^{s-1}W\rangle
+ \langle K\partial_x \partial_x^{s-1} W,\partial_x^{s-1}W\rangle \Big)\\
&\quad + C^{s} \langle \tilde A^0 d_x^{s}W,d_x^{s}W\rangle,
\end{aligned}
\end{equation}
where ${\cal{E}}^{1/2}(\cdot)$ is a norm equivalent to $|\cdot|_{H^s}$,
yielding \eqref{Kenergy} for $s\ge 1$.
We refer to energy estimates of the general type (K3), 
\eqref{kaw}--\eqref{genE} as ``Kawashima-type'' estimates.

\begin{prop}
\label{kawtype}
Assuming (A1)--(A2), (H0) for $W_-=0$, energy estimate \eqref{Kenergy}
is valid for all $q(d)\le s\le q$, $q(d)$ and $q$ as defined
in (H0),
so long as $|W|_{H^s}$ remains sufficiently small.
\end{prop}

\begin{proof}
In the linear, constant-coefficient case, \eqref{Krel}
together with a calculation analogous to that of
\eqref{kaw} and \eqref{Eexp} yields
\begin{equation}
\label{multiEexp1}
(d/dt){\cal{E}}(t)\le -\theta ( |\partial_x W|_{H^{s-1}}^2
+ |\partial_x w^{II}|_{H^s}^2)
\end{equation}
for
\begin{equation}
\begin{aligned}
\label{multigenE}
{\cal{E}}(W)&:=(1/2)\Big(\langle \tilde A^0_- W,W\rangle 
+\langle K_-(\partial_x) W,W\rangle \Big)\\
&\quad
+C\Big(\langle \tilde A^0_- \partial_x W,\partial_x W\rangle 
+ \langle K_-(\partial_x) \partial_x W,\partial_x W\rangle \Big)
+\cdots +\\
&\quad +C^{s-1}\Big(
\langle \tilde A^0_- \partial_x^{s-1}W,\partial_x^{s-1}W\rangle
+ \langle K_-(\partial_x) \partial_x^{s-1} W,\partial_x^{s-1} W\rangle \Big)\\
&\quad
+ C^{s} \langle \tilde A^0_- \partial_x^{s}W,\partial_x^{s}W\rangle,
\end{aligned}
\end{equation}
where operator $K_-(\partial_x)$ is defined by
\begin{equation}
\widehat{K_-(\partial x f)}(\xi):= iK_-(\xi) \hat f(\xi)
\end{equation}
with $K_-(\cdot)$ as in \eqref{Krel}, where $\hat g$ denotes
Fourier transform of $g$, for all $s\ge 1$.

In the general (nonlinear, variable-coefficient) case, 
we obtain by a similar calculation
\begin{equation}
\label{multiEexp}
(d/dt){\cal{E}}(t)\le -\theta ( |\partial_x W|_{H^{s-1}}^2 
+ | \partial_x w^{II}|_{H^s}^2)
+ C|W|_{L^2}^2
\end{equation}
for $C>0$ sufficiently large, for
\begin{equation}
\begin{aligned}
\label{multigenEvc}
{\cal{E}}(W)&:=(1/2)\Big(\langle \tilde A^0(W) W,W\rangle 
+\langle K_-(\partial_x) W,W\rangle \Big)\\
&\quad +C\Big(\langle \tilde A^0(W) \partial_x W,\partial_x W\rangle 
+ \langle K_-(\partial_x) \partial_x W,\partial_x W\rangle \Big)
+\cdots +\\
&\quad 
+C^{s-1}\Big( \langle \tilde A^0(W) \partial_x^{s-1}W,\partial_x^{s-1}W\rangle
+ \langle K_-(\partial_x) \partial_x^{s-1} W,\partial_x^{s-1} W\rangle \Big)\\
&\quad
+ C^{s} \langle \tilde A^0(W) \partial_x^{s}W,\partial_x^{s}W\rangle,
\end{aligned}
\end{equation}
provided $s\ge q(d):= [d/2]+2$ and the coefficient functions
possess sufficient regularity $C^s$, i.e., $s\le q$.

Here, we estimate time-derivatives of $\tilde A^0$ terms 
in straightforward fashion, using integration by parts
as in the proof of Proposition \ref{hypwellposed}.
We estimate terms involving the nondifferential operator
$K_-(\partial_x)$ in the frequency domain as
\begin{equation}
\label{Keq}
\begin{aligned}
(d/dt)(1/2)\langle K(\partial_x)\partial_x^r W, \partial_x^r W\rangle&=
(d/dt)(1/2)\langle iK(\xi) (i\xi)^r \hat W, (i\xi)^r \hat W\rangle\\
&= \langle iK(\xi) (i\xi)^r \hat W, (i\xi)^r \hat W_t\rangle\\
&= \langle (i\xi)^r \hat W, -K(\xi)(\tilde A^0_-)^{-1}\big(\sum_j \xi_j
\tilde A^j_- \big)\hat (i\xi)^r \hat W\rangle \\
&\quad + \langle (iK(\xi) i\xi)^r \hat W, (i\xi)^r \hat H\rangle\\
\end{aligned}
\end{equation}
using Plancherel's identity together with the equation, written in the
frequency domain as
\begin{equation}
\hat W_t = -\sum_j i\xi_j (\tilde A^0_-)^{-1}\tilde A^j_- \hat W + \hat H,
\end{equation}
where 
\begin{equation}
\label{H}
\begin{aligned}
H&:= \sum_j
\big((\tilde A^0_-)^{-1}\tilde A_-^j- (\tilde A^0)^{-1} \tilde A^j(W))W_{x_j}\big)\\
& +\sum_{j,k} (\tilde A^0)^{-1} (\tilde B^{jk}W_{x_k})_{x_j} 
+ (\tilde A^0)^{-1} \tilde G(W_x, W_x)\Big).\\
\end{aligned}
\end{equation}

By a calculation similar to those in the proofs of Propositions
\ref{hypwellposed} and \ref{parabwellposed}
(exercise, using smallness of $ |\tilde A_-^j- \tilde A^j(W)| \sim |W|$
and the Moser inequality \eqref{moser}),
we obtain
\begin{equation}
\label{Hbd}
| \partial_x^{r} H|_{L^2}\le
C|\partial_x^{r+2}w^{II}|_{L^{2}}
+ 
C|W|_{H^{r+1}}\big(|W|_{H^{[d/2]+2}}+
|W|_{H^{[d/2]+2}}^r\big).
\end{equation}
Thus, using homogeneity, $|K(\xi)|\le C|\xi|$, 
together with the Cauchy--Schwartz
inequality and Plancherel's identity,
we may estimate the final
term on the righthand side of \eqref{Keq} as 
\begin{equation}
\label{Keqfinal}
\begin{aligned}
\langle (iK(\xi) i\xi)^r \hat W, &(i\xi)^r \hat H\rangle
\le C|\partial_x^{r+1}  W|_{L^2} || \partial_x^{r} H|_{L^2}\\
&\le 
 C|\partial_x^{r+1} W|_{L^2} |\partial_x^{r+2}w^{II}|_{L^{2}}
+ \epsilon |W|_{H^{r+1}}^2,
\end{aligned}
\end{equation}
any $\epsilon>0$,
for $|W|_{H^{[d/2]+2}}$
sufficiently small: that is, 
a term of the same form arising in the constant-coefficient case
plus an absorbable error.

Combining $\tilde A^0$- and $K(\xi)_-$-term estimates, 
and using \eqref{Krel} together with
\begin{equation}
\label{calcmisc}
\begin{aligned}
\sum_{j,k} \langle \partial_x^r W_{x_j}, \tilde B^{jk}(W) 
\partial_x^r W_{x_k}\rangle&=
\langle  (i\xi)^r \hat W,  \sum_{j,k}\xi_j\xi_k \tilde B^{jk}_- 
(i\xi)^r \hat W\rangle\\
&\quad + {\CalO}(|\partial_x^{r+1} W|_{L^2}^2|W|_{L^\infty})\\
&\ge
\langle  (i\xi)^r \hat W,  \sum_{j,k}\xi_j\xi_k \tilde B^{jk}_- 
(i\xi)^r \hat W\rangle
-\epsilon|\partial_x^{r+1} W|_{L^2}^2
\end{aligned}
\end{equation}
for any $\epsilon>0$, for $|W|_{L^\infty}$ sufficiently small,
we obtain the result, similarly as in the constant-coefficient case. 
\end{proof}

\begin{ex}
\label{interp}
(Alternative proof)
Show using \eqref{vsfriedrichs} 
together with \eqref{Keq}--\eqref{calcmisc} that
\begin{equation}
\label{Ecalc}
(d/dt)\CalE\le  C|W|_{H^{s-1}}^2 - 
\theta_1 (|\partial_x^s W|_{L^2}^2 + |\partial_x^{s+1}w^{II}|_{L^2})
\end{equation}
for $C>0$ sufficiently large and $\theta_1>0$ sufficiently small, for
\begin{equation}
\CalE(W):=  
(1/2)\sum_{r=0}^{s}\langle \partial_x^r W, \tilde A^0 \partial_x^r W\rangle
+ (1/2C) \langle K_-(\partial_x) \partial_x^{s-1} W,\partial_x^{s-1}
W\rangle.
\end{equation}
Using the $H^s$ interpolation formula
\begin{equation}
\label{Hsinterp}
|f|_{H^{s_*}}^2\le \beta C^{1/\beta}|f|_{H^{s_1}}^2+ 
(1-\beta)C^{-1/(1-\beta)}|f|_{H^{s_2}}^2,
\quad
\beta= (s_2-{s_*})/(s_2-s_1), 
\end{equation}
valid for $s_1\le s_* \le s_2$,
taking $s_1=0$, $s_2=s$, and $C$ sufficiently large,
show that \eqref{Ecalc} implies \eqref{multiEexp}.
\end{ex}

{\bf 3.4.2.  Alternative formulation.}
A basic but apparently new observation is that energy estimate
\eqref{multigenEvc} implies not only \eqref{Kenergy} but
also the following considerably stronger estimate.

\begin{prop}\label{ccKenergy}
Under the assumptions of Proposition \ref{kawtype},
\begin{equation}
\label{ccexp}
|W(t)|^2_{H^s}\le C\big(e^{-\theta_1 t}|W(0)|^2_{H^s}
+\int_0^t e^{-\theta_1 (t-s)}|W(s)|_{L^2}^2\, ds\big),
\qquad
\theta_1>0, 
\end{equation}
so long as $|W|_{H^s}$ remains sufficiently small.
\end{prop}

\begin{proof}
Rewriting \eqref{multiEexp} using equivalence of
$\CalE^{1/2}$ and $H^s$ as
\begin{equation}
(d/dt){\cal{E}}(t)\le -\theta_1 {\cal{E}}(t) + C|W(t)|_{L^2}^2,
\end{equation}
we obtain by Gronwall's inequality
$e^{\theta_1 t}\calE |^t_0\le C \int_0^t e^{\theta_1 s}|W(s)|_{L^2}^2\, ds$
and thereby the result (exercise).
\end{proof}

\begin{proof}[Proof of Proposition \ref{globalwellposed}]
Exercise, using \eqref{ccexp} together with 
the fact that $|U(t)|_{L^2}\le C|U(0)|_{L^2}$ 
for $|W|_{L^\infty}$ sufficiently small (Exercise \ref{gposed}).
\end{proof}

In Appendix B, we sketch also a simple proof of Proposition
\ref{conststab} using \eqref{ccexp} together with linearized
stability estimates.
This may be helpful for the reader in motivating the nonlinear
stability argument carried out in Section 4 for the more complicated
variable-coefficient situation of a viscous shock profile.

\bigbreak
\clearpage

\section{Reduction to low frequency}

We now begin our main analysis,
carrying out high- and mid-frequency estimates
reducing
the nonlinear long-time stability problem (LT) 
to the study of the linear resolvent equation $(\lambda -L_\Txi)U=f$
in the low-frequency regime $|(\Txi,\lambda)|$ small
described in Proposition \ref{ZS}.
The novelty of the present treatment lies in the simplified argument
structure based entirely on energy estimates.
More detailed linearized estimates were obtained in [Z.4] using 
energy estimates together with the pointwise machinery of [MaZ.3].

\medbreak

{\bf 4.1. Nonlinear estimate.}
Define the nonlinear perturbation
\begin{equation}
\label{pert1}
U:= \tilde U - \bar U,
\end{equation}
where $\tilde U$ denotes a solution of \eqref{viscous} with
initial data $\tilde U_0$ close to $\bar U$.
Our first step is to establish the following large-variation
version of Proposition \ref{ccKenergy}.

\begin{prop}[{[MaZ.4, Z.4]}]
\label{auxenergy}
Given  (A1)--(A2), (H0)--(H3), and Assumption \ref{lax},
\begin{equation}
\label{auxexp}
|U(t)|^2_{H^s}\le C\big(e^{-\theta_2 t}|U(0)|^2_{H^s}
+\int_0^t e^{-\theta_2 (t-s)}|U(s)|_{L^2}^2\, ds\big),
\quad \theta_2>0,
\end{equation}
so long as $|U|_{H^s}$ remains sufficiently small, for
all $q(d)\le s\le q$, $q(d)$ and $q$ as defined in (H0).
\end{prop}

\medbreak

{\bf 4.1.1. ``Goodman-type'' weighted energy estimate.}
The obvious difficulty in proving Proposition \ref{auxenergy}
is that (A2) is assumed only at endpoints $W_\pm$, hence we
cannot hope to establish smoothing in the key  hyperbolic 
modes $w^I$ through the circle of ideas discussed in Section 3.4,
except in the far field $x_1 \to \pm\infty$.
The complementary idea needed to treat the ``near field'' or
``interior'' region $x_1\in [-M,M]$ is that propagation in
the hyperbolic modes thanks to assumption (H1) is
uniformly transverse to the shock profile, so that signals 
essentially spend only finite time in the near field,
at all other times experiencing smoothing properties of the far field.

The above observation may be conveniently quantified using a type
of weighted-norm energy estimate introduced by Goodman [Go.1--2]
in the study of one-dimensional stability of small-amplitude
shock waves.
Consider a linear hyperbolic equation 
$\tilde A^0W_t + \sum_{j=1}^d \tilde A^j W_{x_j}=0$ 
with large- but finite-variation coefficients depending
only on $x_1$,
\begin{equation}
\label{coeffvar}
|A^j_{x_1}|\le \Theta, \quad \int \Theta(y)\, dy < \infty,
\quad j=0, 1, \dots d,
\end{equation}
$\tilde A^j$ symmetric, $\tilde A^0>0$,
under the ``upwind'' assumption 
\begin{equation}
\label{upwind}
\tilde A^1\ge \theta>0.
\end{equation}
Define the scalar weight
$\alpha(x_1):= e^{\int_0^{x_1} -(2C\Theta/\theta)(y)\, dy}$,
positive and bounded above and below by \eqref{coeffvar},
where $C>0$ is sufficiently large.

Then, the zero-order Goodman's estimate, in this simple setting, is just
\begin{equation}
\label{goodmaneg}
\begin{aligned}
(d/dt)(1/2)\langle W, \alpha \tilde A^0 W\rangle&=
\langle W, \alpha\tilde A^0 W_t \rangle\\
&=
-\langle W, \alpha \sum_j \tilde A^j W_{x_j} \rangle\\
&=
(1/2)\langle W, (\alpha \tilde A^1)_{x_1} W \rangle\\
&=
(1/2)\langle W, \alpha (-(2C\Theta/\theta) \tilde A^1
+ \tilde A^1_{x_1}) W \rangle\\
&\le
-(C/2)\langle W, \alpha \Theta W\rangle, 
\end{aligned}
\end{equation}
yielding time-exponential decay in $|W|_{L^2}$ on any 
finite interval $x_1\in [-M,M]$.

\begin{ex}
\label{goodman}
Verify the corresponding $s$-order estimate
\begin{equation}
\label{sorder}
(d/dt)(1/2)\sum_{r=0}^s
\langle \partial_x^r W, \alpha \tilde A^0 \partial_x^r W\rangle
\le
-(C/2)\sum_{r=0}^s 
\langle \partial_x^r W, \alpha \Theta \partial_x^r W\rangle
\end{equation}
for $C>0$ sufficiently large, yielding time-exponential decay
in $|W|_{H^s}$ on any finite interval $x_1\in [-M,M]$.
\end{ex}

\begin{rem}
\label{trans}
\textup{
Proper accounting of the favorable effects of transverse propagation 
is a recurring theme in the analysis of stability of 
viscous waves; see, e.g., [Go.1--2, L.4, LZ.1--2, LZe.1, SzZ, 
L.1, ZH, BiB, GrR, M\'eZ.1, MGWZ.1--4].
}
\end{rem}

\medbreak
{\bf 4.1.2. Large-amplitude hyperbolic--parabolic smoothing.}
Combining the weighted energy estimate just described above with
the techniques of Section 3.4, we are now ready to 
establish Proposition \ref{auxenergy}.

\begin{proof}[Proof of Proposition \ref{auxenergy}]
Equivalently, we establish
\begin{equation}
\label{Wauxexp}
|W(t)|^2_{H^s}\le C\big(e^{-\theta_1 t}|W(0)|^2_{H^s}
+\int_0^t e^{-\theta_1 (t-s)}|W(s)|_{L^2}^2\, ds\big),
\quad \theta_1>0,
\end{equation}
where
\begin{equation}
\label{W}
W(x,t):= \tilde W(x,t)-\bar W(x_1):= W(\tilde U(x,t))-W(\bar U(x_1)).
\end{equation}
(Exercise: Show $|U|_{H^r}\sim |W|_{H^r}$ for all
$0\le r\le s$, in particular $0$ and $s$, using
$s\ge [d/2]+2$ and same estimates used
to close the energy estimate in the proof of Proposition
\ref{parabwellposed}.)
%

Noting that both $\tilde W$ and $\bar W$ satisfy \eqref{symmhypparab},
we obtain by straightforward calculation [MaZ.4, Z.4] 
the nonlinear perturbation equation
\begin{equation}
\label{Wpert}
\begin{aligned}
\tilde A^0 W_t + &\sum_j \tilde A^j W_{x_j}= 
\sum_{j,k}(\tilde B^{jk}W_{x_k})_{x_j}\\
&\quad+
M_1\bar W_{x_1} + \sum_j(M_2^j \bar W_{x_1})_{x_j}
+\sum_j M_3^j W_{x_j} 
+ G(x, \partial_x W),
\end{aligned}
\end{equation}
where
\begin{equation}
\label{Wcoeffs}
\tilde A^j(x,t):= \tilde A^j(\tilde W(x,t)), 
\quad \tilde B^{jk}(x,t):= \tilde B^{jk}(\tilde W(x,t)),
\end{equation}
\begin{equation}
\label{M1}
M_1=M_1(W, \Bar W):=\tilde A^1(\tilde W)-\tilde A^1(\bar W)
=
\left( \int_0^1 d\tilde A^1(\bar W + \theta W )\, d\theta  \right) W
\end{equation}
\begin{equation}
\label{M2}
M_2^j=M_2^j(W, \Bar W):= \tilde B^{j1}-\bar B^{j1}
=
\begin{pmatrix}
0&0\\
0& 
( \int_0^1 db^{j1}(\bar W + \theta W )\, d\theta  ) W
\end{pmatrix}.
\end{equation}
\begin{equation}
\label{M3}
M_3^j=M_3^j(x_1):= 
(\partial \tilde G/\partial W_{x_j})|_{\partial_x \bar W} 
=
\begin{pmatrix}
0&0\\
\CalO(|\bar W_{x_1}|) &\CalO(|\bar W_{x_1}|) 
\end{pmatrix}.
\end{equation}
and 
\begin{equation}
\label{WG}
\begin{aligned}
G= \begin{pmatrix} 0\\g \end{pmatrix} 
&:= \begin{pmatrix} 0 \\ g(\partial_x \tilde W, \partial_x \tilde W)-
g(\partial_x \bar W, \partial_x \bar W)- dg(\bar W)\partial_x W \end{pmatrix}\\
&=\begin{pmatrix} 0 \\ \CalO(|\partial_x W|^2)\end{pmatrix}.
\end{aligned}
\end{equation}
(In the case that there exists a global convex viscosity-compatible
entropy, we may further arrange that $G\equiv 0$.)

To clarify the argument, we drop terms $M$ and $G$ in \eqref{Wpert}.
The reader may verify that these generate harmless,
absorbable error terms in our energy estimates
(exercise; see [MaZ.4, Z.4]).  
With this simplification, we reduce to the familiar situation
\begin{equation}
\label{redpert}
\tilde A^0 W_t + \sum_j \tilde A^j W_{x_j}= 
\sum_{j,k}(\tilde B^{jk}W_{x_k})_{x_j},
\end{equation}
with the difference that 
$\tilde A^j=\tilde A^j(\bar W+W)$ and
$\tilde B^{jk}=\tilde B^{jk}(\bar W+W)$ 
are no longer approximately constant, but for small $|W|_{H^s}$
approach the possibly rapidly varying values
$\tilde A^j(\bar W)$ and $\tilde B^{jk}(\bar W)$ along
the background profile.

By Corollary \ref{odecor} and Assumption (H1), respectively, we have
\begin{equation}
|(d/d_{x_1})^k \tilde A^j(\bar W)|, \,
|(d/d_{x_1})^k \tilde B^{jk}(\bar W)| \le
\Theta(x_1):= Ce^{-\theta|x|} 
\end{equation}
and
\begin{equation}
\label{upwind2}
A^1\ge \theta>0,
\end{equation}
for some $C$, $\theta>0$.
Define
\begin{equation}
\label{alpha}
\alpha(x_1):= e^{\int_0^{x_1} -(2C_*\Theta/\theta)(y)\, dy},
\end{equation}
where $C_*>0$ is a sufficiently large constant to be determined later.

Decompose $W$ into the sum of near- and far-field parts $W_N$ and $W_F^\pm$,
\begin{equation}
\label{NFdecomp}
\begin{aligned}
W&= W_{N}+ W_{F}^++ W_F^-\\
&:= \Big(1-\chi^+(x_1)-\chi^-(x_1)\Big) W
+ \chi^+(x_1) W + \chi^-(x_1)W,
\end{aligned}
\end{equation}
where $\chi^\pm$ are $C^\infty$ cutoff functions 
supported on $x_1\ge M$ and $x_1\le -M$ 
and one on $x_1\ge 2M$ and $x_1\le -2M$, respectively,
for $M>0$ sufficiently large, with
\begin{equation}
\label{chibds}
|d^r \chi^\pm(x_1) |\le \epsilon 
\min \{|\chi^\pm|, |1-\chi^\pm|\}, 
\qquad 0\le r\le s+1,
\end{equation}
where $\epsilon>0$ is a sufficiently small constant to be determined later.
Up to lower-derivative error terms that are harmless for $\epsilon$
sufficiently small, both $W_N$ and $W_F^\pm$ satisfy the same equation
\eqref{redpert} satisfied by $W$.
Again, we'll drop these terms, leaving it to the reader to
verify that they can indeed be absorbed in the estimates.
(The main point is to notice that first-order contributions
$\sum_{k} d\chi \tilde B^{1k} W_{x_k}$ have zero first row,
hence do not affect the structural assumption $\tilde A^j_{11}$
symmetric.)

On the respective supports $x_1\ge M$, $x_1\le -M$ of $W_F^\pm$, 
coefficients $\tilde A^j$
and $\tilde B^{jk}$ remain arbitrarily close in $H^s$ to $\tilde A^j_\pm$
and $\tilde B^{jk}_\pm$ 
provided $|W|_{H^s}$ is taken sufficiently small and $M$ sufficiently large.
Thus, using the argument of Proposition \ref{kawtype}, we obtain
\begin{equation}
\label{EFexp}
(d/dt){\cal{E}_F^\pm}(W_F^\pm(t))\le -\theta ( |\partial_x W_F^\pm|_{H^{s-1}}^2 
+ | \partial_x w_F^{II\pm}|_{H^s}^2)
+ C|W_F^\pm|_{L^2}^2
\end{equation}
for $C>0$ sufficiently large, for
\begin{equation}
\begin{aligned}
\label{EF}
{\cal{E}_F^\pm}(W)&:=(1/2)\Big(\langle \tilde A^0 W,W\rangle 
+\langle K_\pm(\partial_x) W,W\rangle \Big)
+C\langle \tilde A^0 \partial_x W,\partial_x W\rangle \\
&\quad + C\langle K_\pm(\partial_x) \partial_x W,\partial_x W\rangle \Big)
+C^2 \langle \tilde A^0 \partial_x^2 W,\partial_x^2 W\rangle +  \dots\\
&\quad
+ C^{s} \langle \tilde A^0\partial_x^{s}W,\partial_x^{s}W\rangle,
\end{aligned}
\end{equation}
where operator $K_\pm(\partial_x)$ is defined by
\begin{equation}
\widehat{K_\pm(\partial x f)}(\xi):= iK_\pm(\xi) \hat f(\xi)
\end{equation}
with $K_\pm(\cdot)$ as in \eqref{Krel}, where $\hat g$ denotes
Fourier transform of $g$.

To estimate $W_N$, we use a different technique, substituting
Goodman- for Kawashima-type estimates in controlling $w_N^I$ terms.
Combining calculations \eqref{v0friedrichs} and \eqref{sorder},
we obtain
\begin{equation}
\label{v0pert}
\begin{aligned}
\frac{1}{2}&\langle W_N, \alpha \tilde A^0 W_N\rangle_t 
=
\frac{1}{2}\langle W_N, \alpha \tilde A^0_t W_N\rangle+
 \langle W_N, \alpha \tilde A^0 (W_N)_t \rangle\\
&=
\frac{1}{2}\langle W_N, \alpha \tilde A^0_t W_N\rangle 
 -\langle W_N, \sum_j\alpha \tilde A^j (W_N)_{x_j} 
 + \sum_{j,k}\alpha (\tilde B^{jk}(W_N)_{x_k})_{x_j} \rangle\\
&=
\Big\{
\frac{1}{2} 
\langle w_N^I, (\alpha \tilde A^1_{11}(\bar W))_{x_1} w_N^I\rangle 
- \sum_{j,k} \langle (w_N^{II})_{x_j}, 
\alpha \tilde B^{jk}(w_N^{II})_{x_k}\rangle \Big\}\\
&\quad
+\frac{1}{2}\langle W_N, \alpha \tilde A^0_t W_N\rangle 
+ \langle w_N^I, \sum_j \big(\alpha 
\big(\tilde A_{11}^j- \tilde A^j_{11}(\bar W)\big) \big)_{x_j} w_N^I\rangle \\
&\quad
+ 
\frac{1}{2} 
\langle w_N^I, \sum_j \alpha \tilde A_{12}^j (w_N^{II})_{x_j}\rangle 
-\frac{1}{2} 
\langle (\alpha w_N^{II})_{x_j}, \sum_j \tilde A_{21}^j w_N^{I}\rangle 
\\
&\quad
+\frac{1}{2} 
\langle w_N^{II}, \sum_j \alpha \tilde A_{22}^j (w_N^{II})_{x_j}\rangle 
- \sum_{k} \langle w_N^{II}, 
\alpha_{x_1} \tilde B^{1k}(w_N^{II})_{x_k}\rangle \Big\} \\
&\le
-\Big\{C_* \langle w_N^I, \alpha \Theta w_N^I\rangle 
+\theta \langle \partial_x w_N^{II},
\alpha \partial_x w_N^{II}\rangle \Big\}\\
&\quad
+ C(|W_t|_\infty +|W_x|_\infty)\int \alpha |W|^2
+ C \int \alpha |\partial_x w_N^{II}||W| \\
&\quad
+ C C_* \int \Theta \alpha |w_N^{II}|\big(|w_N^I|+ |\partial_x w_n^{II}|\big),
\\
\end{aligned}
\end{equation}
for some $c>0$, $C>0$ independent of $C_*$.
Choosing $C_*$ sufficiently large that
$C_*>>C(1+|W|_{H^s})$ 
and $C_*>>C^2/\theta \Theta$ on the support $[-2M,2M]$ of $W_N$,
and estimating the unbracketed terms in the final line of \eqref{v0pert}
using Young's inequality, we find that
\begin{equation}
\label{v0pert2}
\begin{aligned}
\frac{1}{2}\langle W_N, \alpha \tilde A^0 W_N\rangle_t 
&\le
-\frac{C_*}{2}\langle w_N^I, \alpha \Theta w_N^I\rangle 
-\frac{\theta}{2} \langle \partial_x w_N^{II},
\alpha \partial_x w_N^{II}\rangle \\
&\quad + C(C_*) \langle w_N^{II}, \alpha w_N^{II}\rangle.
\end{aligned}
\end{equation}

More generally, defining
\begin{equation}
\label{EN}
\CalE_N(W):=
(1/2)\sum_{r=0}^s
c^{-r}\langle \partial_x^r W, \alpha \tilde A^0 \partial_x^r W\rangle,
\end{equation}
for $c=c(C_*)\ge 4C_*/\theta>0$ sufficiently large, 
we obtain by telescoping sum that
\begin{equation}
\label{ENexp}
\begin{aligned}
(d/dt)&\CalE_N(W_N(t))\le  
\sum_{r=0}^s c^{-r}
\Big(-\frac{C_*}{2}\langle \partial_x^r w_N^I, 
\alpha \Theta \partial_x^r w_N^I\rangle \Big) \\
&\quad +
\sum_{r=0}^s c^{-r}
\Big(-\frac{\theta}{2} \langle \partial_x^{r+1} w_N^{II},
\alpha \partial_x^{r+1} w_N^{II}\rangle 
+ C(C_*) \langle \partial_x^{r}w_N^{II}, 
\alpha \partial_x^{r} w_N^{II}\rangle\Big)\\
& \le -\theta_1 \big( |\partial_x W_N|_{H^{s-1}}^2 
+ | \partial_x w_N^{II}|_{H^s}^2\big)
+ C_1|W_N|_{L^2}^2,
\end{aligned}
\end{equation}
$\theta_1>0$ and $C_1>0$ depending on $C_*$,
where in the second inequality of
\eqref{ENexp} we have used the fact that $\alpha$ and $\Theta$
are bounded above and below on the support $[-2M, 2M]$
of $W_N$.

Summing \eqref{EFexp} and \eqref{ENexp}, and using the fact
that $(\CalE_F^\pm)^{1/2}$ and $(\CalE_N)^{1/2}$ 
are each equivalent to $H^s$, we obtain
\begin{equation}
(d/dt){\cal{E}}(W(t))\le -\theta_2 {\cal{E}}(W(t)) + C_2|W(t)|_{L^2}^2,
\end{equation}
where 
\begin{equation}
\CalE(W):= \CalE_F^-(W_F^-)+ \CalE_N(W_N) + \CalE_F^+(W_F^+)
\end{equation}
with $\CalE^{1/2}(W)$ equivalent to $|W|_{H^s}$.
Applying Gronwall's inequality, 
similarly as in the proof of Proposition \ref{ccKenergy},
we obtain
$e^{\theta_2 t}\calE |^t_0\le C \int_0^t e^{\theta_2 s}|W(s)|_{L^2}^2\, ds$
and thereby the result.
\end{proof}

\begin{rem}
\label{ptwiserem}
\textup{
Note that \eqref{auxexp} gives pointwise-in-time control
on the $H^s(x)$ norm of $W$ in terms of a time-weighted
$L^2(t)$ average of the $L^2(x)$ norm.
We will take advantage of this fact in the linearized
analysis of Section 4.3, where we establish $L^2$-time-averaged 
high-frequency estimates by a simple Parseval argument,
thus avoiding the complicated pointwise bounds of [MaZ.4, Z.4].
%
}
\end{rem}

\medbreak
{\bf 4.2. Linearized estimate.}
We next establish the following estimate on the
linearized inhomogeneous problem
\begin{equation}
\label{inhom2}
U_t-LU=f_1 + \partial_x f_2, \qquad U(0)=U_0,
\end{equation}
where $L$ is defined as in \eqref{linearized},
assuming the low-frequency bounds of [Z.3].

\begin{prop}
\label{linest}
Given (A1)--(A2), (H0)--(H5), assumption \eqref{lax}, 
and strong spectral, structural, and strong refined
dynamical stability, the solution
$U(t):=\int_0^te^{L(t-s)}f(s)ds$ of \eqref{inhom2} satisfies
\begin{equation}
\label{linesteq}
\begin{aligned}
\int_0^T e^{-8\theta_1 (T-s)}|U(s)|_{L^2}^2\, ds &\le
C (1+T)^{-(d-1)/2+ 2\epsilon}|U_0|_{L^1}^2\\
&\quad+
C e^{-2\theta_1 T} \big(|U_0|_{H^1}^2+|LU_0|_{H^1}^2\big)
\\
&\quad+
C \int_0^T e^{-2\theta_1(T-s)}
\Big(|f_1(s)|_{H^1}^2+ |\partial_x f_2|_{H^1}^2\Big)
\, ds\\
&\quad +
C \Big(\int_0^T (1+T-s)^{-(d-1)/4+ \epsilon}
|f_1(s)|_{L^1}\, ds\Big)^2 \\
&\quad +
C \Big(\int_0^T (1+T-s)^{-(d-1)/4+ \epsilon-1/2}
|f_2(s)|_{L^1}\, ds\Big)^2 \\
\end{aligned}
\end{equation}
for any fixed $\epsilon>0$, for some $C=C(\epsilon)>0$ sufficiently large,
for any $U_0$, $f_1$, $f_2$ such that the righthand side
is well-defined.
In the case of strong (inviscid, not refined) dynamical stability,
we may take $\epsilon=0$ and $C>0$ fixed.
\end{prop}

\begin{rem}
\label{ptwiserem2}
\textup{
Somewhat stronger linearized bounds 
were obtained in [Z.4] by a more detailed analysis
(pointwise in time, with no loss in
regularity, and for the initial-data rather than the zero-data
inhomogeneous problem).
However, these will suffice for our argument and are much easier
to obtain.
}
\end{rem}

\medbreak
{\bf 4.2.1.  High- and mid-frequency resolvent bounds.}
Our first step is to estimate solutions of resolvent equation \eqref{res}.

\begin{prop}\label{HF} 
(High-frequency bound)
Given (A1)--(A2) and (H0)--(H1),
\begin{equation}
\label{HFres}
|(\lambda-L_\Txi)^{-1}|_{\hat H^1(x_1)}\le C
\quad 
\text{\rm for $|(\tilde \xi, \lambda)|\ge R$ and $\R \lambda \ge -\theta$},
\end{equation}
for some $R$, $C>0$ sufficiently large and $\theta >0$ sufficiently small,
where $|\cdot|_{\hat H^1}$ denotes $\hat H^1$ operator norm,
$|f|_{\hat H^1}:= |(1+ |\partial_{x_1}|+|\Txi|)f|_{L^2}$.
\end{prop}
\begin{proof}
A Laplace--Fourier transformed version (now with respect to
$(t, \tilde x)$) of nonlinear estimate \eqref{ENexp}, $s=1$, 
carried out on the linearized equations \eqref{linearized} 
written in $W$-coordinates, $W:=(\partial W/\partial U)(\bar U)U$, 
yields
\begin{equation}
\label{auxHF}
\begin{aligned}
\R \lambda \Big( (1+ |\Txi|^2)|W|^2 + &|\partial_{x_1} W|^2\Big)\le
-\theta_1 \Big( |\Txi|^2|W|^2 + |\partial_{x_1} W|^2\Big) \\
&\quad + 
C_1\Big(|W|^2 + (1+ |\Txi|^2)|W||f| + |\partial_{x_1}W||\partial_{x_1}f|\Big)
\end{aligned}
\end{equation}
for some $C_1>0$ sufficiently big and $\theta_1>0$ sufficiently small,
where $|\cdot|$ denotes $|\cdot|_{L^2(x_1)}$ and
\begin{equation}
(\lambda-L_\Txi)U=f.
\end{equation}
We omit the proof, which is just the translation into frequency
domain of the proof of Proposition \ref{auxenergy}, carried out
in the much simpler, linearized setting.
For related calculations, see the proof of Proposition \ref{specres}
or the proof of (K3) in Lemma \ref{KSh}.

Rearranging, and dividing by factor $|(1+|\Txi|+|\partial_{x_1})|W|$,
we have
\begin{equation}
\label{auxHF2}
\begin{aligned}
(\R \lambda +\theta_1)
|( (1+ |\Txi|+ |\partial_{x_1}|)W| &\le
C_1|( (1+ |\Txi|+ |\partial_{x_1}|)f| + C_1|W|,
\end{aligned}
\end{equation}
which yields \eqref{HFres} with $\theta:=\theta_1/2$ 
(i.e, $\R \lambda +\theta_1 \ge \theta_1/2>0$)
and $C:=2C_1/\theta$,
for $|(\Txi, \lambda)|\ge 2C_1$ 
(i.e., sufficiently large that $C_1W$ may be absorbed in the lefthand side
of \eqref{auxHF2}).
\end{proof}
%

\begin{prop}\label{MF}
(Mid-frequency bound)
Given (A1)--(A2), (H0)--(H1), and strong spectral stability
\eqref{strongevans}, 
\begin{equation}
\label{MFres}
|(\lambda-L_\Txi)^{-1}|_{\hat H^1(x_1)}\le C
\quad 
\text{\rm for $ R^{-1}\le |(\tilde \xi, \lambda)|\le R$ 
and $\R \lambda \ge -\theta$},
\end{equation}
for any $R>0$, and $C=C(R)>0$ sufficiently large and $\theta=\theta(R) >0$ 
sufficiently small,
where $|\cdot|_{\hat H^1}$ is as defined in Proposition \ref{HF}.
\end{prop}
\begin{proof}
This follows by compactness of the set of frequencies under
consideration together with the 
fact that the resolvent $(\lambda-L_{\Txi})^{-1}$ is analytic with
respect to $\hat H^1$ in 
$(\Txi,\lambda)$ for $\lambda$ in the resolvent set $\rho_{\hat H^1}(L)$
of $L_\Txi$ with respect to $\hat H^1$
once we establish that $\lambda$ lies in the 
$\rho_{\hat H^1}(L_\Txi)$ whenever
$R^{-1}\le |(\tilde \xi, \lambda)|\le R$ 
and $\R \lambda \ge -\theta$.
Here, the $\hat H^1$-resolvent set 
$\rho_{\hat H^1}(L_\Txi)$ is defined as the set of $\lambda$ such 
that $(\lambda-L_\Txi)$
has a bounded inverse with respect to the $\hat H^1$-operator norm
[Kat, Yo, Pa], and similarly for $\rho_{L^2}(L_\Txi)$,
to which we refer elsewhere simply as $\rho(L_\Txi)$; see
Definition \ref{spectrum}, Appendix A.
Analyticity in $\lambda$ on the resolvent set of $(\lambda-L)^{-1}$ 
holds for general operators $L$ (Lemma \ref{analytic}, Appendix A), 
while analyticity in $\Txi$ on the resolvent set of $(\lambda -L_\Txi)^{-1}$
follows using standard properties of asymptotically constant-coefficient
ordinary differential operators [He, Co, ZH] from
analyticity in $\Txi$ of the limiting, constant-coefficient symbol
together with convergence at integrable rate of the coefficients
of $L_\Txi$ to their limits as $x_1\to \pm \infty$.
Alternatively, we may establish this directly for $\theta$
sufficiently small, using the smoothing properties induced
by (A1)--(A2), (H0)--(H2);
see Exercises \ref{relative}.1--2, Appendix A.

Under assumptions (A1)--(A2), (H0)--(H2), 
$\rho_{\hat H^1}(L_\Txi)$ and 
$\rho_{L^2}(L_\Txi)$ agree on the set
$R^{-1}\le |(\tilde \xi, \lambda)|\le R$, $\R \lambda \ge -\theta$,
provided $\theta$ is chosen
sufficiently small relative to $R^{-1}$.
For, 
standard considerations yield that, 
on the ``domain of consistent
splitting'' comprised of the intersection of the rightmost components of the
resolvent sets of the limiting, constant-coefficient operators $L_\Txi^\pm$
as $x_1\to \pm \infty$, 
the spectra of the asymptotically constant-coefficient ordinary differential
operators $L_\Txi$ with respect to $\hat H^s$,
for any $0\le s\le q(d)$, $q(d)$ as defined in (H0) 
(indeed, with respect to any reasonable norm), consist entirely of
isolated $L^2$ eigenvalues, i.e., 
$\lambda-L_\Txi$ is Fredholm in each $\hat H^s$;
for further discussion, see [He, GZ, ZH, ZS, Z.3--4] or Section 5, below.
By (K3), the domain of consistent splitting includes
the domain under consideration for $\theta$ sufficiently small,
giving $\rho_{\hat H^1}(L_\Txi)= \rho_{L^2}(L_\Txi)$ as claimed.
Alternatively, we may establish this by the direct but more
special argument that \eqref{auxHF2} together with
$|(\lambda-L_\Txi)^{-1}|_{L^2(x_1)}\le C_1$ implies 
$|(\lambda-L_\Txi)^{-1}|_{\hat H^1(x_1)}\le C_2$.

With this 
observation,
the result follows by the 
strong spectral stability assumption \eqref{strongspec},
since analyticity of the resolvent in $(\Txi, \lambda)$
implies that $\{(\Txi,\lambda): \, \lambda\in \rho_{L^2}(L_\Txi)\}$
is open, and therefore contains an open neighborhood
of $\{(\Txi,\lambda):\, \R \lambda \ge 0 \} \setminus \{(0,0)\}$.
(Recall, $\sigma_{L^2}(L):= \rho_{L^2}(L)^{c}$.)
\end{proof}

\begin{rems}
\label{linrems}
\textup{
1. As the argument of Proposition \ref{MF} suggests, uniform 
mid-frequency bounds are essentially automatic in the context of 
asymptotically constant-coefficient ordinary differential operators, 
following from compactness/continuity of the resolvent alone.
}
\smallbreak
\textup{
2.  Resolvent bounds \eqref{HFres}
for $\lambda$ with negative real part, correspond roughly to time-exponential 
decay in high-frequency modes, similarly as in \eqref{auxexp}.
Indeed, a uniform bound $|(\lambda-L_\Txi)^{-1}|_{\hat H^1}\le C$
for all $\Txi$, $\R \lambda \ge -\theta<0$, or equivalently
$|(\lambda-L)^{-1}|_{H^1}\le C$
for all $\R \lambda \le -\theta<0$ would imply time-exponential
linearized stability $|e^{Lt}|_{H^1}\le C(\tilde \theta) e^{-\tilde \theta t}$
for any $\tilde \theta <\theta$, by a general Hilbert-space
theorem of Pr\"uss [Pr]; see Remark \ref{pruss}, Appendix A.
Contrasting with estimate \eqref{res} corresponding to local
well-posedness, what we have done here is to use the additional structure
(A2), (H1) to push $\gamma_0$ into the negative half-plane.
Note that we made no stability assumptions in the 
statement of
Proposition \ref{HF}.
}

\textup{
3.  High-frequency linearized resolvent bounds may alternatively be 
obtained by Kreiss symmetrizer techniques applying to a 
much more general class of systems; see [GMWZ.4].
In particular, we may dispense with the symmetry assumptions
in (A1)--(A2) connected with energy estimates using integration by parts, 
replacing them everywhere with sharp, spectral conditions
of hyperbolicity/parabolicity.
By corresponding pseudodifferential estimates, it might be possible 
to recover also the nonlinear time-evolutionary estimate \eqref{Wauxexp}
for this more general class of equations.
This would be a very interesting direction for future investigation.
}
\end{rems}

\medbreak
{\bf 4.2.1. Splitting the solution operator.}
Letting $L$ as usual denote the linearized operator
defined in \eqref{linearized}, define a ``low-frequency cutoff''
\begin{equation}
\label{S1def}
\begin{aligned}
S_1&(t)f(x):=\\
&\int_{|\Txi|^2\le \theta_1+\theta_2} 
\oint_{ \R \lambda=\theta_2 - |\Txi|^2-|\Im \lambda|^2 \ge -\theta_1}
e^{i\tilde \xi \cdot \tilde x + \lambda t}
(\lambda- L_\Txi)^{-1}\hat f(x_1, \Txi) \, d\tilde \xi d\lambda 
\end{aligned}
\end{equation}
of the solution operator $e^{Lt}$ described in \eqref{ILFT},
similarly as in the constant-coefficient argument carried
out in Appendix B, where $\theta_1$, $\theta_2>0$ are
to be determined later.
From Propositions \ref{HF}--\ref{MF}, we obtain 
the following decompositions justifying this analogy.

\begin{cor}\label{decomposition}
Given (A1)--(A2), (H0)--(H1), and strong spectral stability,
let $|e^{L t}|\le Ce^{\gamma_0 t}$ denote the $C^0$ semigroup on $L^2$ 
generated by $L$, with domain $\CalD(L):=\{ U: \, U,\, LU \in L^2\}$.
Then, for arbitrary $\theta_2>0$, $\theta_1>0$ sufficiently
small relative to $\theta_2$, and all $t\ge 0$,
\begin{equation}
\label{cutoffILFT}
\begin{aligned}
&e^{Lt}f(x)= S_1(t)f(x)\\
&\quad +
\pv \int_{-\theta_1-i\infty}^{-\theta_1+i\infty}
\int_{\RR^{d-1}} 
\II_{|\Txi|^2+|\Im \lambda|^2|\ge \theta_1+\theta_2}
e^{i\tilde \xi \cdot \tilde x + \lambda t}
(\lambda- L_\Txi)^{-1}\hat f(x_1, \Txi) \, d\tilde \xi d\lambda \\
\end{aligned}
\end{equation}
for all $f\in \CalD(L)\cap H^1$,
where $\II_P$ denotes the indicator function for
logical proposition $P$ (one for $P$ true and zero for $P$ false)
and $\hat f$ denotes Fourier transform of $f$.
Likewise, for $0\le t\le T$,
\begin{equation}
\label{cutoffinhomILFT}
\begin{aligned}
&\int_0^t e^{L(t-s)}f(s) \, ds =
\int_0^T S_1(t-s) f(s) \, ds \\
&+\pv \int_{-\theta_1-i\infty}^{-\theta_1+i\infty}
\int_{\RR^{d-1}} 
\II_{|\Txi|^2+|\Im \lambda|^2|\ge \theta_1+\theta_2}
e^{i\tilde \xi \cdot \tilde x + \lambda t}
(\lambda- L_\Txi)^{-1}
{\widehat {\widehat {f^T}}}(x_1, \Txi, \lambda) \, d\tilde \xi d\lambda \\
\end{aligned}
\end{equation}
for all $f \in L^1([0,T]; \CalD(L)\cap H^1(x))$, 
where $\hat{\hat f}$ denotes Laplace--Fourier transform of $f$ and
truncation $f^T$ is as in \eqref{truncation}.
\end{cor}

\begin{proof}
Using analyticity on the resolvent set $\rho(L)$
of the resolvent $(\lambda-L)^{-1}$ (Lemma
\ref{analytic}, Appendix A), 
the bound 
\begin{equation}
|(\lambda-L)^{-1}f|
\le C(|(\lambda-L)^{-1}||Lf| + |f|)|\lambda|^{-1}
\end{equation}
coming from resolvent identity \eqref{residentbd} (Lemma \ref{resident}),
and the results of Propositions \ref{HF}--\ref{MF}, we 
may deform the contour in \eqref{ILT1} using Cauchy's Theorem
to obtain
\begin{equation}
\label{1step}
\begin{aligned}
&e^{Lt}f(x)= 
\oint_{ \R \lambda=\theta_2 - |\Im \lambda|^2;  \,
|\Im \lambda|^2\le \theta_1+\theta_2}
e^{\lambda t}
(\lambda- L)^{-1} f(x) \, d\lambda \\
&\quad +
\pv \int_{-\theta_1-i\infty}^{-\theta_1+i\infty}
\II_{|\Im \lambda|^2|\ge \theta_1+\theta_2}
e^{\lambda t}
(\lambda- L)^{-1} f(x) \, d\lambda \\
&=\oint_{ \R \lambda=\theta_2 - |\Im \lambda|^2;  \,
|\Im \lambda|^2\le \theta_1+\theta_2}
\int_{\Txi\in \RR^{d-1}}
e^{i\Txi\cdot \tilde x+\lambda t}
(\lambda- L_\Txi)^{-1}\hat f(x_1,\Txi) \, d\Txi d\lambda \\
&\quad +
\pv \int_{-\theta_1-i\infty}^{-\theta_1+i\infty}
\II_{|\Im \lambda|^2|\ge \theta_1+\theta_2}
\int_{\Txi\in \RR^{d-1}}
e^{i\Txi\cdot \tilde x+\lambda t}
(\lambda- L_\Txi)^{-1}\hat f(x_1,\Txi) \, d\Txi d\lambda. \\
\end{aligned}
\end{equation}
Using exercise \ref{ftconv} below, we may 
exchange the order of integration to write the 
first term of the last line as
\begin{equation}
\label{2step}
\begin{aligned}
\int_{\Txi\in \RR^{d-1}}
\oint_{ \R \lambda=\theta_2 - |\Im \lambda|^2;  \,
|\Im \lambda|^2\le \theta_1+\theta_2}
e^{i\Txi\cdot \tilde x+\lambda t}
(\lambda- L_\Txi)^{-1}\hat f(x_1,\Txi) \, d\Txi d\lambda. \\
\end{aligned}
\end{equation}
Deforming individual contours in the second integral, 
using analyticity of $(\lambda-L_\Txi)^{-1}$ on the resolvent
set $\rho(L_\Txi)$
together with the results of Propositions \ref{HF}--\ref{MF}, and
exchanging orders of integration using the result of Exercise
\ref{ftconv}, we may rewrite this further as
\begin{equation}
\label{3step}
\begin{aligned}
&\int_{|\Txi|^2 \le \theta_1+ \theta_2}
\oint_{ \R \lambda = \theta_2-|\Im \lambda|^2- |\Txi|^2}
e^{i\Txi\cdot \tilde x+\lambda t}
(\lambda- L_\Txi)^{-1}\hat f(x_1,\Txi) \, d\Txi d\lambda \\
&\quad + \int_{\Txi\in \RR^{d-1}}
\int_{-\theta_1-i\sqrt{\theta_1+\theta_2}}^{-\theta_1+i\sqrt{\theta_1+\theta_2}}
\II_{|\Im \lambda|^2|+ |\Txi|^2\ge \theta_1+\theta_2}
e^{i\Txi\cdot \tilde x+\lambda t}
(\lambda- L_\Txi)^{-1}\hat f(x_1,\Txi) \, d\Txi d\lambda \\
&=\int_{|\Txi|^2 \le \theta_1+ \theta_2}
\oint_{ \R \lambda = \theta_2-|\Im \lambda|^2- |\Txi|^2}
e^{i\Txi\cdot \tilde x+\lambda t}
(\lambda- L_\Txi)^{-1}\hat f(x_1,\Txi) \, d\Txi d\lambda \\
&\quad + 
\int_{-\theta_1-i\sqrt{\theta_1+\theta_2}}^{-\theta_1+i\sqrt{\theta_1+\theta_2}}
\int_{\Txi\in \RR^{d-1}}
\II_{|\Im \lambda|^2|+ |\Txi|^2\ge \theta_1+\theta_2}
e^{i\Txi\cdot \tilde x+\lambda t}
(\lambda- L_\Txi)^{-1}\hat f(x_1,\Txi) \, d\Txi d\lambda. \\
\end{aligned}
\end{equation}
Substituting back into \eqref{1step} and recombining terms, we obtain
\eqref{cutoffILFT}.
A similar computation, together with the observation that
\begin{equation}
{\widehat{\widehat {f^T}}}(x_1,\Txi, \lambda):=
\int_0^T e^{-\lambda s} \hat f(x_1,\Txi, s)\, ds
\end{equation}
is analytic in $(\Txi, \lambda)$ as the absolutely convergent integral
of an analytic integrand, yields \eqref{cutoffinhomILFT}.
\end{proof}

\begin{ex}\label{ftconv}
For $g\in H^{1}(x)$, show using Parseval's identity and the
Fourier transform definition
$|g|_{H^s(x)}:= |(1+ |\xi|^2)^{s/2} \hat g(\xi)|_{L^2(\xi)}$
that
\begin{equation}
\label{ftbd}
\int_{|\Txi|\ge L}
|\hat g|_{L^2(x_1)}^2 \, d\Txi \le C(1+L^2)^{-1})|g|_{H^{1}(x)}^2
\end{equation}
converges uniformly to zero as $L\to +\infty$, where $\hat g(x_1,\Txi)$ 
denotes the Fourier transform of $g$ with respect to $\tilde x$.
\end{ex}

\begin{rem}
\label{causality}
\textup{
Defining the ``high-frequency cutoff'' $S_2(t):=e^{Lt}-S_1(t)$,
we have $\int_0^t e^{L(t-s)}ds= \int_0^t S_1(t-s) ds
+ \int_0^t S_2(t-s) ds$, where, by \eqref{cutoffinhomILFT},
\begin{equation}
\label{causeq}
\begin{aligned}
& \int_0^t S_2(t-s)f(s)\, ds=
\int_t^T S_1(t-s) f(s) \, ds \\
&+\pv \int_{-\theta_1-i\infty}^{-\theta_1+i\infty}
\int_{\RR^{d-1}} 
\II_{|\Txi|^2+|\Im \lambda|^2|\ge \theta_1+\theta_2}
e^{i\tilde \xi \cdot \tilde x + \lambda t}
(\lambda- L_\Txi)^{-1}
{\widehat {\widehat {f^T}}}(x_1, \Txi, \lambda) \, d\tilde \xi d\lambda. \\
\end{aligned}
\end{equation}
Note that $S_1(t)= -S_2(t)\ne0$ for $t<0$, in violation of causality.
``Causality error'' $\int_t^T S_1(t-s)f(s)\, ds$ 
is the price for fixing the cutoff $T$ in the second term of 
\eqref{causeq}.
}
\end{rem}

\medbreak
{\bf 4.2.3. Proof of the main estimate.}
Now, assume the following low-frequency estimate of [Z.3],
to be discussed further in Section 5.
Similar results hold for profiles of nonclassical, 
over- or undercompressive type, as described in [Z.3--4].

\begin{prop}[{[Z.3]\footnote{
For a more detailed exposition, see [Z.4].}}]
\label{S1}
Under the assumptions of Proposition \ref{linest},
\begin{equation}
\label{S1eq}
\begin{aligned}
|S_1(t)(f_1 + \partial_x f_2)|_{L^2(x)}&\le
C(1+t)^{-(d-1)/4+ \epsilon} |f_1|_{L^1(x)}\\
&\quad + C (1+t)^{-(d-1)/4+ \epsilon-1/2} |f_2|_{L^1(x)}
\end{aligned}
\end{equation}
for all $t\ge 0$,
for any fixed $\epsilon>0$, for some $C=C(\epsilon)>0$ sufficiently large,
for any $f_1\in L^1(x)$ and
$f_2\in L^1(x)$.
In the case of strong (inviscid, not refined) dynamical stability,
we may take $\epsilon=0$ and $C>0$ fixed.
\end{prop}

\begin{cor}
\label{S1cor}
Under the assumptions of Proposition \ref{linest},
\begin{equation}
\label{S1eqv}
\begin{aligned}
V(x,t):=
\int_0^t
S_1(T-t)f(x,s) ds
\end{aligned}
\end{equation}
satisfies
\begin{equation}
\label{S1bd}
\int_0^T e^{-8\theta_1 (T-t)}|V(t)|_{L^2(x)}^2 \, dt \le
C \Big(\int_0^T (1+T-s)^{-(d-1)/4+ \epsilon-1/2}
|f(s)|_{L^1}\, ds\Big)^2 \\
\end{equation}
for all $T\ge 0$, for any $f(x,t)\in L^1([0,T]; L^1(x))$.
\end{cor}

\begin{proof}
Using \eqref{S1eq} and the inequality
\begin{equation}
e^{-2\theta_1 (T-t)} (1+t-s)^{-(d-1)/4+ \epsilon-1/2}\le C
(1+T-s)^{-(d-1)/4+ \epsilon-1/2}
\end{equation}
for $0\le s\le t\le T$, we obtain
\begin{equation}
\begin{aligned}
\int_0^T &e^{-8\theta_1 (T-t)}|V(t)|_{L^2(x)}^2 \, dt =
\int_0^T e^{-8\theta_1 (T-t)}\Big(
\int_0^t S_1(t-s)|f(s)|_{L^1(x)} ds\Big)^2 \, dt \\
&\le
C\int_0^T e^{-8\theta_1 (T-t)}\Big(
\int_0^t (1+t-s)^{-(d-1)/4+ \epsilon-1/2}
|f(s)|_{L^1(x)} ds\Big)^2 \, dt \\
&\le
C\int_0^T e^{-4\theta_1 (T-t)}\Big(
\int_0^T (1+T-s)^{-(d-1)/4+ \epsilon-1/2}
|f(s)|_{L^1(x)} ds\Big)^2 \, dt \\
&\le
C_2\Big(\int_0^T (1+T-s)^{-(d-1)/4+ \epsilon-1/2} |f(s)|_{L^1(x)} ds\Big)^2.
\end{aligned}
\end{equation}
\end{proof}

\begin{lem}
\label{noncausal}
Under the assumptions of Proposition \ref{MF},
\begin{equation}
\label{noncausalS1eq}
\begin{aligned}
|S_1(t)f)|_{L^2(x)}&\le
Ce^{-\theta_1 t} |f|_{H^1}
\end{aligned}
\end{equation}
for all $t< 0$, for any $f\in H^1(x)$.
\end{lem}

\begin{proof}
Immediate, from representation \eqref{S1def} and Proposition \ref{MF}.
\end{proof}

\begin{cor}
\label{noncausalcor}
Under the assumptions of Proposition \ref{MF},
\begin{equation}
\label{noncausalinhomS1eq}
\begin{aligned}
W_1(x,t):=
\int_t^T
S_1(T-t)f(x,s) ds
\end{aligned}
\end{equation}
satisfies
\begin{equation}
\label{S2bd1}
\int_0^T e^{-8\theta_1 (T-t)}|W_1(t)|_{L^2(x)}^2 \, dt \le
C\int_0^T e^{-2\theta_1 (T-t)}|f(t)|_{H^1(x)}^2 \, dt
\end{equation}
for all $T\ge 0$, for any $f(x,t)\in L^2([0,T], H^1(x))$.
\end{cor}

\begin{proof}
Using the triangle inequality, and \eqref{noncausalS1eq}, we obtain
\begin{equation}
\begin{aligned}
\int_0^T e^{-8\theta_1 (T-t)}|W_1(t)|_{L^2(x)}^2 \, dt &=
\int_0^T e^{-8\theta_1 (T-t)}\Big|
\int_t^T S_1(t-s)f(s) ds\Big|_{L^2(x)}^2 \, dt \\
&\le
C\int_0^T e^{-8\theta_1 (T-t)}\int_t^T
e^{-2\theta_1 (t-s)}|f(s)|_{H^1(x)}^2 \, ds dt \\
&\le
C\int_0^T e^{-4\theta_1 (T-t)}\int_0^T
e^{-2\theta_1 (T-t)}|f(s)|_{H^1(x)}^2 \, ds dt \\
&\le
C_2\int_0^T
e^{-2\theta_1 (T-t)}|f(s)|_{H^1(x)}^2 \, ds. \\
\end{aligned}
\end{equation}
\end{proof}

\begin{lem}
\label{S2}
Under the assumptions of Proposition \ref{MF},
\begin{equation}
\label{W2}
\begin{aligned}
&W_2(x,t):=\\
&\quad
\pv \int_{-\theta_1-i\infty}^{-\theta_1+i\infty}
\int_{\RR^{d-1}} 
\II_{|\Txi|^2+|\Im \lambda|^2|\ge \theta_1+\theta_2}
e^{i\tilde \xi \cdot \tilde x + \lambda t}
(\lambda- L_\Txi)^{-1}{\widehat{\widehat {f^T}}}(x_1, \Txi,\lambda) \, d\tilde \xi d\lambda \\
\end{aligned}
\end{equation}
satisfies
\begin{equation}
\label{S2bd2}
\int_0^T e^{-2\theta_1 (T-t)}|W_2(t)|_{L^2(x)}^2 \, dt \le
C\int_0^T e^{-2\theta_1 (T-t)}|f(t)|_{H^1(x)}^2 \, dt
\end{equation}
for all $T\ge 0$, for any $f(x,t)\in L^2([0,T],H^1(x))$.
\end{lem}

\begin{proof}
Using the relation 
\begin{equation}
\label{ftltrel}
\lt h(\lambda)=\ft \big( e^{-\R \lambda t}h\big)(-\Im \lambda)
\end{equation}
between Fourier transform $\ft h$ and Laplace transform
$\lt h$ of a function $h(t)$, 
we have 
\begin{equation}
\ft e^{\theta_1 t}W_2(t)(-k)= \lt W(-\theta_1+ik)=
(\lambda-L)^{-1}\lt f(-\theta_1+ik).
\end{equation}
Together with
Parseval's identity and Propositions \ref{HF}--\ref{MF}, this yields
\begin{equation}
\begin{aligned}
&
\int_0^T e^{2\theta_1 t} \big|W_2(t) \big|^2_{L^2(x)} \, dt 
=
|e^{\theta_1 t}W_2(t)|_{L^2(x,[0,T])}^2=\\
&\quad 
\int_{-\theta_1-i\infty}^{-\theta_1+i\infty}
\Big|\II_{|\Txi|^2+|\Im \lambda|^2|\ge \theta_1+\theta_2}
(\lambda- L_\Txi)^{-1}{\widehat{\widehat{f^T}}}\Big|^2_{L^2(x_1,\Txi)} 
\, d\lambda \\
&\le
C\int_{-\theta_1-i\infty}^{-\theta_1+i\infty}
\big|{\widehat{\widehat{f^T}}}\big|^2_{H^1(x_1,\Txi)} 
\, d\lambda \\
&=
C\int_0^T e^{2\theta_1 t}
\big|f(t) \big|^2_{H^1(x)}
\, dt \\
\end{aligned}
\end{equation}
yielding the claimed estimate upon multiplication by $e^{-2\theta_1 T}$.
\end{proof}

\begin{lem}
\label{S2ivp}
Under the assumptions of Proposition \ref{MF},
\begin{equation}
\label{W0}
\begin{aligned}
&W_0(x,t):=\\
&\quad
\pv \int_{-\theta_1-i\infty}^{-\theta_1+i\infty}
\int_{\RR^{d-1}} 
\II_{|\Txi|^2+|\Im \lambda|^2|\ge \theta_1+\theta_2}
e^{i\tilde \xi \cdot \tilde x + \lambda t}
(\lambda- L_\Txi)^{-1}{\widehat {U_0}}(x_1, \Txi) \, d\tilde \xi d\lambda \\
\end{aligned}
\end{equation}
satisfies
\begin{equation}
\label{S2bd0}
\int_0^T e^{-2\theta_1 (T-t)}|W_0(t)|_{L^2(x)}^2 \, dt \le
C e^{-2\theta_1 T}\big(|U_0|_{H^1(x)}^2+ |LU_0|_{H^1}^2\big) 
\end{equation}
for all $T\ge 0$, for $U_0$, $LU_0 \in H^1(x,[0,T])$.
\end{lem}

\begin{proof}
Without loss of generality taking $\R\lambda \ne 0$, we have by
\eqref{ftltrel},
Parseval's identity, Propositions \ref{HF}--\ref{MF}, 
and the resolvent identity 
\begin{equation}
(\lambda-L_\Txi)^{-1}\widehat{U_0}=
\lambda^{-1}\big((\lambda-L_\Txi)^{-1}L_\Txi \widehat{U_0} + U_0\big)
\end{equation}
(\eqref{residentbd}, Appendix A), that
\begin{equation}
\begin{aligned}
&
\int_0^\infty e^{2\theta_1 t} \big|W_0(t) \big|^2_{L^2(x)} \, dt 
=
|e^{\theta_1 t}W_0(t)|_{L^2(x,t)}^2=\\
&\quad 
\int_{-\theta_1-i\infty}^{-\theta_1+i\infty}
\Big|\II_{|\Txi|^2+|\Im \lambda|^2|\ge \theta_1+\theta_2}
(\lambda- L_\Txi)^{-1}{\widehat{U_0}}\Big|^2_{L^2(x_1,\Txi)} 
\, d\lambda \\
&\le
C\int_{-\theta_1-i\infty}^{-\theta_1+i\infty}
|\lambda|^{-2} \, d\lambda 
\big(|L_\Txi \widehat{U_0}|_{\hat H^1(x_1, \Txi)} 
+ |U_0|_{L^2(x_1, \Txi)}\big)\\
&\le
C \big(|L U_0|_{H^1(x)} + |U_0|_{L^2(x)}\big),\\
\end{aligned}
\end{equation}
from which the claimed estimate follows 
upon multiplication by $e^{-2\theta_1 T}$.
\end{proof}

\begin{proof}[Proof of Proposition \ref{linest}]
Immediate, combining the estimates of Proposition
\ref{S1}, Corollaries
\ref{S1cor} and \ref{noncausalcor},
and Lemmas \ref{S2ivp} and \ref{S2},
and using representation 
\begin{equation}
U(T)=e^{LT}U_0+\int_0^T e^{L(T-t)}(f_1+\partial_x f_2)
\, dt=:S_1(T)U_0 +V+W_0+ W_1+W_2
\end{equation} 
coming from \eqref{cutoffILFT}--\eqref{cutoffinhomILFT}, \eqref{causeq}.
\end{proof}

%

\medbreak
{\bf 4.3. Nonlinear stability argument.}
We are now ready to establish the result of Theorem
\ref{mainsuf}, assuming for the moment the bounds asserted
in Proposition \ref{S1} on the low-frequency solution operator $S_1(t)$.
Define the nonlinear perturbation
\begin{equation}
\label{pert2}
\begin{aligned}
U&:= \tilde U - \bar U, 
\end{aligned}
\end{equation}
where $\tilde U$ denotes a solution of \eqref{viscous} with
initial data $\tilde U_0$ close to $\bar U$.
Then, Taylor expanding about $\bar U$, we may rewrite \eqref{viscous} as
\begin{equation}
\begin{aligned}
U_t - LU&= 
\partial_x Q(U,\partial_x U)
=:  \partial_x f_2,
\qquad
U(0)&=U_0,\\
\end{aligned}
\end{equation}
where 
\begin{equation}
\label{Qbds}
\begin{aligned}
|Q(U,\partial_x U|&\le C(|U||\partial_x U|+|U|^2),\\
|\partial_x Q(U,\partial_x U|&\le C(|U||\partial_x^2 U|+ |\partial_x U|^2
+|U||\partial_x U|),\\
|\partial_x^2 Q(U,\partial_x U|&\le 
C(|U||\partial_x^3 U|+ |\partial_x U||\partial_x^2U|+ |\partial_x U|^2),\\
\end{aligned}
\end{equation}
so long as $U$ remains uniformly bounded.

\begin{lem}\label{fbds}
So long as $U$ remains uniformly bounded,
\begin{equation}
\label{fbdseq}
\begin{aligned}
|f_2(t)|_{L^1(x)}&\le 
 C|U(t)|_{H^1(x)}^2,\\ 
|\partial_x f_2(t)|_{H^1(x)}&\le 
C|U(t)|_{H^{[d/2]+2}(x)}^2.\\
\end{aligned}
\end{equation}
\end{lem}
\begin{proof}
Immediate, by \eqref{Qbds} and 
$|\partial_x U|_{L^\infty(x)}\le |U|_{H^{[d/2]+2}(x)}$.
\end{proof}

\begin{proof}[Proof of Theorem \ref{suff}]
Set $8\theta_1:= \theta_2$, $\theta_2>0$ as in 
the statement of Proposition \ref{auxenergy}, and define
\begin{equation}
\label{zeta}
\zeta(t):=\sup_{0\le \tau \le t}
|U(\tau)|_{L^2}(1+\tau )^{\frac{(d-1)}{4}-\epsilon},
\end{equation}
where $\epsilon>0$ is as in the statement 
of Theorem \ref{mainsuf}. 
By local well-posedness in $H^s$, Proposition \ref{parabwellposed},
and the standard principle of continuation, 
there exists a solution  $U\in H^s(x)$, $U_t\in H^{s-2}(x)\subset L^2(x)$
on the open time-interval for which $|U|_{H^s}$ remains bounded.
On this interval, $\zeta$ is well-defined and continuous. 

Now, let $[0,T)$ be the maximal interval on which $|U|_{H^s(x)}$
remains bounded by some fixed, sufficiently small constant $\delta>0$.
By Proposition \ref{auxenergy},
\begin{equation}
\label{2calc}
\begin{aligned}
|U(t)|_{H^s}^2&\le C |U(0)|^2_{H^s}e^{-\theta t}
+C\int_0^t e^{-\theta_2 (t-\tau )}|U(\tau)|_{L^2}^2\, d\tau \\
&\le C_2\big(|U(0)|^2_{H^s}+ \zeta(t)^2\big) (1+\tau )^{-(d-1)/2+ 2\epsilon}.
\end{aligned}
\end{equation}
Combining \eqref{linesteq} with the bounds of Lemma \ref{fbds},
we thus obtain
\begin{equation}
\label{1calc}
\begin{aligned}
\int_0^t &e^{-\theta_2 (t-\tau )}|U(\tau )|_{L^2}^2\, d\tau  \le
C(1+t)^{-(d-1)/4+ \epsilon}|U_0|_{L^1}^2\\
&\qquad +
Ce^{-2\theta_1 t}\big(|U_0|_{L^2}+|LU_0|_{H^1}\big)\\
&\qquad +
C \int_0^t e^{-2\theta_1(t-\tau )}
|\partial_x f_2(\tau)|_{H^1}^2 \, d\tau  \\ 
&\qquad +
C \Big(\int_0^t (1+t-\tau)^{-(d-1)/4+ \epsilon-1/2}
|f_2(\tau)|_{L^1}\, d\tau\Big)^2 \\
&\quad\le
C (1+t)^{-(d-1)/2+ 2\epsilon} |\tilde U_0-\bar U|_{L^1\cap H^3}^2\\
&\qquad + 
C\zeta(t)^4
 \Big(\int_0^t (1+t-\tau)^{-(d-1)/4+\epsilon -1/2}(1+\tau)^{-(d-1)/2+ 2\epsilon}
\,d\tau\Big)^2\\
&\quad\le
C_2 \Big(|\tilde U_0-\bar U|_{L^1\cap H^3}^2 + C\zeta(t)^4\Big)
(1+\tau)^{-(d-1)/2+ 2\epsilon}
\end{aligned}
\end{equation}
for $d\ge 2$ and $\epsilon=0$, or $d\ge 3$ and $\epsilon$ sufficiently small
(exercise).

Applying Proposition \ref{auxenergy} a second time, we obtain
\begin{equation}
\label{2calc2}
\begin{aligned}
|U(t)|_{H^s}^2&\le 
Ce^{-\theta_2 t}|\tilde U_0-\bar U|_{H^s}^2+
C \int_0^t e^{-\theta_2 (t-\tau)}|U(\tau)|_{L^2}^2\, d\tau\\
&\le C_2 \Big(|\tilde U_0-\bar U|_{L^1\cap H^s}^2 + C\zeta(t)^4\Big)
(1+\tau)^{-(d-1)/2+ 2\epsilon}
\end{aligned}
\end{equation}
and therefore
\begin{equation}
\label{suffclaim}
\zeta(t)\le 
C \Big(|\tilde U_0-\bar U|_{L^1\cap H^s} + C\zeta(t)^2\Big)
\end{equation}
so long as $|U|_{H^s}$ and thus $|U|$ remains uniformly bounded.

From \eqref{suffclaim}, it follows by continuous induction (exercise)
that 
\begin{equation}
\label{continduct}
\zeta(t) \le 2C |U(0)|_{L^1\cap H^s} 
\end{equation}
for $|U(0)|_{L^1\cap H^s}$ sufficiently small,
and thus also
\begin{equation}
|U(t)|_{H^s}\le 2C(1+t)^{-(d-1)/4+\epsilon}|U(0)|_{L^{1}\cap H^s}
\end{equation}
as claimed, on the maximal interval $[0,T)$ for which $|U|_{H^s}<\delta$.  
In particular, $|U(t)|_{H^s}<\delta/2$ for 
$|U(0)|_{L^{1}\cap H^s}$ sufficiently small, so that $T=+\infty$,
and we obtain both global existence and $H^s$ decay at the claimed rate.
Applying \eqref{sobolev}, we obtain the same bound for $|U(t)|_{L^\infty}$,
and thus, by $L^p$ interpolation (formula \eqref{Lpinterp}, Appendix B),
for $|U(t)|_{L^p}$, all $2\le p\le \infty$.
\end{proof}
\medbreak
Together with the prior results of [Z.3], the results of
this section complete the analysis.
We'll review those prior results in the low-frequency analysis of
the next section.

\begin{rems}
\label{boltz}
\textup{
1. 
Sharp rates of decay in $L^p$, $2\le p\le \infty$ may be obtained
by a further, bootstrap argument, as in the related analysis of [Z.4].
}
\smallbreak
\textup{
2.  The approach followed here, yielding high-frequency estimates
entirely from resolvent bounds, also greatly simplifies the analysis
of the one-dimensional case, replacing the detailed pointwise high-frequency
estimates of [MZ.3].  Pointwise low-frequency bounds 
for the moment are still required; however, a new method of
shock-tracking introduced in [GMWZ.3] might yield an alternative
approach based on resolvent bounds.  This would be an interesting direction
for further investigation. 
At the same time, the new approach yields slightly more general
results: specifically, we may drop the assumption that the hyperbolic
convection matrix $A_*(\bar U(x))$ have constant multiplicity with
respect to $x$, and we may reduce the regularity requirement
on the coefficients from $C^5$ to $C^3$
and the regularity of the perturbation from 
$L^1\cap H^3$ to $L^1\cap H^2$.\footnote{
More precisely, to $C^4$ and $L^1\cap H^3$
by the arguments of this article alone.
To obtain the weaker regularity stated, we must substitute
for the uniform $H^1\to H^1$ high- and mid-frequency resolvent bounds 
of Propositions \ref{HF} and \ref{MF} the corresponding
$L^2\to L^2$ bounds obtained in [MaZ.3] by pointwise methods
or in [GMWZ.4] by symmetrizer estimates.
}
}
%
%
%
\smallbreak
\textup{
3.  Kawashima-type energy estimates are available also in
the ``dual'', relaxation case, at least for small-amplitude
shocks.  Thus, we immediately obtain by the methods of this
article, together with prior results of [Z.3], a corresponding
result of nonlinear multi-dimensional stability of small-amplitude 
relaxation fronts.
In the one-dimensional case, an energy estimate has recently
been established for large-amplitude profiles in [MZ.5];
in combination with the methods of this article, this recovers
the large-amplitude nonlinear stability result of [MZ.5] without
the restrictive hypothesis of constant multiplicity of
relaxation characteristics.
Besides being mathematically appealing, this improvement
is important for applications to moment-closure models.
}
\smallbreak
\textup{
4.  The analysis so far, with the exception of the large-amplitude
version of the ``magic'' energy estimate \eqref{auxexp}, is all ``soft''.
The novelty lies, rather, in the argument structure, which for the
first time successfully integrates Parseval- and semigroup-type (i.e.,
direct, pointwise-in-time)
bounds in the analysis of long-time viscous shock stability.
See [KK] for an interesting earlier approach 
based entirely on Parseval's identity 
and Hausdorff--Young's inequality, which may be used
in cases of sufficiently fast decay to establish
$\int_0^\infty |U(s)|_{L^2}^2\, ds<+\infty$;
in the shock-wave setting, this occurs for dimensions $d=3$
and higher [GMWZ.2].
For further discussion/comparison, see Remark \ref{deformed} below.
More generally, our mid- and high-frequency analysis could 
be viewed as addressing the larger problem
of obtaining time-algebraic decay estimates from spectral information
in the absence of either a spectral gap or sectorial-type resolvent bounds;
see Remark \ref{pruss}.
}

\textup{
Note in particular that our conclusions are independent of the
particular setting considered here, depending only on the availability
of nonlinear smoothing-type energy estimates favorable in the highest
derivative.
Thus, our methods may be of use in the study of delicate stability phenomena
arising in other equations, for example, stability of traveling pulse- or 
front-type solutions of nonlinear Schr\"odinger's equation.
Moreover, they are truly multi-dimensional, in the sense that
they do not intrinsically depend on planar structure/decoupling 
of Fourier modes, so could be applied also to the case of a
background wave $\bar u$ varying also in the $\tilde x$ directions:
for example to traveling waves in a channel $(x_1,\tilde x)\in
\RR \times X$, $X\subset \RR^{d-1}$ bounded.
}

\textup{
Indeed, there is no requirement that the resolvent equation
even be posable as an (possibly infinite-dimensional) ODE.
Thus our methods may be of use in interesting situations
for which this standard assumption does not hold: in particular,
{\it stability of irrational-speed semi-discrete traveling-waves 
for non-upwind schemes} (see discussion, treatment of the upwind
case in [B, BHR]), {\it or stability of shock profiles for
Boltzmann equations}.
Of course, in each of the mentioned cases (fully multi-dimensional,
semidiscrete, and kinetic equations), there remains the problem
of carrying out a suitable low-frequency linearized analysis;
nonetheless, it is a substantial reduction of the problem.
}
\end{rems}

\bigbreak
\clearpage

\section{Low frequency analysis/completion of the proofs}

We complete our analysis by carrying out the remaining, low-frequency
analysis, establishing Propositions \ref{S1} and \ref{ZS}, 
and Theorem \ref{mainnec}.
Our main tools in this endeavor 
are the conjugation lemma, Lemma \ref{conjugation},
by which we reduce the resolvent equation to a constant-coefficient
equation with implicitly determined transmission condition, 
and the Evans function, through which we in effect obtain asymptotics
for the resulting, unknown transmission condition in the low-frequency limit
(Proposition \ref{ZS}).
\medbreak

{\bf 5.1. Bounds on $S_1(t)$.}
Continuing in our backwards fashion, 
we begin by establishing the deferred low-frequency bounds cited in 
Proposition \ref{S1} of the previous section
under appropriate resolvent bounds,
thus reducing the remaining analysis
to a detailed study of the resolvent and eigenvalue equations, 
to be carried out in the remainder of this section.
The following low-frequency estimates were obtained in [Z.3] by
pointwise estimates on the resolvent kernel;
we sketch a different proof in Section 5.4, based on 
degenerate symmetrizer estimates as in [GMWZ.1, W].

\begin{prop}
[{[Z.3--4]}]
\label{mresbounds}
Assuming (A1)--(A3), (H0)--(H5), for 
$\rho:=|(\Txi,\lambda)|$ sufficiently small, and 
$$
\Re \lambda = 
-\theta (|\Txi|^2+|\Im \lambda|^2) 
$$
for $\theta$ sufficiently small,
there holds
\begin{equation}
\label{resbdeq}
|(L_\Txi-\lambda)^{-1}\partial_{x_1}^{\beta} f|_{L^p(x_1)}
\le C\gamma_1\gamma_2 \rho^{(1-\alpha)|\beta|-1} |f|_{L^1(x_1)}
\end{equation}
for all $2\leq p\leq \infty$, $0\le |\beta|\le 1$,
where 
\begin{equation}
\label{4.13.4a}
\gamma_1 (\Txi,\lambda ):=
\begin{cases}
1 \ \hbox{in case of strong (uniform) dynamical stability,}\\
1+ \sum_{j} 
[\rho^{-1}|\Im \lambda -i\tau_j (\Txi)| +\rho ]^{-1} 
\ {\hbox{otherwise}},
\end{cases}
\end{equation}
\begin{equation}
\label{gamma2}
\gamma_2(\Txi,\lambda):= 
1+ \sum_{j,\pm} [\rho^{-1}|\Im \lambda - \eta_j^\pm(\Txi)| + \rho]^{-t}, 
\qquad 0<t<1,
\end{equation}
$\eta_j(\cdot)$ is as in (H5) and the zero-level
set of $\Delta(\Txi, i\tau)$ for $\Txi$, $\tau$ real
is given by $\cup_{j, \Txi} (\Txi, i\tau_j(\Txi))$, and
\begin{equation}
\label{4.13.5a}
\alpha:=
\begin{cases}
0\ \hbox{for Lax or overcompressive case,}\\
1\ \hbox{for undercompressive case.}
\end{cases}
\end{equation}
(For a definition of overcompressive and undercompressive types,
see, e.g., Section 1.2, [Z.4].)
More precisely, $t:=1-1/K_{\max}$, where $K_{\max}:=\max K^\pm_j=\max s_j^\pm$ 
is the maximum among the orders of all branch singularities 
$\eta_j^\pm(\cdot)$, $s_j^\pm$ and $\eta_j^\pm$ defined as in (H5); 
in particular, $t=1/2$ in the (generic) case that only square-root
singularities occur.
\end{prop}

\begin{rem}\label{comparison}
\textup{
Rewriting the $L^1\to L^2$ resolvent bound \ref{resbdeq} as
\begin{equation}\label{cf}
|(L_\Txi-\lambda)^{-1} f|_{L^2(x_1)}
\le \frac{C\gamma_2 {|f|_{L^1(x_1)}}}{\sqrt{\R \lambda +\rho^2}},
\end{equation}
we may recognize it as essentially a second-order correction of the 
corresponding $L^2\to L^2$ resolvent bound
\begin{equation}\label{inviscidresbd}
|(L_\Txi-\lambda)^{-1} f|_{L^2(x_1)}
\le \frac{C {|f|_{L^2(x_1)}}}{\R \lambda }
\end{equation}
of the inviscid theory, which in turn is approximately the condition
for a bounded $C^0$ semigroup (Appendix A).
The singular factor $\gamma_2$ appearing in the numerator of
\eqref{cf}, and the square root in the denominator, 
are new effects connected with the fact that the bound is taken
between $L^1$ and $L^2$ rather than in $L^2$ alone; for further
discussion, see [Z3].
}
\end{rem}

\begin{proof}[Proof of Proposition \ref{S1}]
For simplicity, restrict to the Lax-type, uniformly
dynamically stable case; other cases are similar.
By analyticity of the resolvent on the resolvent set, we may deform
the contour in \eqref{S1} to obtain 
\begin{equation}
\begin{aligned}
\label{S1new}
S_1&(t)f(x)=\\
&\int_{|\Txi|^2\le \theta_1+\theta_2} 
\oint_{ \R \lambda= \frac{-\theta_1}{\theta_1+\theta_2}
(|\Txi|^2-|\Im \lambda|^2) \ge -\theta_1}
e^{i\tilde \xi \cdot \tilde x + \lambda t}
(\lambda- L_\Txi)^{-1}\hat f(x_1, \Txi) \, d\tilde \xi d\lambda,
\end{aligned}
\end{equation}
where $\theta_j>0$ are as in \eqref{S1def},
and Corollary \ref{decomposition}, with $\theta_1$ sufficiently
small in relation to $\theta_2$.
Bounding
\begin{align}
|\tilde f|_{L^\infty(\xi', L^1(x_1))}\le
|f|_{L^1(x_1,\tilde x)}=|f|_1
\end{align}
using Hausdorff--Young's inequality, and appealing to the $L^1\to L^p$
resolvent estimates of Proposition \ref{mresbounds} with
$\theta:= \frac{-\theta_1}{\theta_1+\theta_2}>0$ taken suffiently small,
we may thus bound 
\begin{align}
|\hat u(x_1,\xi',\lambda)|_{L^2(x_1)}\le |f|_1 
b(\Txi,\lambda),
\end{align}
for $\hat u:=(L_\Txi -\lambda)f$ and 
$\rho:=|\Txi|+|\lambda|$ sufficiently small,
where 
\begin{equation}
\label{bdef}
b(\Txi, \lambda):=\rho^{-1} \gamma_2(\Txi,\lambda).
\end{equation}

Denoting by $\Gamma(\Txi)$ the arc
\begin{equation}
\label{Gammadef}
\R \lambda= \frac{-\theta_1}{\theta_1+\theta_2}
(|\Txi|^2+|\Im \lambda|^2) \ge -\theta_1
\end{equation}
and using in turn Parseval's identity, Fubini's Theorem,
the triangle inequality, and our $L^1 \to L^2$ resolvent
bounds, we may estimate
\begin{align}
\begin{split}
|S_1(t)f|_{L^2(x_1,\tilde x)}(t)&=
\Big(
\int_{x_1}
\int_{\Txi\in \RR^{d-1}} 
\Big| 
\oint_{\lambda\in \tilde \Gamma(\Txi)}
e^{\lambda t} \hat u(x_1,\Txi,\lambda) d\lambda
\Big|^2
d\Txi \, dx_1
\Big)^{1/2}\\
&=
\Big(
\int_{\Txi\in \RR^{d-1}} 
\Big| 
\oint_{\lambda\in \tilde \Gamma(\Txi)}
e^{\lambda t} \hat u(x_1,\Txi,\lambda) d\lambda
\Big|_{L^2(x_1)}^2 
d\Txi
\Big)^{1/2}\\
&\leq
\Big(
\int_{\Txi\in \RR^{d-1}} 
\Big| 
\oint_{\lambda\in \tilde \Gamma(\Txi)}
|e^{\lambda t}| 
|\hat u(x_1,\Txi,\lambda)|_{L^2(x_1)} d\lambda
\Big|^2 
d\Txi
\Big)^{1/2}\\
\leq
&|f|_1
\Big(
\int_{\Txi\in \RR^{d-1}} 
\Big| 
\oint_{\lambda\in \tilde \Gamma(\Txi)}
e^{\Re \lambda t} 
b(\Txi,\lambda) d\lambda
\Big|^2 
d\Txi
\Big)^{1/2},\\
\end{split}
\end{align}
from which we readily obtain the claimed bound on $|S_1(t) f|_2$
using \eqref{bdef} and the definition of $\gamma_2$.
Specifically , parametrizing $\Gamma(\Txi)$ by
$$
\lambda(\Txi,k)= ik - \theta_1(k^2 + |\Txi|^2), \qquad k\in \RR,
$$
and observing that in nonpolar coordinates
\begin{align}
\begin{split}
\rho^{-1}\gamma_2
&\le \Big[(|k|+ |\Txi|)^{-1}(1+ \sum_{j\ge 1}
\Big(\frac{|k- \tau_j(\Txi)|}{\rho}\Big)^{\frac{1}{s_j}-1}
\Big]\\
&\le \Big[|k|+ |\Txi|)^{-1}(1+ \sum_{j\ge 1}
\Big(\frac{|k- \tau_j(\Txi)|}{\rho}\Big)^{\varepsilon-1}
\Big]
,\\
\end{split}
\end{align}
where $\varepsilon:= \frac{1}{\max_j s_j}$
($0<\varepsilon<1$ chosen arbitrarily if there are no singularities),
we obtain a contribution bounded by
\begin{align}
\begin{split}
&C|f|_1
\Big(
\int_{\Txi\in \RR^{d-1}} 
\Big| 
\int_{-1}^{+1}
e^{-\theta(k^2 + |\Txi|^2) t} 
(\rho)^{-1}\gamma_2 dk
\Big|^2 
d\Txi
\Big)^{1/2}\\
&\leq C
|f|_1
\int_{\Txi\in \RR^{d-1}} 
\Big(e^{-2\theta |\Txi|^2t}
|\Txi|^{-2 \varepsilon}
\Big| 
\int_{-\infty}^{+\infty}
e^{-\theta |k|^2  t} 
|k|^{\varepsilon -1} dk
\Big|^2 
d\Txi
\Big)^{1/2}\\
+ C\sum_{j\ge 1}
|f|_1
&\int_{\Txi\in \RR^{d-1}} 
\Big(e^{-2\theta |\Txi|^2t}
|\Txi|^{-2 \varepsilon}
\Big| 
\int_{-\infty}^{+\infty}
e^{-\theta |k|^2  t} 
|k-\tau_j(\Txi)|^{\varepsilon -1} dk
\Big|^2 
d\Txi
\Big)^{1/2}\\
&\leq C
|f|_1
\int_{\Txi\in \RR^{d-1}} 
\Big(e^{-2\theta |\Txi|^2t}
|\Txi|^{-2 \varepsilon}
\Big| 
\int_{-\infty}^{+\infty}
e^{-\theta |k|^2  t} 
|k|^{\varepsilon -1} dk
\Big|^2 
d\Txi
\Big)^{1/2}\\
&\leq C |f|_1 t^{-(d-1)/4},
\end{split}
\end{align}
yielding the asserted bound for $t\ge 1$; for $t\le 1$ 
on the other hand, we obtain the asserted uniform bound
by local integrability together with boundedness of
the region of integration.
Derivative bounds, $\beta=1$, follow similarly.
\end{proof}

\begin{rem}
\label{deformed}
\textup{
In the proof of Proposition \ref{S1}, we have used
in a fundamental way the semigroup representation and
the autonomy of the underlying equations,
specifically in the parabolic deformation of contours revealing
the stabilizing effect of diffusion.
The resulting bounds are not available through Parseval's
identity or Hausdorff--Young's inequality, corresponding to the choice of 
straight-line contours parallel to the imaginary axis.
Such ``parabolic,'' semigroup-type estimates do not immediately 
translate to the small-viscosity context, 
for which the linearized equations are variable-coefficient in time.
Moreover, it is not clear in this context how one could obtain
resolvent bounds between different norms, since standard
pseudo- or paradifferential techniques are based on 
Fourier decompositions respecting $L^2$ but not necessarily
other norms.
The efficient accounting of diffusive effects in this setting is an
interesting and fundamental problem that appears to be of 
general mathematical interest [M\'eZ.1].
}
\end{rem}

\begin{rem}
\label{newboltz}
\textup{
As noted in the acknowledgements, 
all estimates on the linearized solution operator,
both on $S_1(t)$ and $S_2(t)$,
have been obtained through resolvent bounds alone.
These could be obtained in principle by a variety
of methods, 
an observation that may be useful in the more general
situations (e.g., semidiscrete or Boltzmann shock profiles)
described in Remark \ref{boltz}.4.
}
\end{rem}
\medbreak
{\bf 5.2. Link to the inviscid case.}
It remains to study the resolvent (resp. eigenvalue) equation
\begin{equation}
\label{resevalue}
(L_\Txi-\lambda )U=\begin{cases} f,\\ 0, \end{cases}
\end{equation}
where
\begin{equation}
\label{resevaluedef}
\begin{aligned}
(L_\Txi-\lambda )U&=
\overbrace{ (B^{11}U')'-(A^1U)'}^{L_0U} 
 -i \sum_{j\not= 1}A^j \xi_j U \\
&\quad + i\sum_{j\not= 1} B^{j1}\xi_j U' 
 + i\sum_{k\not= 1}(B^{1k}\xi_k U)' 
-\sum_{j,k\not= 1} B^{jk}\xi_j \xi_k U - \lambda U.
\end{aligned}
\end{equation}
Up to this point our analysis has been rather general;
from here on, we make extensive use of the property,
convenient for the application of both Evans function 
and Kreiss symmetrizer techniques,
that the \eqref{resevalue} may be written as an ODE
\begin{equation}
\label{ODE5}
W'- \mA(\Txi,\lambda, x_1)W= 
\begin{cases}  F,\\ 0, \end{cases}
\end{equation}
in the phase variable 
$W:=\begin{pmatrix}U\\B^{11}_{II}U'\end{pmatrix}\in \CC^{n+r}$,
where $F\sim f$.
This is a straightforward consequence of block-diagonal structure,
(A1), and invertibility of $\tilde A^1_{11}$, (H1); 
see [MZ.3, Z.4] for further details.
(There is ample reason to relax this condition, however; see
Remark \ref{boltz}.4.)
Our first task, carried out in this subsection,
is to make contact with the inviscid case.

\medbreak
{\bf 5.2.1. Normal modes.}  By the conjugation lemma, Lemma
\ref{conjugation}, we may reduce \eqref{resevaluedef} 
by a change of coordinates $W=P_\pm Z_\pm$, $f=P_\pm \tilde f$, with 
$P_\pm \to I$ as $x_1\to \pm \infty$,  $P$ analytic in $(\Txi, \lambda)$,
to a pair of constant-coefficient equations
\begin{equation}
\label{cc5}
Z_\pm'-\mA_\pm(\Txi,\lambda)Z_\pm= 
\begin{cases}
\tilde f\\0,
\end{cases}
\end{equation}
on the half-lines $x_1 \gtrless 0$, coupled by the implicitly determined
transmission conditions 
\begin{equation}
\label{bc}
P_+Z_+(0)- P_-Z_-(0)=0 
\end{equation}
at the boundary $x_1=0$,  where
\begin{equation}
\mA_\pm(\Txi,\lambda):=\lim_{x_1\to \pm \infty} \mA(\Txi,\lambda,x_1):
\end{equation}
that is, a system that at least superficially resembles that
arising in the inviscid theory.
The following lemma generalizes a standard result of Hersch [H]
in the inviscid theory; see Exercise \ref{hersch}.

\begin{lem}[{[ZS, Z.3]}]\label{vhersch}
Assuming (A1)--(A2), (H1) (and, implicitly, existence of a profile),
the matrices $\mA_\pm$ have no center subspace on 
\begin{equation} \label{partialdiss}
\Lambda:=\{\lambda:\,
\R \lambda > \frac{-\theta (|\Txi|^2+ |\Im \lambda|^2)}
{1+ |\Txi|^2+ |\Im \lambda|^2}\},
\end{equation}
for $\theta>0$ sufficiently small;
moreover, the dimensions of their stable and unstable subspaces
agree, summing to full rank $n+r$.
\end{lem}

\begin{proof}
The fundamental modes of \eqref{cc5} are of form $e^{\mu x}V$,
where $\mu$, $V$ satisfy the {\it characteristic equation}
\begin{equation}
\label{char2}
\begin{aligned}
&\left[\mu^2 B^{11}_\pm +\mu(-A^1_\pm + i\sum_{j\not=
1}B^{j1}_\pm \xi_j
+i\sum_{k\not= 1}B^{1k} \xi_k)
\right.\\
&\qquad \left.-(i\sum_{j\not= 1}A^j\xi_j
+ \sum_{jk\not= 1}B^{jk}\xi_j \xi_k + \lambda
I)\right]\V=0.
\end{aligned} 
\end{equation}
The existence of a center manifold thus corresponds with existence of
solutions $\mu=i\xi_1$, $V$ of \eqref{char2}, $\xi_1$ real, i.e., solutions
of the {\it dispersion relation}
\begin{equation}
\label{disp}
(-\sum_{j,k}B^{jk}_\pm\xi_j \xi_k-i\sum_j A^j_\pm \xi_j - \lambda I)\V=0.
\end{equation}
But, $\lambda\in \sigma (-B^{\xi \xi} -iA^{\xi}_\pm )$ implies, by 
\eqref{multidiss}, that 
$$
\R \lambda \le -\theta_1 |\xi|^2/(1+ |\xi|^2).
$$
Noting that for low frequencies $|\Im \lambda |=\CalO(|\xi|)$,
we thus obtain 
$$
\R \lambda \le \frac{-\theta (|\Txi|^2+ |\Im \lambda|^2)}
{1+ |\Txi|^2+ |\Im \lambda|^2},
$$
in contradiction of \eqref{partialdiss}.

Nonexistence of a center manifold, together with 
connectivity of the set described in \eqref{partialdiss}, 
implies that the dimensions of stable and unstable manifolds 
at $+\infty$/$-\infty$ are constant on $\Lambda$.
Taking $\lambda\to +\infty$ along the real axis, with $\Txi\equiv 0$,
we find that these dimensions sum to the full dimension $n+r$ as claimed.
For, Fourier expansion about $\xi_1 =\infty$ of the one-dimensional
($\Txi=0$) dispersion relation (see, e.g.,  Appendix A.4 [Z.4]) yields 
$n-r$ ``hyperbolic'' modes
$$
\lambda_j = -i\xi_1 a_j^* + \dots, \quad j=1, \dots, n-r,
$$
where $a_j^*$ denote the eigenvalues of $A^1_*:=A^1_{11}-(b^{11}_2)^{-1}
b^{11}_1 A^1_{12}$
and $r$ ``parabolic'' modes
$$
\lambda_{n-r +j} = -b_j\xi_1^2 + \dots, \quad j=1, \dots r,
$$
where $b_j$ denote the eigenvalues of $b^{11}_2$; here, we have
suppressed the $\pm$ indices for readability.
Inverting these relationships to solve for $\mu:=i\xi_1$, we find,
for $\lambda \to \infty$, that there are $n-r$ hyperbolic
roots $\mu_j\sim -\lambda/a_j^*$, and $r$ parabolic roots
$\mu_{n-r+j}^\pm \sim \sqrt{\lambda/b_j}$.
By assumption (H1), $\det A^1_*(x_1) \ne 0$ for all $x_1$,
and so the former yield a fixed
number $k$/$(n-r-k)$ of stable/unstable roots, independent of $x_1$, and
thus of $\pm$.  Likewise, ($\tilde{\text\rm H1}$(i) implies that
the latter yields $r$ stable, $r$ unstable roots.
(Note: we have here used the existence of a connecting profile: i.e.,
equality is not a consequence of algebraic structure alone.)
Combining, we find the desired consistent splitting, with $(k+r)$/ $(n-k)$
stable/unstable roots at both $\pm \infty$.  
\end{proof}

\begin{cor}\label{normalmodes}
Assuming (A1)--(A2) and (H1), 
on the set $\Lambda$ defined in \eqref{partialdiss},
there exists an analytic choice of basis vectors
$V_j^\pm$, $j=1, \dots k$ (resp. $V_j^\pm$, $j=k+1,\dots, n+r$) 
spanning the stable (resp. unstable) subspace of $\mA_\pm$.
\end{cor}

\begin{proof}
By spectral separation of stable and unstable subspaces, the associated
(group) eigenprojections are analytic.
The existence of analytic bases
then follows by a standard result of Kato; see [Kat], pp. 99--102.
\end{proof}

\begin{defi}\label{normal}
\textup{
For $x_1\gtrless 0$, define {\it normal modes} for the variable-coefficient
ODE \eqref{ODE5} as
\begin{equation}\label{wmodes}
W_j^\pm(x_1):= P_\pm e^{\mA_\pm x_1}V_j^\pm.
\end{equation}
}
\end{defi}

\medbreak
In particular, $W_1^+, \dots, W_k^+$ span the manifold of solutions
of the eigenvalue equation decaying as $x_1\to +\infty$, and
$W_{k+1}^+, \dots, W_{n+r}^+$ span the manifold of solutions
of the eigenvalue equation decaying as $x_1\to -\infty$.
\medbreak
{\bf 5.2.2. Low-frequency asymptotics.}  
For comparison with the inviscid case, it is convenient to 
introduce polar coordinates 
\begin{equation}
\label{polar}
(\Txi, \lambda)=:(\rho \Txi_0, \rho\lambda_0),
\end{equation}
$\rho \in \RR^1$, $(\Txi_0, \lambda_0) \in \RR^{d-1}\times 
\{\R \lambda \ge 0\} \setminus \{(0,0)\}$,
and consider $W_j^\pm$, $V_j^\pm$ as functions of $(\rho, \Txi_0, \lambda_0)$.
With this notation, the resolvent (resp. eigenvalue) equation
\eqref{resevaluedef} may be viewed as a singular perturbation
as $\rho \to 0$ of the corresponding inviscid equation.
The following results quantify this observation, separating
normal modes into slow-decaying modes asymptotic to those of 
the inviscid theory (the ``slow manifold'' in the singular perturbation
limit)
and fast-decaying transient modes associated with the viscous regularization.


\begin{lem}[{[K, M\'e.4]}]\label{3.4}  
Let there hold (A1), (H2), and (H4), or, more generally, 
$\sigma(df^\xi(u_\pm))$ real, semisimple, and of constant multiplicity
for $\xi \in \RR^d\setminus \{0\}$ and $\det df^1\ne 0$.   
Then, the vectors $\{r^-_1,\cdots, r_{n-i_-}^-\}, \{r^+_{i_+ +1},\cdots,
r^+_n\}$ defined as in \eqref{Delta}--\eqref{CalAhyp}, Section 2.2
(and thus also $\Delta (\Txi,\lambda )$)
may be chosen to be homogeneous degree zero (resp. one), 
analytic on $\Txi\in \RR^{d-1}$, $\R \lambda>0$
and continuous at the boundary $\Txi \in \RR^{d-1}\setminus \{0\}$,
$\R \lambda =0$.
\end{lem}

\begin{proof}
See Exercises 4.23--4.24 and
Remark 4.25, Section 4.5.2 of [Z.3]
for a proof in the general case
(three alternative proofs, based respectively
on [K], [CP], and [ZS]).
In the case of main interest, when (H5) holds
as well, this will be established through
the explicit computations of Section 5.4.
\end{proof}

\begin{lem} [{[ZS, Z.3, M\'eZ.2]}] \label{vcont} 
Under assumptions (A1)--(A2) and (H0)--(H4), the functions
$V_j^\pm$ may be chosen within groups of $r$ ``fast'', or
``viscous'' modes bounded away from the center subspace 
of coefficient $\mA_\pm$, analytic in $(\rho, \Txi_0, \lambda_0)$ 
for $\rho \ge 0$, $\Txi_0\in \RR^{d-1}$, $\R \lambda_0\ge 0$,
and $n$ ``slow'', or ``inviscid''  modes approaching the center 
subspace as $\rho \to 0$, analytic in $(\rho, \Txi_0, \lambda_0)$
for $\rho>0$, $\Txi_0\in \RR^{d-1}$, $\R \lambda_0\ge 0$ and 
continuous at the boundary $\rho=0$, with limits
\begin{equation}
\label{Vr}
V_j^\pm(0, \Txi_0, \lambda_0)=
\begin{pmatrix} (A^1_\pm)^{-1}r_j^\pm (\Txi_0, \lambda_0)\\0\end{pmatrix},
\end{equation}
$r_j^\pm$ defined as in \eqref{Delta}--\eqref{CalAhyp}.
\end{lem}

\begin{proof}
We carry out here the simpler case $\R \lambda_0>0$.
In the case of interest, that (H5) also holds, 
the case $\R \lambda = 0$ follows by the 
detailed computations in Section 5.4; see Remark \ref{limitrem}.
For the general case (without (H5)), see [M\'e.2].  

Substituting $U=e^{\mu x_1}\V$ into the limiting
eigenvalue equations 
written in polar coordinates, we obtain the polar characteristic equation,
\begin{equation}
\label{2.17}
\begin{aligned}
&\Big[\mu^2 B^{11}_\pm +\mu(-A^1_\pm + i\rho \sum_{j\not=
1}B^{j1}_\pm \xi_j
+i\rho \sum_{k\not= 1}B^{1k} \xi_k) \\
&-(i\rho \sum_{j\not= 1}A^j\xi_j
+\rho ^2 \sum_{jk\not= 1}B^{jk}\xi_j \xi_k + \rho \lambda
I)\Big]\V=0,
\end{aligned} 
\end{equation}
where for notational convenience we have dropped subscripts 
from the fixed parameters $\Txi_0$, $\lambda_0$.
At $\rho=0$, this simplifies to
$$
\left(\mu^2 B^{11}_\pm -\mu A^1_\pm  \right)\V=0,
$$
which, by the analysis in Appendix A.2
of the linearized traveling-wave ordinary differential equation
$\left(\mu B^{11}_\pm - A^1_\pm  \right)\V=0$,
has $n$ roots $\mu=0$, and $r$ roots $\R \mu \ne 0$.
The latter, ``fast'' roots correspond to stable and unstable
subspaces,  which extend analytically as claimed by their
spectral separation from other modes; thus, we need only focus on the
bifurcation as $\rho $ varies near zero
of the $n$-dimensional center manifold associated with ``slow'' roots $\mu=0$. 

Positing a first-order Taylor expansion 
\begin{equation}
\label{2.18}
\begin{cases}
\mu=0+\mu^1 \rho +o(\rho),\\
\V=\V^0 + \V^1 \rho + o(\rho),
\end{cases} 
\end{equation}
and matching terms of order $\rho $ in \eqref{2.17}, we obtain
\begin{equation}
\label{2.19}
(-\mu^1 A^1_\pm - i\sum_{j\not= 1}A^j\xi_j - \lambda I)\V^0=0,
\end{equation}
or equivalently $-\mu_1$ is an eigenvalue of $(A^1)^{-1}(\lambda + iA^\Txi) $
with associated eigenvector $\V^0$.

For $\R \lambda>0$, $(A^1)^{-1}(\lambda + iA^\Txi)$ has no center
subspace. 
For, substituting $\mu^1=i\xi_1$ in \eqref{2.19}, we obtain $\lambda \in
\sigma(iA^\xi)$, pure imaginary, a contradiction.  Thus, the
stable/unstable spectrum {\it splits} to first order, and we 
obtain the desired analytic extension by standard matrix
perturbation theory, though not in fact the analyticity of 
individual  eigenvalues $\mu$.
\end{proof}
\medbreak

\begin{rem}
\label{bif1}
\textup{
The first-order approximation \eqref{2.19} is exactly the 
matrix perturbation problem arising in the inviscid theory [K, M\'e.1].
}
\end{rem}

\begin{cor}\label{2.5}  
Denoting $W=(U,b^{11}_1u^I+b^{11}_2u^{II})$ as above, we may
arrange at $\rho =0$ that all $W^\pm_j$ satisfy the linearized 
traveling-wave ODE 
$$
(B^{11}U')'-(A^1U)'=0,
$$
with constant of integration
\begin{equation}
\label{Weq2}
B^{11} U^{\pm'}_j-A^1U_j^\pm\equiv 
\begin{cases} 0, \quad \text{\rm for fast modes,}\\
r_j^\pm, \quad \text{\rm for slow modes,}\\
\end{cases}
\end{equation}
$r^\pm_j$ as above, with fast modes analytic at the $\rho=0$ boundary
and independent of $(\Txi_0,\lambda_0)$ for $\rho=0$,
and slow modes continuous at $\rho=0$.
\end{cor}

\begin{proof}  Immediate \end{proof}

\medbreak

{\bf 5.2.4. The Evans function and its low-frequency limit.} 
We now complete the analogy with the inviscid case, 
introducing the Evans function $D$ and establishing its
relation to the Lopatinski determinant $\Delta$ 
in the limit as frequency goes to zero, Proposition \ref{ZS}. 

\begin{defi}\label{evans}
\textup{
Following the standard construction of, e.g., 
[E.1--5, AGJ, PW, GZ, ZS], we define on the set $\Lambda$ 
an {\it Evans function}
\begin{equation}
\label{evanseq}
\begin{aligned}
D(\Txi, \lambda)&:=
\det
\Big(
 W_{1}^+, \dots,  W_k^+,
 W_{k+1}^-, \dots, W_{N}^-
\Big)_{|x=0, \lambda} \\
&=
\det
\Big(
P_+ V_{1}^+, \dots, P_+ V_k^+,
P_- V_{k+1}^-, \dots, P_- V_{N}^-
\Big)_{|x=0, \lambda},
\end{aligned}
\end{equation}
measuring the (solid) angle of intersection
between the manifolds of solutions of the eigenvalue
equation decaying as $x_1\to +\infty$ and $x_1\to -\infty$, respectively,
where $P_\pm$, $V_j^\pm$, $W_j^\pm$ are as in Definition \ref{normal}.
}
\end{defi}

\medbreak
Evidently, $D$ vanishes if and only if there exists an exponentially
decaying solution of the eigenvalue equation, i.e., $\lambda$ is
an eigenvalue of $L_\Txi$, or, equivalently, there exists a solution
of the boundary-value problem \eqref{cc5} satisfying boundary
condition \eqref{bc}.
That is, $D$ is the analog for the viscous problem of the Lopatinski
determinant $\Delta$ for the inviscid one.

%

\medbreak

\begin{proof}[Proof of Proposition \ref{ZS}]
We carry out the proof in the Lax case only.  The proofs 
in the under- and overcompressive cases are quite similar;
see [ZS].
We are free to make any analytic choice of bases, and any
nonsingular choice of coordinates, since these affect the
Evans function only up to a nonvanishing analytic multiplier
which does not affect the result.
Choose bases $W_j^\pm$ as in Lemma \ref{2.5}, 
$W=(U, z_2')$, $z_2=b^{11}_1u^{I}+b^{11}_2u^{II}$.
Noting that $L_0\bU'=0$, by translation invariance, we have that
$\bU'$ lies in both $\Span\{U^+_1,\cdots, U^+_{K}\}$ and
$\Span(U^-_{K+1},\cdots, U^-_{n+r})$ for $\rho =0$, 
hence without loss of generality 
\begin{equation}
\label{2.22b}
U^+_1=U^-_{n+r}=\bU',
\end{equation}
independent of $\Txi$, $\lambda$.
(Here, as usual, $``\, '\, $'' denotes $\partial/\partial x_1)$.

More generally, we order the bases so that
$$
W_1^+, \dots, W_k^+ \hbox{ and  }
W_{n+r-k-1}^-, \dots, W_{n+r}^-
$$
are fast modes (decaying for $\rho=0$) and
$$
W_{k+1}^+, \dots, W_K^+ \hbox{ and }
W_{K+1}^-, \dots, W_{n+r-k-2}^-
$$
are slow modes (asymptotically constant for $\rho=0$),
fast modes analytic and slow modes continuous at $\rho=0$
(Corollary \ref{2.5}).

Using the fact that $U_1^+$ and $U_{n+r}^-$ are analytic,
we may express 
\begin{equation}
\begin{aligned}
U_1^+(\rho)&=U_1^+(0)+ U_{1,\rho}^+\rho + o(\rho),\\
U_{n+r}^-(\rho)&=U_{n+r}^-(0)+ U_{{n+r},\rho}^-\rho + o(\rho),
\end{aligned}
\end{equation}
Writing out the eigenvalue equation
\begin{equation}
\label{2.25}
\begin{aligned}
(B^{11}w')&=(A^1w')-i\rho \sum_{j\not= 1} B^{j1}\xi_j w'\\
&\quad -i\rho (\sum_{k\not= 1}B^{1k}\xi_k w)' + i\rho \sum_{j\not= 1} A^j
\xi_j w+\rho \lambda w\\
&\quad -\rho ^2 \sum_{j,k\not= 1} B^{jk}\xi_j \xi_k w,
\end{aligned}
\end{equation}
in polar coordinates, we find that $Y^+:= U_{1\rho}^+(0)$ and
$Y^-:= U_{n+r,\rho}^-(0)$ satisfy the variational equations 
\begin{equation}
\label{2.26}
\begin{aligned}
(B^{11}Y')'&=(A^1Y')-\Big[i\sum_{j\not= 1}B^{j1}\xi_j \bU'\\
&\quad + i\Big(\sum_{k\not= 1}B^{1k}\xi_k \bU'\Big) +i
\sum_{j \not= 1} A^j \xi_j \bU'+\lambda \bU'\Big],\\
\end{aligned}
\end{equation}
with boundary conditions $Y^+(+\infty)= Y^-(-\infty)=0$.
Integrating from $+\infty$, $- \infty$ respectively, we obtain therefore
\begin{equation}
\label{2.27}
\begin{aligned}
B^{11}Y{\pm'}-A^1Y^{\pm'}&=if^\Txi(\bU)-iB^{1\Txi}(\bU)\bU'\\
&\quad -iB^{\Txi 1}(\bU)\bU'+\lambda \bU-\Big[if^\Txi (u_\pm)+\lambda
u_\pm\Big], 
\end{aligned}
\end{equation}
hence $\tilde Y:=(Y^- -Y^+)$ satisfies
\begin{equation}
\label{Yeq}
B^{11}\tilde Y'-A^1\tilde Y=i[f^\Txi]+\lambda [u].
\end{equation}

By (A1) together with (H1), 
$\begin{pmatrix} A^1_{11}& A^1_{12}\\b^{11}_1&b^{11}_2\end{pmatrix} $
is invertible, hence 
\begin{equation}
\label{coord}
\begin{aligned}
(U,z_2')&\to (z_2, -z_1, z_2'+ (A^1_{21}-b^{11'}_1,A^1_{22}-b^{11'}_2)U)\\
&\quad =(z_2, B^{11}U'-AU)
\end{aligned}
\end{equation}
is a nonsingular coordinate change, where 
$
\begin{pmatrix} z_1\\z_2\end{pmatrix}:=
\begin{pmatrix} A^1_{11}& A^1_{12}\\b^{11}_1&b^{11}_2\end{pmatrix} U.
$

Fixing $\Txi_0$, $\lambda_0$, and using $W_1^+(0)=W_{n+r}^-(0)$,
we have
$$
\begin{aligned}
D(\rho)
&= \det 
\Big(
W_1^+(0)+ \rho W_{1_\rho }^+(0)+o(\rho),\, \cdots,W_{K}^+(0)+o(1),\\
&\qquad \qquad \qquad
\, W_{K+1}^-(0)+o(1),\, \cdots, \, W^-_{n+r}(0) \rho W_{{n+r}_\rho }^-(0)+o(\rho)
\Big)\\
&= \det 
\Big(
W_1^+(0)+ \rho W_{1_\rho }^+(0)+o(\rho),\, \cdots,W_{K}^+(0)+o(1),\\
&\qquad \qquad \qquad
\, W_{K+1}^-(0)+o(1),\, \cdots, \, \rho \tilde Y(0)+o(\rho)
\Big)\\
&= \det 
\Big(
W_1^+(0) \, \cdots,W_{K}^+(0)
\, W_{K+1}^-(0)\, \cdots, \, \rho \tilde Y(0)
\Big)  + o(\rho)\\
\end{aligned}
$$

Applying now \eqref{coord} and using \eqref{Weq2} and \eqref{Yeq}, we obtain
\setbox0=\hbox{$r^+_{i_+ +1},\cdots, r^+_n,r^-_1, \cdots, r^-_{n-i_-}$}

\setbox1=\hbox {$i[f^\Txi(u)]+\lambda [u]$}

\setbox2=\vbox{\hsize.46\hsize
$$\begin{aligned}
&\hbox to \wd1{\hfil
$\overbrace{z^+_{2,1} \cdots, z^+_{2,k}}^{\hbox{fast}},$\hfil} 
\hbox to \wd0{\hfil 
$\overbrace{*, \cdots, *,*,\cdots,*}^{\hbox{slow}},$\hfil } \\
&\hbox to \wd1{\hfil $0,\cdots,0,$ \hfil }
\hbox to \wd0{ \hfil $ 
r^+_{i_+ +1},\cdots, r^+_n,r^-_1, \cdots, r^-_{n-i_-}
$\hfil } \\
\end{aligned}
         $$
}

\setbox3=\vbox{\hsize.76\hsize
$$\begin{aligned}
& \hbox to \wd0{\hfil$
\overbrace{z^-_{2,n+r-k-1}, \cdots, z^-_{2, n+r-1}}^{\hbox{fast}},
$\hfil } 
\hbox to \wd1{\hfil $* $\hfil } \\
& \hbox to \wd0{\hfil$0,\cdots,0,$\hfil } 
\hbox to \wd1{\hfil $i[f^\Txi(u)]+\lambda [u]$\hfil } \\
\end{aligned}  
	 $$
}

$$\displaylines{
D(\rho)= C
\det \left(\vcenter{\box2}\right.\hfill\cr
\hfill
\left.\vcenter{\box3}\right)_{|_{x_1=0}} + o(\rho)\cr
=\gamma  \Delta (\Txi,\lambda )+ o(\rho) \qquad \qquad \qquad 
\qquad \qquad \qquad \qquad
\hfil
}
$$
as claimed, where
$$
\gamma:= C
\det \left( z^+_{2,1}, \cdots, z^+_{2,k},  z^-_{2,n+r-k-1},\cdots,
z^-_{2,n+r-1} \right)_{x_1=0}.
$$
Noting that $\{z^+_{2,1},\cdots, z^+_{2,k}\}$ and
$\{z^-_{2,n+r-k-1},\cdots, z^-_{2,n+r}\}$ span the tangent manifolds
at $\bU(\cdot)$ of the stable/unstable manifolds of traveling wave
ODE \eqref{II} at
$U_+/U_-$, respectively, with $z^+_{2,1}=z^-_{2,n+r}
=(b^{11}_1,b^{11}_2)\bU'$ 
in common, we see that $\gamma $ indeed measures transversality 
of their intersection;  
moreover, $\gamma$ is constant, by Corollary \ref{2.5}.
\end{proof}

\begin {rem}
\textup{
The proof of Proposition \ref{ZS} may be recognized as a generalization
of the basic Evans function calculation pioneered by Evans [E.4], 
relating behavior near the origin to geometry of the phase space 
of the traveling wave ODE and thus giving an explicit link between PDE
and ODE dynamics.
The corresponding one-dimensional result was established in [GZ];
for related calculations, see, e.g., [J, AGJ, PW].
}
\end{rem}

\medbreak
{\bf 5.3.  Spectral bounds and necessary conditions for stability.}
Using Proposition \ref{ZS}, we readily obtain the
stated necessary conditions for stability, Theorem \ref{mainnec}.
Define the {\it reduced Evans function} as
\begin{equation}
\label{3.1}
\bDelta(\Txi,\lambda ):= \lim_{\rho \to 0} \rho ^{-\ell} D(\rho
\Txi,\rho \lambda ).
\end{equation}
By the results of the previous section, the limit $\bDelta$ exists
and is analytic, with
\begin{equation}
\label{3.2}
\bDelta = \gamma \Delta,
\end{equation}
for shocks of pure type (indeed, such a limit exists for all
types).
Evidently, $\bDelta (\cdot , \cdot)$ is homogeneous, degree $\ell$.\footnote{
Here, and elsewhere, homogeneity is with respect to the positive reals,
as in most cases should be clear from the context. 
Recall that $\Delta$ (and thus $\bar \Delta$) 
is only defined for real $\Txi$, and $\R\lambda\ge 0$.
}

\begin{lem}[{[ZS]}]\label{tangency}  Let $\bDelta(0,1)\not= 0$.
Then, near any root $(\Txi_0,\lambda _0)$ of 
$\bDelta (\cdot , \cdot)$, there exists a continuous branch
$\lambda (\Txi)$, homogeneous degree one, of solutions of
\begin{equation}
\label{3.5}
\bDelta (\Txi,\lambda (\Txi))\equiv 0
\end{equation}
defined in a neighborhood $V$ of $\Txi_0$, with
$\lambda(\Txi_0)=\lambda_0$.  Likewise, there exists a continuous branch
$\lambda_*(\Txi)$ of roots of
\begin{equation}
\label{3.6}
D(\Txi, \lambda_* (\Txi))\equiv 0,
\end{equation}
defined on a conical neighborhood 
$V_{\rho_0}:=\{\Txi=\rho \tilde \eta:\tilde \eta \in V, \, 0<\rho <\rho_0 \}$, 
$\rho _0 >0$ sufficiently small,  ``tangent'' to $\lambda(\cdot)$ 
in the sense that
\begin{equation}
\label{3.6.1}
|\lambda _*(\Txi)-\lambda (\Txi)|=o(|\Txi|)
\end{equation}
as $|\Txi|\to 0$, for $\Txi\in V_{\rho_0}$.
\end{lem}

\begin{proof}
Provided that $\bDelta (\Txi_0,\cdot)\not \equiv 0$, 
the statement \eqref{3.5} follows by Rouche's Theorem,
since $\bDelta (\Txi_0,
\cdot)$ are a continuous family of analytic functions.
But, otherwise, restricting $\lambda$ to the positive real axis, we have
by homogeneity that
$$
\begin{aligned}
0 = \lim_{\lambda \to +\infty}\bDelta (\Txi_0,\lambda )
&= \lim_{\lambda \to +\infty} \bDelta (\Txi_0/\lambda ,1)
= \bDelta (0,1),
\end{aligned}
$$
in contradiction with the hypothesis.
Clearly we can further choose $\lambda(\cdot)$ homogeneous degree one,
by homogeneity of $\bDelta$.
Similar considerations yield existence of a branch of roots
$\bar \lambda
(\Txi,\rho )$ of the family of analytic functions
\begin{equation}
\label{3.7}
g^{ \Txi,\rho} (\lambda ):= \rho ^{-\ell}D(\rho \Txi,\rho \lambda ),
\end{equation}
for $\rho $ sufficiently small, since $g^{\Txi,0}=\bDelta
(\Txi,\cdot)$.
Setting $\lambda _*(\Txi):= |\Txi|\bar \lambda (\Txi/|\Txi|,|\Txi|)$, we have
$$
D(\Txi,\lambda _*(\Txi))=|\Txi|^\ell g^{\Txi/|\Txi|,|\Txi|}(\lambda _*)
\equiv 0,
$$
as claimed.  ``Tangency,'' in the sense of \eqref{3.6.1}, follows
by continuity of $\bar \lambda$ at $\rho=0$, the definition
of $\lambda_*$, and the
fact that 
$\lambda(\Txi)=|\Txi|\lambda(\Txi/|\Txi|)$, by homogeneity of $\bDelta$.
\end{proof}
\medbreak


\begin{proof}[Proof of Theorem \ref{mainnec}]
Weak spectral stability is clearly necessary for viscous stability.
For, unstable ($L^p$) spectrum of $L_\Txi$ (necessarily point
spectrum, by TODO)
corresponds to unstable ($L^p$) essential spectrum of the operator $L$
for $p<\infty$, by a standard limiting argument (see e.g.
[He, Z.1]), and unstable point spectrum for $p=\infty$.
This precludes $L^p\to L^p$ stability, by the generalized
Hille--Yosida theorem, (Proposition \ref{semigpcriteria}, Appendix A). 
Moreover, standard spectral continuity results [Ka, He, Z.2] 
yield that instability, if it occurs, occurs for a band of $\Txi$ values, 
from which we may deduce by inverse Fourier transform the exponential
instability of \eqref{linearized} for test function initial
data $U_0\in C^\infty_0$, with respect to any $L^p$, $1\le p\le\infty$.  

Thus, it is sufficient to establish 
that failure of weak refined dynamical stability implies 
failure of weak spectral stability, i.e.,
existence of a zero $D(\tilde \xi, \lambda)=0$
for $\xi\in \RR^{d-1}$, $\R \lambda >0$.
Failure of weak inviscid stability, or
$\Delta(\Txi,\lambda )= 0$ for $\Txi \in \RR^{d-1}$, $\R\lambda >0$,
implies immediately the existence of such a root, by tangency of
the zero-sets of $D$ and $\Delta$ at the origin, Lemma \ref{tangency}.
Therefore, it remains to consider the case that weak inviscid stability
holds, but there exists a root $D(\xi, i\tau)$ for $\xi$, $\tau$ real,
at which $\Delta$ is analytic, $\Delta_\lambda\ne 0$, and 
$\
 \beta(\xi, i\tau)<0$, where $\beta$ is defined as in \eqref{beta}.

Recalling that $D(\rho\Txi,\rho \lambda)$ vanishes
to order $\ell$ in $\rho$ at $\rho=0$, we find by L'Hopital's rule that
$$
(\partial/\partial \rho )^{\ell+1}D(\rho\Txi,\rho \lambda)
|_{\rho=0,\lambda=i\tau}
=
(1/\ell !)
(\partial/\partial \rho) g^{\Txi,i\tau}(0)
$$
and
$$
(\partial/\partial \lambda\ ) D(\rho\Txi,\rho i\tau)
|_{\rho=0,\lambda=i\tau}
=
(1/\ell !)
(\partial/\partial \lambda) g^{\Txi,i\tau}(0),
$$
where $g^{\Txi,\rho }(\lambda ):=\rho ^{-\ell}D(\rho \Txi,\rho \lambda
)$ as in \eqref{3.7}, whence
\begin{equation}
\label{3.10}
\beta =\frac{(\partial/\partial \rho) g^{\Txi,\lambda}(0)}
{(\partial/\partial \lambda\ ) g^{\Txi,\lambda}(0)},
\end{equation}
for $\beta$ defined as in Definition \ref{refstab},
with $(\partial/\partial \lambda\ ) g^{\Txi,\lambda}(0)\ne0$.
  
By the (analytic) Implicit Function Theorem, therefore,
$\lambda(\Txi,\rho)$ is analytic in $\Txi,\rho$ at $\rho=0$,
with 
\begin{equation}
\label{3.11}
(\partial/\partial \rho)  \ \lambda (\Txi,0)=-\beta,
\end{equation}
where $\lambda_(\Txi,\rho)$ 
as in the proof of Lemma \ref{tangency}
is defined implicitly by $g^{\Txi,\lambda}(\rho)=0$, $\lambda(\Txi,0):=i\tau$. 
We thus have, to first order,
\begin{equation}
\label{3.12}
\lambda (\Txi,\rho )=i\tau-\beta \rho +\CalO(\rho ^2).
\end{equation}

Recalling the definition 
$\lambda _*(\Txi):= |\Txi|\bar \lambda (\Txi/|\Txi|,|\Txi|)$, we have
then, to second order, the series expansion
\begin{equation}
\label{taylor}
\lambda_*(\rho \Txi)= i \rho \tau - \beta \rho^2 + \CalO(\rho ^3),
\end{equation}
where $\lambda_*(\Txi)$ is the root of $D(\Txi,\lambda)=0$ defined
in Lemma \ref{tangency}.
It follows that
there exist unstable roots of $D$ for small $\rho >0$ unless
$\R \beta \ge 0$.  
\end{proof}

\medbreak

{\bf 5.4. Low-frequency resolvent estimates.}
It remains to establish the low-frequency resolvent
bounds of Proposition \ref{mresbounds}.
Accordingly, we restrict attention to arcs
\begin{equation}
\Gamma^\Txi: \,
\R \lambda= \theta (|\Txi|^2+|\Im \lambda|^2),
\quad 
0< |(\Txi, \Im \lambda)|\le \delta,
\end{equation}
with $\theta>0$ and $\delta$ taken sufficiently small.

\medbreak
{\bf 5.4.1. Second-order perturbation problem.}
We begin by deriving
a second-order, viscous correction of the central matrix perturbation 
problem \eqref{2.19} underlying the inviscid stability analysis 
[K, Ma.1--3, M\'e.1].
Introduce the curves
\begin{equation}
\label{gammacurve}
(\Txi,\lambda)(\rho ,\Txi_0,\tau_0):=
\big(\rho \Txi_0,\rho i\tau_0-\theta_1\rho^2 \big),
\end{equation}
where $\Txi_0\in \Bbb{R}^{d-1}$ and $\tau_0\in \Bbb{R}$ are
restricted to the unit sphere $S^d: |\Txi_0|^2+ |\tau_0|^2=1$.
Evidently, as $(\Txi_0,\tau_0,\rho)$ range in the compact
set $S^d\times [0, \delta]$, $(\Txi,\lambda)$ traces out
the surface $\cup_{\Txi} \Gamma^\Txi$ of interest.

Making as usual the Ansatz 
$U=: e^{\mu x_1}\V$,
and substituting $\lambda=i\rho \tau_0 - \theta_1 \rho^2$
into \eqref{2.17}, we obtain the characteristic equation
\begin{equation}
\label{curvechar}
\begin{aligned}
&\left[\mu^2 B^{11}_\pm +\mu(-A^1_\pm + i\rho \sum_{j\not=
1}B^{j1}_\pm \xi_{0_j}
+i\rho \sum_{k\not= 1}B^{1k} \xi_{0_k})
\right.\\
&\left.-(i\rho \sum_{j\not= 1}A^j\xi_{0_j}
+\rho ^2 \sum_{jk\not= 1}B^{jk}\xi_{0_j} \xi_{0_k}
 + (\rho i\tau_0 - \theta_1 \rho^2)
I)\right]\V=0.
\end{aligned} 
\end{equation}
Note that this agrees with \eqref{2.17} up to second order in
$\rho$, hence the (first-order) matrix bifurcation analysis
of Lemma \ref{vcont} applies for any fixed $\theta$.
We focus on ``slow'', or ``inviscid'' modes $\mu \sim \rho$.

Positing the Taylor expansion
\begin{equation}
\label{2.18second}
\begin{cases}
\mu=0+\mu^1 \rho +\cdots,\\
\V=\V^0 + \cdots
\end{cases} 
\end{equation}
(or Puisieux expansion, in the case of a branch singularity) as
before, and matching terms of order $\rho $ in \eqref{curvechar}, we obtain
\begin{equation}
\label{inv5}
(-\mu^1 A^1_\pm - i\sum_{j\not= 1}A^j\xi_{0_j} - i\tau_0 I)\V=0,
\end{equation}
just as in \eqref{2.19}, or equivalently
\begin{equation}
\label{slow}
[(A^1)^{-1}(i\tau_0  +i A^{\Txi_0}) - \alpha_0 I]\V=0,
\end{equation}
with $\mu_1=:-\alpha_0$, which can be recognized as the equation
occurring in the inviscid theory on the imaginary boundary $\lambda=i\tau_0$.

In the inviscid stability theory, solutions of \eqref{inv5}
are subcategorized into ``elliptic'' modes, for which
$\alpha_0$ has nonzero real part, 
``hyperbolic'' modes, for which $\alpha_0$ is pure imaginary
and locally analytic in $(\Txi, \tau)$, and ``glancing''
modes lying on the elliptic--hyperbolic boundary, 
for which $\alpha_0$ is pure imaginary with a branch
singularity at $(\Txi_0, \tau_0)$.
Elliptic modes admit a straightforward treatment, both in
the inviscid and the viscous theory;
however, hyperbolic and glancing modes require a more detailed
matrix perturbation analysis.

Accordingly, we now restrict to the case of a
pure imaginary eigenvalue $\alpha_0=:i\xi_{0_1}$;
here, we must consider quadratic order terms in $\rho$,
and the viscous and inviscid theory part ways.
Using $\mu=-i\rho \xi_{0_1} +o(\rho)$, we obtain at second order
the modified equation:
\begin{equation}
\label{superslow}
[(A^1)^{-1}\big(i\tau_0 + \rho(B^{\xi_0 \xi_0}-\theta_1)  +i A^{\Txi_0}\big) 
- \tilde \alpha I]\tilde \V=0,
\end{equation}
where $\tilde \alpha$ is the next order correction to $\alpha\sim -\mu/\rho$,
and $\tilde \V$ the next order correction to $\V$.
(Note that this derivation remains valid near branch singularities,
since we have only assumed continuity of $\mu/\rho$ and not analyticity
at $\rho=0$).
Here, $B^{\xi_0 \xi_0}$ as usual denotes $\sum B^{jk}\xi_{0_j} \xi_{0_k}$,
where $\xi_0:=(\xi_{0_1},\Txi_0)$.
Equation \eqref{superslow} generalizes the perturbation equation 
\begin{equation}
\label{inviscidpert}
[(A^1)^{-1}(i\tau_0 + \gamma +i A^{\Txi_0}) - \tilde \alpha I]\tilde\V=0,
\end{equation}
$\gamma:=\R \lambda \to 0^+$, 
arising in the inviscid theory near the imaginary boundary $\lambda=i\tau_0$
[K, M\'e.4].

Note that $\tau_0$ is an eigenvalue of $A^{\xi_0}$, as can be seen
by substituting $\alpha_0=i\xi_{0_1}$ in \eqref{slow}, hence
$|\tau_0| \le C |\xi_0|$
and therefore (since clearly also $|\xi_0|\ge |\Txi_0|$)
\begin{equation}
|\xi_0|\ge (1/C)(|(\Txi_0,\tau_0)|= 1/C.
\end{equation}
Thus, for $B$ positive definite, and $\theta_1$ sufficiently small,
perturbation $\rho(B^{\xi_0 \xi_0}-\theta_1)$, roughly speaking,
enters \eqref{superslow}
with the same sign as does $\gamma I$ in \eqref{inviscidpert},
and the same holds true for semidefinite $B$ under the genuine
coupling condition \eqref{gc}; see (K1), Lemma \ref{KSh}.
Indeed, for identity viscosity $B^{jk}:=\delta^j_k I$,
\eqref{superslow} reduces for fixed
$(\Txi_0$, $\tau_0)$ exactly to \eqref{inviscidpert} for
$\theta_1$ sufficiently small,
by the rescaling $\rho \to \rho/(|\xi_0|^2 - \theta_1)$.
This motivates the improved viscous resolvent bounds
of Proposition \ref{mresbounds}; see Remark \ref{cf}.

\medbreak{\bf 5.4.2. Matrix bifurcation analysis.}
The matrix bifurcation analysis for \eqref{superslow} 
goes much as in the inviscid case, but with additional technical
complications due to the presence of the additional parameter $\rho$.
The structural hypothesis (H5), however, allows us to reduce
the calculations somewhat, at the same time imposing additional
structure to be used in the final estimates of Section 5.4.3.

For hyperbolic modes, 
$(\Txi_0,\tau_0)$ bounded away from the set of branch singularities
$\cup (\Txi,\eta_j(\Txi))$, we may treat \eqref{superslow} as a continuous
family of single-variable matrix perturbation problem in $\rho$, indexed
by $(\Txi_0,\tau_0)$; the resulting continuous family of
analytic perturbation series will then yield uniform bounds by compactness.
For glancing modes, $(\Txi_0,\tau_0)$ near a branch singularity, 
on the other hand, we
must vary both $\rho$ and $(\Txi_0,\tau_0)$, in general a complicated
multi-variable perturbation problem.
Using homogeneity, however, and the uniform structure assumed in (H5),
this can be reduced to a two-variable perturbation problem 
that again yields uniform bounds.
For, noting that $\eta_j(\Txi)\equiv 0$,
we find that $\Txi_0$ must be bounded away from the origin at
branch singularities;  thus, we may treat the direction 
$\Txi_0/|\Txi_0|$ as a fixed parameter and vary only $\rho$ and the
ratio $|\tau_0|/|\Txi_0|$.  Alternatively, relaxing the restriction
of $(\Txi_0,\tau_0)$ to the unit sphere, we may fix $\Txi_0$
and vary $\rho$ and $\tau_0$, obtaining after some rearrangement
the rescaled equation
\begin{equation}
\label{modsuperslow}
[(A^1)^{-1}\big(i\tau_0 + 
[i\sigma + \rho(B^{\xi_0 \xi_0}-\theta_1(|\Txi_0|^2+ 
|\tau_0+\sigma|^2)]  
+i A^{\Txi_0}\big) 
- \tilde \alpha I]\tilde \V=0,
\end{equation}
where $\sigma$ denotes variation in $\tau_0$.

%
%
%

\begin{prop} [{[Z.3--4]}] \label{branch} 
Under the hypotheses of Theorem \ref{mainsuf}, let
$\alpha_0=i\xi_{0_1}$ be a pure imaginary root of the
inviscid equation \eqref{slow} for some given
$\Txi_0$, $\tau_0$, i.e.  $\det(A^{\xi_0}+ \tau_0)=0$.
Then, associated with the corresponding root $\tilde \alpha$ 
in \eqref{superslow}, we have the following behavior, for
some fixed $\epsilon$, $\theta >0$ independent of $(\Txi_0,\tau_0)$:

(i) For $(\Txi_0$, $\tau_0)$ bounded distance $\epsilon$ away
from any branch singularity $(\Txi,\eta_j(\Txi))$ involving $\alpha$,
$\eta_j$ as defined in (H5),
the root $\tilde \alpha(\rho)$ in \eqref{superslow} such that
$\tilde \alpha(0):=\alpha_0$ bifurcates smoothly into $m$ roots
$\tilde \alpha_1,\dots, \tilde \alpha_m$, 
where $m$ is the dimension of $\ker (A^{\xi_0}+ \tau_0)$,
satisfying
\begin{equation}
\label{evalue}
\R \tilde \alpha_j\ge \theta \rho
\hbox{ or }
\R \tilde \alpha_j\le -\theta \rho;
\end{equation}
moreover, there is an analytic choice of eigenvectors spanning
the associated group eigenspace $\CalV$, in which coordinates the
restriction $\mA_\CalV^\pm$ of the limiting coefficient matrices 
for \eqref{ODE5} satisfy
\begin{equation}
\label{condition}
\theta\rho^2\le \R \mA_\CalV^\pm \le C \rho^2
\hbox{ or }
-C\rho^2\le  \R \mA_\CalV^\pm \le -\theta \rho^2,
\end{equation}
$C>0$,
in accordance with \eqref{evalue}, for $0<\rho \le \epsilon$.

(ii) For $(\Txi_0$, $\tau_0)$ lying at
a branch singularity $(\Txi,\eta_j(\Txi))$ involving $\alpha$,
the root $\tilde \alpha(\rho,\sigma)$ in \eqref{modsuperslow}
such that $ \tilde \alpha(0,0)=\alpha_0$ 
bifurcates (nonsmoothly) into $m$ groups of $s$ roots each:
\begin{equation}
\label{groups}
\{\tilde \alpha^1_1,\dots,\tilde \alpha^1_s\},
\dots, 
\{\tilde \alpha^m_1,\dots,\tilde \alpha^m_s\},
\end{equation}
where $m$ is the dimension of $\ker (A^{\xi_0}+ \tau_0)$ and
$s$ is some positive integer, such that, 
for $0\le \rho\le \epsilon$ and $|\sigma| \le \epsilon$,
\begin{equation}
\label{bevalues}
\tilde \alpha^j_k= \alpha +
\pi^j_k
+o(|\sigma|+|\rho|)^{1/s},
\end{equation}
and, in appropriately chosen (analytically varying) coordinate system,
the associated group eigenspaces $\CalV^j$ are spanned by
\begin{equation}
\label{bvectors}
\V^j_k=
\begin{pmatrix}
0\\
\vdots\\
0\\
e_k\\
\Pi^j e_k\\
(\Pi^j)^2 e_k\\
\vdots\\
(\Pi^j)^{s-1} e_k\\
0\\
\vdots\\
0 \\
\end{pmatrix}
+o(|\sigma|+|\rho|)^{1/s},
\end{equation}
where 
\begin{equation}
\label{pijk}
\pi^j_k:= \varepsilon^j
i(p\sigma - iq_k \rho )^{1/s}, 
\end{equation}
$\varepsilon:= 1^{1/s}$,
for $p(\Txi_0)$  real-valued and
uniformly bounded both above and away from zero and
$q_k$ denoting the eigenvalues of $Q\in \CC^{m\times m}$
such that $\sgn p Q\ge \theta >0$, 
\begin{equation}
\label{Pij}
\Pi^j:= \varepsilon^j
i(\sigma p I_m - iQ \rho )^{1/s}
\end{equation}
(well-defined, by definiteness of $p$, $Q$),
$e_k$ denote the standard basis elements in $\CC^{m}$, and,
moreover, the
restrictions $\mA_{\CalV^j}^\pm$ of the limiting coefficient matrices 
for \eqref{ODE5} to the invariant subspaces $\CalV^j_\pm$ are of form
$\Pi^j+o(|\sigma|+|\rho|)^{1/s}$, satisfying
\begin{equation}
\label{bcondition}
\theta \rho\R\tilde \pi^j\le  \R \mA_\CalV^\pm \le C \rho\R\tilde \pi^j,
\hbox{ or }
-C\rho\R\tilde \pi^j\le  \R \mA_\CalV^\pm \le -\theta \rho\R\tilde \pi^j,
\end{equation}
$\tilde \pi^j:= i(\sigma  - i \rho )^{j/s}$,
$C>0$, for $0<\rho \le \epsilon$.
\end{prop}

\begin{proof}
We here carry out the proof in the much simpler strictly hyperbolic case,
which permits a direct and relatively straightforward treatment.
A proof of the general case is given in Appendix C.
In this case, the dimension of $\ker (A^{\xi_0}+\tau_0)$ is one, 
hence $m$ is simply one.
Let $l(\Txi_0,\tau_0)$ and $r(\Txi_0,\tau_0)$ denote left and right 
zero eigenvectors of $ (A^{\xi_0}+\tau_0)=1$,  spanning co-kernel
and kernel, respectively; these are necessarily real, since 
$ (A^{\xi_0}+\tau_0)$ is real.
Clearly $r$ is also a right (null) 
eigenvector of $(A^1)^{-1} (i\tau_0 +iA^{\xi_0})$,
and $lA^1$ a left eigenvector.

Branch singularities are signalled by the relation
\begin{equation}
\label{vanish}
l A^1 r=0,
\end{equation}
which indicates the presence of a single Jordan chain 
of generalized eigenvectors of $(A^1)^{-1} (i\tau_0 + iA^{\xi_0})$
extending up from the genuine eigenvector $r$; we denote the
length of this chain by $s$.

\medbreak
\begin{obs}[{[Z.3]}] \label{B}
\textup{
Bound \eqref{multidiss} implies that
\begin{equation}
\label{Beq}
l B^{\xi_0 \xi_0} r \ge \theta>0,
\end{equation}
uniformly in $\Txi$.
}
\end{obs}

\begin{proof}[Proof of Observation]
In our present notation, \eqref{multidiss} can be written as
\begin{equation}
\R \sigma( -iA^{\xi_0} - \rho B^{\xi_0\xi_0}) \le -\theta_1 \rho,
\end{equation}
for all $\rho>0$, some $\theta_1>0$.
(Recall: $|\xi_0|\ge \theta_2>0$, by previous discussion).
By standard matrix perturbation theory [Kat],
the simple eigenvalue $\gamma=i\tau_0$ of $-iA^{\xi_0}$ perturbs analytically
as $\rho$ is varied around $\rho=0$, with perturbation series
\begin{equation}
\gamma(\rho)= i\tau_0 - \rho  l B^{\xi_0\xi_0}r + o(\rho).
\end{equation}
Thus, 
\begin{equation}
\R \gamma(\rho) = - \rho  l B^{\xi_0\xi_0}r + o(\rho)
\le -\theta_1 \rho,
\end{equation}
yielding the result. 
\end{proof}

In case (i), $\tilde\alpha(0)=\alpha$ is a simple eigenvalue of
$(A^1)^{-1} (i\tau_0 + iA{\Txi_0})$, and so perturbs analytically
in \eqref{superslow} as $\rho$ is varied around zero, with
perturbation series
\begin{equation}
\label{perti}
\tilde \alpha(\rho)= \alpha + \rho \mu^1 
+o(\rho),
\end{equation}
where $\mu^1=\tilde l (A^1)^{-1} \tilde r$,
$\tilde l$, $\tilde r$ denoting left and right eigenvectors
of $(A^1)^{-1} (i\tau_0 + A{\Txi_0})$.
Observing by direct calculation that $\tilde r= r$, 
$\tilde l=  l A^1/lA^1 r$, we find that
\begin{equation}
\mu^1= 
l B^{\xi_0\xi_0}r/lA^1r
\end{equation}
is real and bounded uniformly away from zero, by Observation \ref{B},
yielding the result
\eqref{evalue} for any fixed $(\Txi_0,\tau_0)$,
on some interval $0\le \rho\le \epsilon$, where $\epsilon$
depends only on a lower bound for $\mu^1$ and the maximum
of $\gamma''(\rho)$ on the interval $0\le \rho\le \epsilon$.
By compactness, we can therefore make a uniform choice of $\epsilon$
for which \eqref{evalue} is valid on the entire set of
$(\Txi_0,\tau_0)$ under consideration.
As $\mA_\CalV^\pm$ are scalar for $m=1$,
\eqref{condition} is in this case identical to \eqref{evalue}.
\medbreak
In case (ii),  $\tilde \alpha(0,0)=\alpha$
is an $s$-fold eigenvalue of $(A^1)^{-1} (i\tau_0 + iA{\Txi_0})$, 
corresponding to a single $s\times s$ Jordan block.
By standard matrix perturbation theory, the 
corresponding $s$-dimensional invariant subspace
(or ``total eigenspace'') varies analytically with $\rho$
and $\sigma$, and admits an analytic choice of basis
with arbitrary initialization at $\rho$, $\sigma=0$ [Kat].
Thus, by restricting attention to this subspace we
can reduce to an $s$-dimensional perturbation problem;
moreover, up to linear order in $\rho$, $\sigma$, the
perturbation may be calculated with respect to the
fixed, initial coordinization at $\rho$, $\sigma=0$.

Choosing the initial basis as a real, Jordan chain
reducing the restriction (to the subspace of interest)
of $(A^1)^{-1} (i\tau_0 + iA{\Txi_0})$ to $i$ times a
standard Jordan block, we thus reduce \eqref{modsuperslow}
to the canonical problem
\begin{equation}
\label{sslowcanonical}
\big(iJ+i\sigma M + \rho N-(\tilde \alpha-\alpha)\big)\V_I=0,
\end{equation}
where
\begin{equation}
\label{jordan}
J:=
\begin{pmatrix}
0 & 1 & 0 & \cdots &0 \\
0 & 0 & 1 & 0 & \cdots \\
0 & 0 & 0 & 1 & \cdots \\
\vdots &  \vdots & \vdots &  \vdots & \vdots\\
0 & 0 & 0 & \cdots&0 \\
\end{pmatrix},
\end{equation}
$\V_I$ is the coordinate representation of $\V$ in the
$s$-dimensional total eigenspace, and $M$ and $N$
are given by
\begin{equation}
\label{M}
M:=
\tilde L
(A^1)^{-1}
\tilde R
\end{equation}
and
\begin{equation}
\label{N}
N:=
\tilde L
(A^1)^{-1}
(B^{\xi_0 \xi_0}-\theta_1)  
\tilde R,
\end{equation}
respectively, where $\tilde R$ and $\tilde L$ are the
initializing (right) basis, and its corresponding (left) dual.

Now, we have only to recall that, as may be readily seen by 
the defining relation 
\begin{equation}
\tilde L 
(A^1)^{-1} (i\tau_0 + iA{\Txi_0})
\tilde R=J,
\end{equation}
or equivalently
$(A^1)^{-1} (i\tau_0 + iA{\Txi_0}) \tilde R=\tilde R J$ and
$\tilde L (A^1)^{-1} (i\tau_0 + iA{\Txi_0}) =J \tilde L$,
the first column of $\tilde R$ and the last row of $\tilde L$
are genuine left and right eigenvectors $\tilde r$ and $\tilde l$ of
$(A^1)^{-1} (i\tau_0 + iA{\Txi_0})$, hence without loss of generality
\begin{equation}
\tilde r=r, \quad \tilde l= pl A^1
\end{equation}
as in the previous (simple eigenvalue) case, 
where $p$ is an appropriate nonzero real constant.
Applying again Observation \ref{B}, we thus find that the crucial
$s,1$ entries of the perturbations $M$, $N$, namely 
$p$ and $pl(B^{\xi_0 \xi_0} - \theta_1)r=:q$, respectively,
are real, nonzero and of the same sign.
Recalling, by standard matrix perturbation theory, that this
entry when nonzero is the only significant one, we have reduced
finally (modulo $o(|\sigma|+|\rho|)^{1/s})$ errors)
to the computation of the eigenvalues/eigenvectors of
\begin{equation}
\label{computation}
i\begin{pmatrix}
0 & 1 & 0 & \cdots &0 \\
0 & 0 & 1 & 0 & \cdots \\
0 & 0 & 0 & 1 & \cdots \\
\vdots &  \vdots & \vdots &  \vdots & \vdots\\
p\sigma -iq\rho & 0 & 0 & \cdots&0 \\
\end{pmatrix},
\end{equation}
from which results \eqref{bevalues}--\eqref{pijk} 
follow by an elementary calculation,
for any fixed $(\Txi_0,\tau_0)$, and some choice of $\epsilon >0$;
as in the previous case, the corresponding global results 
then follow by compactness.
Finally, bound \eqref{bcondition} 
follows from \eqref{bevalues} and \eqref{pijk}
by direct calculation.
(Note that the addition of further $\CalO(|\sigma|+|\rho|)$
perturbation terms in entries other than the lower lefthand
corner of \eqref{computation} does not affect the result.
Note also that $\mA_{\CalV^j}^\pm$ are scalar in the strictly
hyperbolic case $m=1$, and $\CalV^j_\pm$ are simply eigenvectors
of $\mA_\pm$.)
This completes the proof in the strictly hyperbolic case.
\end{proof}

\begin{rem}\label{limitrem}
\textup{The detailed description of hyperbolic and glancing
modes given in Proposition \ref{branch}
readily yield the result of Lemma \ref{vcont} in the deferred
case $\R \lambda_0=0$, under the additional hypothesis (H5)
(exercise).
}
\end{rem}

\medbreak
{\bf 5.4.3. Main estimates.}
Combining the Evans-function estimates of Proposition \ref{ZS}
(and, in the case that refined but not uniform dynamical
stability holds, also those established in the course of
the proof of Proposition \ref{mainnec})
with the matrix perturbation analysis of Proposition \ref{branch},
we have all of the ingredients needed to carry out the basic
$L^1\to L^2$ resolvent estimates of Proposition \ref{mresbounds}.
In particular, in the uniformly dynamically stable case, 
they may be obtained quite efficiently
by Kreiss symmetrizer estimates generalizing those of the inviscid theory.
We refer to [GMWZ.1] or the article of M. Williams [W] in this volume 
for a presentation of the argument in the strictly hyperbolic,
Laplacian viscosity case.
With the results of Proposition \ref{branch}, a block version of
the same argument
applies in the general case, substituting invariant subspaces $\CalV^j$
for individual eigenvectors $\CalV$. 
We omit the details, which are beyond the scope of this article.

In the refined, but not uniformly dynamically stable case, the estimates
may be obtained instead as in [Z.3--4] by detailed pointwise estimates
on the resolvent kernel, using the second-order Evans function 
estimates carried out in the course of the proof of Proposition \ref{mainnec}
and the explicit representation formula for the resolvent kernel of
an ordinary differential operator [MZ.3, Z.4].
We refer to [Z.3--4] for an account of these more complicated arguments.

\medbreak
{\bf 5.4.4. Derivative estimates.}
Improved derivative estimates, $|\beta|=1$, may now easily be obtained
by a method introduced by Kreiss and Kreiss [KK] in the one-dimensional case;
see [GMWZ.1] or the notes of M. Williams [W] in this volume.
Specifically, given a differentiated source $f=\partial_{x_1}g$, 
in resolvent equation $(L_\Txi-\lambda)U=f$, consider first
the {\it auxiliary equation}
\begin{equation}
\label{auxiliary}
L_0 W= (B^{11}W_{x_1})_{x_1} - (A^1 W)_{x_1}= \partial_{x_1}g.
\end{equation}

Using the conservative (i.e., divergence-form) structure of
$L_0$, we may integrate \eqref{auxiliary} from $-\infty $ to $x_1$
to obtain a reduced ODE 
\begin{equation}\label{redode}
z_2'-\alpha(x_1)z_2= \begin{cases}0,\\G, \end{cases},
\qquad W=\Phi(z_2) 
\end{equation}
analogous to that obtained in [KK, GMWZ.1] for the strictly parabolic case,
where $z_2:= B^{11}_{II}W$ and $G\sim g$, $W\sim z_2$.
More precisely, the inhomogeneous version 
$z_2'-\alpha(x_1)z= 0$ is the linearization about
$\bU$ of the (integrated) traveling-wave ODE \eqref{wode},
from which we obtain by transversality $\gamma\ne 0$ that
the solution of \eqref{redode} is unique modulo $\bU'$.
That a solution exists follows easily from the fact \eqref{indrel}
relating the signatures of $\alpha_\pm$ to those of $A^1_\pm$
(a consequence of transversality, together with
our assumptions on the profile; see Remark \ref{ell}),
together with standard arguments for asymptotically constant ODE as
in, e.g., [He, Co, CL]; 
see, in particular, the argument of  section 10.2, [GMWZ.1] in
the strictly parabolic case, for which $\alpha$ reduces to $A^1$.
See also [Go.2, KK], or the article [W] by M. Williams
in this volume.

Indeed, imposing an additional condition 
\begin{equation}\label{add}
\langle \ell, W\rangle=0,
\end{equation}
where $\ell$ is any constant vector satisfying $\langle \ell, \bU'\rangle
=\ell \cdot [U]\ne 0$, we have [He, Co, CL, GMWZ.1, Go.2, KK, W] the bound
$|z_2|_{W^{1,p}}\le C|G|_{L^p}$ for any $p$, yielding in particular
\begin{equation}\label{auxbd}
|W|_{L^1}+|B^{11} W_{x_1}|_{L^1}\le C|g|_{L^1}.
\end{equation}
The reduction to form \eqref{redode} goes similarly as the
reduction of the nonlinear traveling-wave ODE in Section 3.1;
we leave this as an exercise.

Setting now $U=W+Y$, and substituting the auxiliary equation
into the eigenvalue equation, we obtain equation
\begin{equation}
\label{residual}
\begin{aligned}
(L_\Txi-\lambda)Y&= \CalO(\rho)(|W|_{L^1}+|B^{11}W_{x_1}|_{L^1})\\
&=\CalO(\rho |g|_{L^1}),
\end{aligned}
\end{equation}
for the residual $Y$, from which we obtain the desired bound
from the basic $L^1\to L^2$ estimate of Section 5.4.3 above.

Alternatively, one may obtain the same bounds as
in [Z.3--4] by direct computation on the original resolvent equation.  
Both methods are based ultimately on the
fact that, for Lax-type shocks, the only $L^1$ time-invariants
of solutions of the linearized equations are those determined
by conservation of mass.
This property is shared by over- but not undercompressive
shocks, hence the degraded bounds in the latter case;
for further discussion, see [LZ.2, Z.2].
This completes the proof of Proposition \ref{mresbounds}, and the 
analysis.

%
%

\bigbreak
\clearpage




\appendix

\section{Appendix A.  Semigroup facts}


\begin{defi} 
\label{closed}
\textup{
Given a Banach space $X$, and a linear operator
$L: \CalD(L) \subset X \to X$, we say that 
$L$ is densely defined in $X$ if $\CalD(L)$ is dense in $X$.
We say that $L$ is closed if $u_n\to u$ and $x_n:=Lu_n \to x$
(with respect to $|\cdot|_X$) 
for $u_n\in \CalD(L)$ and $x_n\in X$ 
implies that $u\in \CalD(L)$
and $Lu=x$.  
$\CalD(L)$ is called the domain of $L$.
We define associated domains $\CalD(L^n)$ by induction as the set of 
$x\in \CalD(L^{n-1})$ such that $L^{n-1}x\in \CalD(L)$.
}
\end{defi}

\begin{exs}
\label{canonical}
1.  If $L$ is closed, show that $\CalD(L)$ is a Banach space
under the canonical norm $|u|_{\CalD(L)}:= |u|_X+ |Lu|_X$,
i.e.,  each Cauchy sequences $u_n \in \CalD(L)$ has a limit
$u\in \CalD(L)$. (First note that $u_n$ and $x_n:=Lu_n$ are Cauchy
with respect to $|\cdot|_X$, hence have limits $u$, $x$ in $X$.)
With this choice of norm, $L: \CalD(L) \to X$
is a bounded operator, hence (trivially) closed in the usual, 
Banach space sense.
\smallbreak
2. Show that $L-\lambda$ is closed if and only if is $L$.
\smallbreak
3.  Show that $|u|\le C|Lu|$ for $L$ closed implies that $\range (L)$ is closed.
\end{exs}

\begin{defi} 
\label{spectrum}
\textup{
Given a Banach space $X$, and a closed, densely defined linear operator
$L: \CalD(L) \subset X \to X$, the resolvent set $\rho(L)$,
written $\rho_X(L)$ when we wish to identify the space,
is defined as the set of $\lambda\in \CC$ for
which $(\lambda-L)$ has a bounded inverse
$(\lambda-L)^{-1}: X \to \CalD(L)$.  The operator $(\lambda-L)^{-1}$
is called the resolvent of $L$.
The spectrum $\sigma(L)$ of $L$,
written $\sigma_X(L)$ when we wish to identify the space,
is defined as the complement
of the resolvent set, $\sigma(L):= \rho(L)^c$.
}
\end{defi}

\begin{ex}
\label{density}
If $L$ is densely defined and $\lambda\in \rho(L)$ for some $\lambda$, 
show that $\CalD(L^n)=\range (\lambda-L)^{-n}$
is dense in $X$.
(Show by induction that each $\CalD(L^{n+1})$ is dense in $\CalD(L^n)$.)
\end{ex}

\begin{lem}\label{analytic} 
For a closed, densely defined
operator $L$ on Banach space $X$,
the resolvent operator $(\lambda- L)^{-1}$
is analytic in $\lambda$ with respect to $|\cdot|_X$ for $\lambda$ in
the resolvent set $\rho(L)$.
\end{lem}
\begin{proof}
For $\lambda$ sufficiently near $\lambda_0\in \rho(L)$, we may expand
\begin{equation}
\begin{aligned}
(\lambda-L)&= (\lambda_0-L)^{-1}
\Big(I- (\lambda_0-\lambda)(\lambda_0-L)^{-1} \Big)^{-1}\\
&=
(\lambda_0-L)^{-1}
\sum_{j=0}^\infty  
\Big((\lambda_0-\lambda)(\lambda_0-L)^{-1}\Big)^j,
\end{aligned}
\end{equation}
using
$|(\lambda_0-\lambda)(\lambda-L)^{-1}|\le C |\lambda_0-\lambda|$
by the definition of resolvent set and the
Neumann expansion $(I-T)^{-1}=\sum_{j=0}^\infty T^j$
for $|T|$ sufficiently small.
\end{proof}

\begin{ex}
\label{relative}
1. Let $L(\alpha)$ be a family of closed, densely defined
operators on a single Banach space $X$, 
such that $L$ is ``relatively analytic'' in $\alpha$ 
with respect to $\big(L(\alpha_0)-\lambda_0\big)$, 
$\lambda_0\in \rho(L(\alpha_0)$,
in the sense that $L(\alpha)(\lambda_0- L(\alpha_0))^{-1}$
is analytic in $\alpha$ with respect to $|\cdot|_X$.
Show that the family of resolvent operators 
$(\lambda- L(\alpha))^{-1}$ is analytic in 
$(\alpha,\lambda)$ in a neighborhood of $(\alpha_0,\lambda_0)$.
%
%
\smallbreak
2. Including in \eqref{auxHF2}
formerly discarded beneficial $w^{II}$ terms in the 
energy estimates from which it derives, we
obtain the sharpened version
\begin{equation}
\label{refinedauxH}
\begin{aligned}
(\R \lambda +\theta_1)
\Big(|W|_{\hat H^1} + |\partial_{x_1}^2 w^{II}| \Big)
 &\le
C_1\Big(|f|_{\hat H^1} + C_1|W|\Big),
\end{aligned}
\end{equation}
revealing smoothing in $w^{II}$ for $\R \lambda > -\theta_1$.
Show that this implies analyticity in $\Txi$ with respect to 
$|\cdot|_{\hat H^1}$ of $L_\Txi(\lambda_0-L_{\Txi_0})^{-1}$
(note: $\partial_\Txi L_\Txi$ is a continuous-coefficient
first-order differential operator, for which the derivative
falls only on $w^{II}$ components), and conclude that
$(\lambda-L_\Txi)^{-1}$ is analytic in $(\Txi, \lambda)$ for
$\lambda \in \rho_{\hat H^1}(L_\Txi)$ and $\R \lambda \ge -\theta_1$
with $\theta_1>0$ sufficiently small.
\end{ex}

\begin{lem} [Resolvent identities]
\label{resident}
For $\lambda$, $\mu$ in the resolvent set of a closed, densely defined 
operator $L: \CalD(L)\to X$ on a Banach space $X$, 
\begin{equation}
\label{2residenteq}
(\lambda-L)^{-1}(\mu-L)^{-1}= 
\frac{(\mu-L)^{-1}-(\lambda-L)^{-1} }
{\lambda-\mu}=
(\mu-L)^{-1} (\lambda-L)^{-1}
\end{equation}
on $X$ and
$L(\lambda-L)^{-1}=
\lambda(\lambda-L)^{-1}-I=
(\lambda-L)^{-1}L$
on $\CalD(L)$: in particular, 
\begin{equation}
\label{residentbd}
(\lambda-L)^{-1}u = \lambda^{-1}\big(u+(\lambda-L)^{-1}Lu\big)
\quad
\text{\rm for $u\in \CalD(L)$.}
\end{equation}
\end{lem}
\begin{proof}
Rearranging $(\mu-L)(\mu-L)^{-1}=I$, we obtain 
$(\lambda-L)(\mu-L)^{-1}=(\lambda-\mu)(\mu-L)^{-1} I$, from
which the first equality of \eqref{2residenteq} follows upon application of
$(\lambda-L)^{-1}$ from the left, and the second by symmetry. 
Rearranging defining relation
$(\lambda-L)(\lambda-L)^{-1}=I= (\lambda-L)^{-1}(\lambda-L)$,
we obtain the second assertion, 
whereupon \eqref{residentbd} follows by further rearrangement 
after multiplication by $\lambda^{-1}$. 
\end{proof}

\begin{ex}
\label{redundant}
Assuming that the resolvent set is open,
recover the result of Lemma \ref{analytic} directly from
the definition of derivative, using resolvent identity
\eqref{2residenteq} to establish 
differentiability, $(d/d\lambda) (\lambda-L)^{-1}=-(\lambda-L)^{-2}$.
\end{ex}

\begin{defi}[{[Pa]}]\label{semigroup}
\textup{
A $C^0$ semigroup on Banach space $X$ is a family of bounded
operators $T(t)$ satisfying the properties
(i) $T(0)=I$, (ii) $T(t+s)=T(t)T(s)$ for every
$t$, $s\ge 0$, and (iii) $\lim_{t\to 0^+} T(t)x=x$
for all $x\in X$. The generator $L$ of the semigroup
is defined as $Lx=\lim_{t\to 0^+}(T(t)x-x)/t$ on the 
domain $\CalD(L)\subset X$ for which the limit exists.
We write $T(t)=e^{Lt}$.
Every $C^0$ semigroup satisfies $|e^{Lt}|\le Ce^{\gamma_0t}$
for some $\gamma_0$, $C$; see [Pa], Theorem 2.2.
}
\end{defi}

\begin{rem}
\label{semijust}
\textup{For a $C^0$ semigroup, $(d/dt)e^{Lt}f=L e^{Lt}f= e^{Lt}Lf$
for all $f\in \CalD(L)$ and $t\ge 0$; see [Pa], Theorem 2.4(c).
Thus, $e^{Lt}$ is the solution operator for initial-value problem
$u_t=Lu$, $u(0)=f$, 
justifying the exponential notation.
}
\end{rem}


\begin{prop}[Generalized Hille--Yosida theorem]
\label{semigpcriteria}
An operator $L: \CalD(L)\to X$ is the generator of a 
$C^0$ semigroup $|e^{Lt}|\le Ce^{\gamma_0 t}$ 
on $X$ with domain $\CalD(L)$ if and
only if:
(i) it is closed and densely defined, and
(ii) 
$\lambda \in \rho(L)$ and
$|(\lambda-L)^{-k}|\le C|\lambda-\gamma_0|^{-k}$ for
sufficiently large real $\lambda$,
in which case also
$|(\lambda-L)^{-k}|\le C|\R\lambda-\gamma_0|^{-k}$ for all
$\R \lambda>\gamma_0$.
\end{prop}

\begin{proof}
($\Rightarrow$) By the assumed exponential decay, the Laplace transform 
$\hat T:=\int_0^\infty e^{-\lambda s}T(s)\, ds$ is
well-defined for $\R \lambda>\gamma_0$, with 
$|\hat T(\lambda)|\le C|\R\lambda-\gamma_0|^{-1}$.
By the properties of the solution operator, we have also $\hat TL=L\hat T=
\widehat{LT}$ on $\CalD(L)$, as well as
$\widehat{LT}= \widehat{\partial_t T}=\lambda \hat T- I$.
Thus, $(\lambda-L)\hat T=\hat T(\lambda-L)=I$ on $\CalD(L)$.
But, also, 
\begin{equation}
\begin{aligned}
h^{-1}(T(h)-I)\hat T&= 
h^{-1}(e^{\lambda h}-I)\hat T - h^{-1}\int_0^h e^{-\lambda (s-h)}T(s)\, ds\\
\end{aligned}
\end{equation}
approaches $\lambda\hat T- I $ as $h\to 0^+$,
yielding $\hat TX\subset \CalD(L)$ and $L\hat T=\lambda \hat T -I$ on $X$,
by definition.  Thus, $\lambda \in \rho(L)$ for $\Re \lambda>\gamma_0$,
and $\hat T=(\lambda-L)^{-1}$, whence 
$|(\lambda-L)^{-k}\le C|\Re \lambda -\gamma_0|^{-k}$ for $k=1$.
The bound for general $k$ follows from the computation
$\hat T^k=\lt \overbrace{(T*T\cdots*T)}^{k \, \text{\rm times}}=
\lt(t^k T/k!)$, where $\lt$ denotes Laplace transform and ``$*$''
convolution, $g*h(t):=\int_0^t g(t-s)h(s)\, ds$, together
with $\int_0^\infty e^{-z}z^k/k! \, dz \equiv 1$ for all $k\ge 0$
(exercise). 

($\Leftarrow$)  Without loss of generality, take
$\gamma_0~0$.  For real numbers $z$, $e^z=(e^{-z})^{-1}=
\lim_{n\to \infty}(1-z/n)^{-n}$ gives a stable approximation
for $z<0$.  This motivates the introduction of approximants
\begin{equation}
T_n:=(I-Lt/n)^{-n},
\end{equation}
which, by resolvent bound (ii), satisfies $|T_n|\le C$.
Restricting to the dyadic approximants $T_{2^j}$, 
and using the elementary difference formula
$a^n-b^n=(a^{n-1}+ a^{n-2}b+\dots + b^{n-1})(a-b)$ for commuting
operators $a$, $b$ together with bound (ii), we find for $x\in \CalD(L^2)$
that
\begin{equation}
\begin{aligned}
|(T_{2^{j+1}}-T_{2^j})x|&\le
2^jC^2
\Big| \Big((I-Lt/2^{j+1})^{-2}-(I-Lt/2^j)^{-1}\Big)x  \Big| \\
&=
2^{j}C^2
\Big|(I-Lt/2^{j+1})^{-2}(I-Lt/2^j)^{-1}\Big|
t^2 |L^2 x|2^{-2j-2} \\
&\le C^4 t^2 |L^2 x|2^{-j-2} 
\end{aligned}
\end{equation}
and thus the sequence converges geometrically to a limit $Tx$ for
all $x\in \CalD(L^2)$.  Since $\CalD(L^2)$ is dense in $X$ (Exercise
\ref{density}), and the $T_n$ are uniformly bounded, this implies
convergence for all $x$, defining a bounded operator $T(t)$.

Clearly, $T(0)=T_n(0)\equiv I$.  Also,
$T_nx\to x$ as $t\to 0$ for $x\in \CalD(L^n)$,
whence we obtain $Tx\to x$ for all $x$ by uniform
convergence of the dyadic approximants and density of 
convergence of the approximants and density of 
each $\CalD(L^n)$ in $X$.
Finally, for $x\in \CalD(L)$, we have $(d/dt)(I-Lt/n)x=-Lx/n$, from which
we obtain $(d/dt)T_nx=-n(I-Lt/n)^{-n-1}(-L/n)x=L_n(t) T_n x$,
where $L_n(t):= L(I-Lt/n)^{-1}$ is bounded for bounded $t$ (exercise).
%
Observing that $|LT_n(t)x|=|T_n Lx|\le C|Lx|$ for all $t$,
we find for $x\in \CalD(L^2)$ that 
\begin{equation}
e_n(r,t):=u(t)-v(t)= T_n(r+t)x- T_n(t)T_n(r)x,
\end{equation}
$u'=L_n(s)u$, $v'=L_n(r+s)v$,
$u(0)=v(0):= T_n(r)x$,
satisfies
\begin{equation}
\begin{aligned}
e_n'=L_n(s+r)e_n + 
f_n(r,s), \quad e(0)=0,
\end{aligned}
\end{equation}
$f_n(r,s):= (L_n(s)-L_n(s+r))T_n(r+s)x$, where 
\begin{equation}
\begin{aligned}
|f_n(r,s)|&=
|\Big((I-Ls/n)^{-1}-(I-L(s+r)/n)^{-1}\Big) T_n Lx|\\
&=|(I-Ls/n)^{-1}(I-L(s+r)/n)^{-1}T_n (r/n)L^2x|\\
&\le |L^2x|r/n
\end{aligned}
\end{equation}
goes to zero uniformly in $n$ for bounded $r$.
Expressing
\begin{equation}
e_n(r,t)=\int_0^t T_n(t+r)T_n(t+r-s)^{-1}f_n(r,s)\, ds
\end{equation}
using Duhamel's principle/uniqueness of solutions of bounded-coefficient
equations, and observing (Exercise \ref{keyest} below)
that $ |T_n(t+r)T_n(t+r-s)^{-1}|\le C$, we thus have
$|e_n(r,t)|\le Ct|L^2x|r/n \to 0$ uniformly in $n$ for bounded $r$, $t$,
verifying semigroup property (ii) Definition \ref{semigroup}
in the limit for $x\in \CalD(L^2)$, and thus, by continuity, for all $x\in X$.
\end{proof}

\begin{rem}[{[Pa]}]\label{euler}
\textup{
The approximants $T_n(t)$ correspond to the finite-difference
approximation obtained by first-order implicit Euler's method with mesh
$t/n$.
In the language of numerical analysis,
boundedness of $|T_n|$ corresponds to ``$\CalA$-stability'' of the scheme,
i.e., suitability for 
``stiff'' ODE.
}
\end{rem}

\begin{ex}\label{keyest}
Show by direct computation that 
\begin{equation}
T_n(r+s)T_n(s)^{-1}=
\Big(\Big(\frac{r}{r+s}\Big)(I-L(r+s)/n)^{-1}+ \Big(\frac{s}{r+s}\Big)\Big)^n.
\end{equation}
Assuming resolvent bound (ii), and thus $|(I-L(r+s)/n)^{-k}|\le C$,
show that $|T_n(r+s)T_n(s)^{-1}|\le C$, using binomial expansion
and the fact that resolvents commute.
\end{ex}

\begin{ex}\label{optest}
Show by careful expansion of
\begin{equation}
\begin{aligned}
T_{n+1}-T_n&=
\Big((I-Lt/(n+1))^{-(n+1)}- 
(I-Lt/(n+1))^{-n)} \Big)\\
&\quad + 
\Big((I-Lt/(n+1))^{-n} - 
(I-Lt/n))^{-n}\Big)
\end{aligned}
\end{equation}
that $|T_{n+1}x-T_n x|\le C(|L^2 x|+|Lx|+|x|)/n^2$, so that the entire
sequence $\{T_n\}$ is convergent for $x\in \CalD(L^2)$.
\end{ex}

\begin{ex}\label{largelambda}
Show directly, using the same Neumann expansion argument 
used to prove analyticity
of the resolvent that $|(\lambda-L)^{-1}|\le (\lambda-\gamma_0)^{-1}$
for sufficiently large real $\lambda$ implies
$|(\lambda-L)^{-1}|\le (\R\lambda-\gamma_0)^{-1}$ for all $\R \lambda>\gamma_0$.
\end{ex}

\begin{defi}\label{dissipative}
\textup{
We say that a linear operator $L$ is {dissipative}
if it satisfies an a priori estimate
\begin{equation}
\label{enkA}
{|\lambda - \gamma_0|^k} |U|_X \le
C |(\lambda- L)^{k}U|_{X}
\end{equation}
for all real $\lambda > \lambda_0$.
Condition \eqref{enkA} together with $\range (\lambda-L)=X$ is equivalent
to condition (ii) of Proposition \ref{semigpcriteria}.
}
\end{defi}

\begin{cor}
\label{adjver}
An operator $L: \CalD(L)\to X$ is the generator of a 
$C^0$ semigroup $|e^{Lt}|\le Ce^{\gamma_0 t}$ 
on $X$ with domain $\CalD(L)$ if and
only if:
(i) it is closed and densely defined, and
(ii') 
both $L$ and $L^*$ are dissipative, where
$L^*: X^* \to \CalD(L)^*$ defined by
$\langle L^*v,u\rangle:=\langle v, Lu\rangle$ denotes the adjoint of $L$.
\end{cor}

\begin{proof}
See Corollary 4.4 of [Pa], or exercise \ref{adjoint}, below.
\end{proof}

\begin{cor}[Generalized Lumer--Phillips theorem]
\label{lumer}
A densely defined operator $L: \CalD(L)\to X$ generates a 
$C^0$ semigroup $|e^{Lt}|\le Ce^{\gamma_0 t}$ 
on $X$ with domain $\CalD(L)$ if and only if $L$ is dissipative
and $\range(\lambda_0-L)=X$ for some real $\lambda_0>\gamma_0$.
In particular, for a densely defined dissipative operator $L$ 
with constants $C$, $\gamma_0$, the real
ray $(\gamma_0, +\infty)$ consists either entirely of spectra,
or entirely of resolvent points.
\end{cor}

\begin{proof}
See Theorem 4.3 of [Pa], or Exercise \ref{lp} below. 
\end{proof}

\begin{rem}\label{dissmot}
\textup{
In the contractive case, $C=1$, dissipativity
is equivalent to $\R\langle v, u\rangle\le 0$ for
some $v\in X^*$ such that $\langle v, u\rangle=|u||v|$;
see Theorem 4.2, [Pa].
If $X$ is a Hilbert space, contractive dissipativity ($C=1$) 
of both $L$ and its adjoint
$L^*$ reduce to the single condition $\R\langle u, Lu\rangle \le 0$,
which is therefore necessary and sufficient that $L$ generate
a $C^0$ contraction semigroup; see Exercise \ref{hilex}.
This corresponds to $(d/dt)(1/2)|u|^2\le 0$, motivating our terminology.
}
\end{rem}

\begin{ex}\label{adjoint}
For a closed operator $L$, show that
(ii) of Proposition \ref{semigpcriteria} is equivalent
to (ii') of Proposition \ref{adjver}.
The forward direction follows by the general facts that $|A|=|A^*|$
and $(A^{-1})^{*}=(A^{*})^{-1}$.
The reverse follows by closure of $\range (\lambda-L)$ (Exercise
\ref{canonical}.3), the Hahn--Banach Theorem,
and the observation that
$0=\langle f,(\lambda- L)u\rangle =\langle(\lambda- L)^*f, u\rangle$ 
for all $u\in \CalD(L)$ implies $f=0$, which together
yield $\range (\lambda-L)=X$.
\end{ex}

\begin{ex}\label{lp}
Let $L$ be a dissipative operator with constants $C$, $\gamma_0$,
such that $\range(\lambda_0-L)=X$ for some real $\lambda_0>\gamma_0$.

1. Using Exercise \ref{canonical}.3, show that $L$ is closed. 

2. If $(\lambda_1-L)u_1=(\lambda_2-L)u_2$ for $\lambda_j>\gamma_0$,
show that $|u_1-u_2|\le C|\lambda_1-\gamma_0|^{-1}|\lambda_1-\lambda_2||u_2|$.

3. Using the result of 2, show that
$\range(\lambda_n-L)=X$ for $\lambda_n\to \lambda>\gamma_0$ implies
$\range(\lambda-L)=X$. (Show that $\{u_j\}$ is Cauchy for 
$(\lambda_j-L)u_j:=x$,
then use closure of $L$.) 
Conclude that $\range(\lambda-L)=X$ for all real $\lambda>\gamma_0$,
by the fact that $\rho(L)$ is open.
\end{ex}

\begin{ex}\label{hilex}
For a densely defined linear operator $L$ on a Hilbert space $X$, and
$u\in \CalD(L)$, show 
by direct inner-product expansion that
$|(\lambda - L)u|^2\ge |(\lambda-\gamma_0) u|^2$ for real
$\lambda >\gamma_0$ if and only if $\langle u, Lu\rangle \le \gamma_0|u|^2$,
if and only if
$\langle u, L^*u\rangle \le \gamma_0|u|^2$.
\end{ex}

\begin{rem}
\label{pruss}
\textup{
Proposition \ref{semigpcriteria} includes the rather
deep stability estimate $|e^{Lt}|\le Ce^{\gamma_0 t}$
converting global spectral information to a sharp rate
of linearized time-exponential decay.
For a Hilbert space, a useful alternative criterion for
exponential decay $|e^{LT}|\le Ce^{\gamma_0 t}$
has been given by Pr\"uss [Pr]:
namely, that, for 
some $\gamma<\gamma_0$, $\{\lambda: \R \lambda \ge \gamma\}\subset \rho(L)$,
and $|(\lambda-L)^{-1}\le M$ on $\R \lambda =\gamma$.
A useful observation of Kapitula and Sandstede 
[KS, ProK] is that, in the case that
there exist isolated, finite-multiplicity eigenvalues of $L$ in
$\{\lambda: \R \lambda \ge \gamma\}\subset \rho(L)$, the same
result may be used to show exponential convergence to the union
of their associated eigenspaces.
}
\smallbreak
\textup{
However, these tools are not available in the case, as here, that
essential spectrum of $L$ approaches the contour $\R \lambda=\gamma_0$
($\gamma_0=0$ in our case).
In this situation, one may take the alternative approach of direct
estimation using the inverse Laplace transform formula,
given just below.
However, notice that we do not get from this ``local'' formula
the Hille-Yosida bound, which is not implied by behavior
on any single contour $\R \lambda=\gamma$ (recall that $\lambda\to \infty$ 
in the proof of Proposition \ref{semigpcriteria}).
Indeed, the spectral resolution formula holds under much weaker bounds than
required for existence of a semigroup; see Exercise \ref{misc}.
A theme of this article is that we can nonetheless get somewhat
weaker, time-averaged bounds by similarly simple criteria, and
that, provided we have available a nonlinear smoothing estimate analogous to 
\eqref{multiEexp},
we can use these to close a nonlinear argument in which the 
deficiencies of our linearized bounds disappear.
The result is effectively a ``Pr\"uss-type'' bound
on the high-frequency part of the solution operator, analogous
to the bounds obtained by Kapitula and Sandstede by spectral decomposition
in the case that slow- and fast-decaying modes are spectrally 
separated.
}
\end{rem}

\begin{prop}\label{ILTprop}
Let $L:\CalD(L)\to X$ be a closed, densely defined operator
on Banach space $X$, generating a $C^0$ semigroup
$|e^{L t}|\le Ce^{\gamma_0 t}$: equivalently, satisfying
resolvent bound \eqref{resk}.
Then, for $f\in \CalD(L)$,
\begin{equation}
\label{basicILT}
e^{Lt}f=
\pv \int_{\gamma -i\infty}^{\gamma+i\infty}
e^{\lambda t}
(\lambda- L)^{-1}f \, d\lambda 
\end{equation}
for any $\gamma>\gamma_0$, with convergence in $L^2([0,\infty);X)$.  
For $f$, $Lf\in \CalD(L)$, \eqref{basicILT} converges pointwise
for $t>0$, with uniform convergence
on compact sub-intervals $t\in [\epsilon, 1/\epsilon]$
at a rate depending only on the bound for $|f|_X+|LF|_{X} +|L(Lf)|_X$.
\end{prop}

\begin{proof}
By the bound $|e^{Lt}|\le Ce^{\gamma_0 t}$, the Laplace transform
$\hat u$ of $u(t):= e^{Lt}f$ is well-defined for $\R \lambda >\gamma_0$,
with 
\begin{equation}
\label{inv}
u(t)=
\pv \int_{\gamma -i\infty}^{\gamma+i\infty}e^{\lambda t}\hat u(\lambda)
\, d\lambda 
\end{equation}
for $\gamma>\gamma_0$, 
by convergence of Fourier integrals on $L^2$, 
and the relation between Fourier and Laplace transform
(see \eqref{ftltrel}). 
Moreover, $u(t)\in C^1([0,\infty);X)$
by Remark \ref{semijust}, with $u_t=Lu$ and 
$|u_t|$, $|Lu|\le Ce^{\gamma_0 t}$,
whence $\widehat {u_t}=\lambda \hat u - u(0)= \lambda \hat u -f$
and also $\widehat{Lu}=L\hat u$.
Thus, $(\lambda-L)\hat u=f$, giving \eqref{basicILT} by \eqref{inv}
and invertibility of $\lambda -L$ on $\R \lambda >\gamma_0$.

Next, suppose that $f$, $Lf\in \CalD(L)$,
without loss of generality taking $\lambda\ne 0$.\footnote{
The restriction $\gamma\ne 0$ may be
removed by considering the shifted semigroup 
$\tilde T(t):=e^{-\gamma_0 t}e^{Lt}$.
(Exercise: verify that $\tilde T$ is generated by $L-\gamma_0$.)}
Expanding
\begin{equation}
\label{expansion}
(\lambda-L)^{-1}f= \lambda^{-2}\big((\lambda-L)^{-1}L\cdot Lf + Lf\big)
+ \lambda^{-1}f
\end{equation}
using \eqref{residentbd}, we may split the righthand side of \eqref{basicILT} 
into the sum of an integral
\begin{equation}
\label{absc}
\pv \int_{\gamma -i\infty}^{\gamma+i\infty}
e^{\lambda t}
\lambda^{-2}\Big((\lambda- L)^{-1}L\cdot Lf + Lf\Big) \, d\lambda 
\end{equation}
that is absolutely convergent on bounded time-intervals $t\in [0,\delta^{-1}]$
and the integral
$\pv  \int_{\gamma-i\infty}^{\gamma+ i\infty} 
e^{\lambda t} \lambda^{-1}\, d\lambda f= f$,
uniformly convergent on compact intervals $t\in [\delta,\delta^{-1}]$;
see Exercise \ref{conv} below.
Note that \eqref{absc} converges to 
\begin{equation}
\label{zeroconv}
\pv  \int_{\gamma-i\infty}^{\gamma+ i\infty} 
\lambda^{-2}\Big((\lambda- L)^{-1}
L\cdot Lf + Lf\Big) \, d\lambda =0
\end{equation}
as $t\to 0$, directly verifying semigroup property (iii).
For an alternative proof, see Corollary 7.5 of [Pa].
\end{proof}

\begin{prop}\label{inhomspecres}
Let $L:\CalD(L)\to X$ be a closed, densely defined operator
on Banach space $X$, generating a $C^0$ semigroup
$|e^{L t}|\le Ce^{\gamma_0 t}$:
equivalently, satisfying resolvent bound \eqref{resk}.
Then, for $f\in L^1([0,T]; \CalD(L))$ 
or $f\in L^2([0,T];X)$, 
\begin{equation}
\label{inhomILFTA}
\begin{aligned}
\int_0^t e^{L(t-s)}f(s) \, ds &=
\pv \int_{\gamma -i\infty}^{\gamma+i\infty}
e^{\lambda t}
(\lambda- L)^{-1}
{\widehat {f^T}}(\lambda) \, d\lambda 
\end{aligned}
\end{equation}
for any $\gamma>\gamma_0$, $0\le t\le T$,
with convergence in $L^2([0,T]; X)$,
where ${\hat g}$ denotes Laplace transform of $g$ and
$f^T(x,s):= f(x,s)$ for $0\le s\le T$ and zero otherwise.
For $f$, $Lf\in L^1([0,T]; \CalD(L))$ or 
$f \in L^q([0,T]; \CalD(L))$, $q>1$, 
the convergence is pointwise, and uniform for $t\in [0,T]$.
\end{prop}

\begin{proof}
Equivalence and
convergence in $L^2([0,T]; X)$ for $f\in L^1([0,T]; \CalD(L))$ 
follows similarly as in the proof of Proposition \ref{basicILT},
observing that $u:=\int_0^t e^{L(t-s)}f(s)\, ds$ satisfies
$u_t-Lu=f^T$ for $0\le t\le T$ with $u(0)=0$, and $\hat u$,
$\hat f$ are well-defined, with $(\lambda -L)\hat u=\hat f^T$.
Likewise, we may obtain
pointwise convergence for $f$, $Lf\in L^1([0,T]; \CalD(L))$ 
using expansion \eqref{expansion} and the Hausdorff--Young inequality
\begin{equation}
|\hat g^T(\lambda)|=|\int_0^T e^{-\lambda s} g^T(s) ds|
\le (1+e^{-\R \lambda T})|g|_{L^1(t)}
\end{equation}
to obtain a uniformly absolutely convergent integral plus the
uniformly convergent integral 
\begin{equation}
\begin{aligned}
\pv \int_{\gamma-i\infty}^{\gamma+i\infty}
\lambda^{-1} e^{\lambda t} \widehat{f^T}(\lambda)\, d\lambda
&=
\int_0^T f(s) ds
\pv \int_{\gamma-i\infty}^{\gamma+i\infty}
\lambda^{-1} e^{\lambda (t-s)} \widehat{f^T}(\lambda)\, d\lambda\\
&=\int_0^t f(s) ds;
\end{aligned}
\end{equation}
see Exercise \ref{conv} below.

Convergence in $L^2$ for $f\in L^2([0,T]; X)$ follows using
Parseval's identity,
\begin{equation}
\pv \int_{\gamma-i\infty}^{\gamma+i\infty} |\widehat{g^T}(\lambda)|^2\, d\lambda
=
\int_0^T |e^{-\R \lambda s}g(s)|^2\, ds
\le (1+ |e^{-\R \lambda T}|)^2 |g|_{L^2[0,T]}^2,
\end{equation}
together with the fact that $|e^{\lambda t}(\lambda-L)^{-1}|\le C(\gamma)$
for $\R \lambda =\gamma$.
We may then conclude equivalence by a limiting argument,
using continuity with respect to $f$ of the left-hand side
and equivalence in the previous case $f\in L^1([0,T]; \CalD(L))$. 
A similar estimate together with expansion \eqref{expansion}
yields pointwise convergence for $f\in L^q([0,T]; \CalD(L))$,
$q>1$, using H\"older's inequality 
$|\lambda^{-1}\widehat{g^T}|_{L^2(\lambda)}
\le C|\lambda^{-1}_{L^q}|\widehat{g^T}|_{L^p}$ and
the Hausdorff--Young inequality
$|\widehat{g^T}|_{L^p}\le C|g|_{L^2[0,T]}$, where $1/p+1/q=1$.
\end{proof}

\begin{ex}
\label{misc}
Let $L:\CalD(L)\to X$ be a closed, densely defined operator
on Banach space $X$, satisfying resolvent bound \eqref{resk}
(equivalently, generating a $C^0$ semigroup $|e^{L t}|\le Ce^{\gamma_0 t}$). 

1.  
For $f\in \CalD(L)$,
show by judicious deformation of the contour that 
$\pv \int_{\gamma -i\infty}^{\gamma+i\infty}
e^{\lambda t}
(\lambda- L)^{-1}f \, d\lambda = 0 $
for any $\gamma>\gamma_0$ and $t<0$,
with uniform pointwise convergence
on compact intervals $[-\epsilon, -\epsilon^{-1}]$, $\epsilon>0$,
at a rate depending only on the bound for $|f|_X+|LF|_{X}$.
Indeed, this holds under the resolvent growth bound
$|(\lambda-L)^{-1}|\le C(1+|\lambda|)^{r}$ for $\R \lambda \ge \gamma$,
for any $r<1$.

2. For $f\in \CalD(L)$, show 
without reference to semigroup theory that
the righthand side of \eqref{basicILT} converges in $L^2(x,t)$ to
a weak solution in the sense of [Sm] of initial value problem
$u_t=Lu$, $u(0)=f$,
using Parseval's identity and resolvent identity \eqref{residentbd}
together with the fact that distributional derivatives and limits commute.
(By moving $\gamma$ to $\omega>0$, show that the $L^2(t;X)$ norm of
$e^{-\omega t}(e^{Lt}f-f)=\int e^{(\lambda-\omega) t}\lambda^{-1}(\lambda-L)^{-1}Lf\,d\lambda$
is order $(\omega-\gamma_0)^{-1}$, yielding a trace at $t=0$.)
Indeed, this holds under the resolvent growth bound
$|(\lambda-L)^{-1}|\le C(1+|\lambda|)^{r}$ for $\R \lambda=\gamma$,
for any $r<1/2$.
\end{ex}

\begin{ex}\label{conv}
Show that 
$\pv \int_{\gamma-i\infty}^{\gamma+i\infty} z^{-1} e^z \, dz$
is uniformly bounded, independent of $\gamma$, and converges uniformly to
$\sgn(\gamma)$ for $\gamma$ bounded away from zero.
\end{ex}

\begin{defi}
\label{evolutiondef}
\textup{
A $C^0$ evolutionary system on Banach space $X$ is a family of bounded
operators $U(s,t)$, $s\le t$, satisfying the properties
(i) $U(s,s)=I$, (ii) $U(s,t+\tau)=U(s,t)U(t,t+\tau)$ for every
$s\le t$, $\tau \ge 0$, and (iii) $\lim_{t\to s^+} U(s,t)x=x$
for all $x\in X$. 
The instantaneous generators $L(s)$ of the system
are defined as $L(s)x=\lim_{t \to s^+}(U(s,t)x-x)/(t-s)$ on the 
domain $\CalD(L(s))\subset X$ for which the limit exists.
}
\end{defi}

\begin{rem}
\label{evojust}
\textup{For a $C^0$ evolutionary system, 
$(d/dt)U(s,s)f=L(s)U(s,s)f$
for all $f\in  \CalD(L(s))$; see [Pa], p. 129.
Thus, if there is a common dense subspace $Y\subset\CalD(L(t))$ for
all $t$ that it is invariant under the flow $U(s,t)$, then, for $f\in Y$,
$U(s,t)$ is a solution operator for 
initial-value problem $u_t=L(t)u$, $u(s)=f$, 
generalizing the notion of semigroup to the nonautonomous case.
}
\end{rem}

\begin{prop}
\label{evolutionprop}
Let $L(s): \CalD(L(s))\to X$ be a family of
closed, densely defined operators on Banach space $X$,
each generating a contraction semigroup $|e^{L(s)t}|\le e^{\gamma_0(s)t}$.
Suppose also that there exists a Banach space $Y\subset \cap_s \CalD(L(s))$
contained in the common intersection of their domains, dense in $X$,
invariant under the flow of each semigroup $e^{L(s)t}$, on which $e^{L(s)t}$
generates a semigroup also in the $Y$-norm, and for which
$L(s):Y\to X$ is continuous in $s$ with respect to the $X\to Y$ operator
norm.
Then, $L(s)$ generates a $C^0$ evolutionary system on $X$ with
domains $\CalD(L(s))$.
\end{prop}

\begin{proof}[Sketch of proof]
Substituting for $L(t)$ the piecewise constant approximations
$L^n(t):= L(k/n)$ for $t\in [k/n,(k+1)/n)$, we obtain a sequence
of approximate evolutionary systems $U^n(s,t)$ satisfying the uniform bound
$|U(s,t)|\le e^{\gamma_0 (t-s)}$ for $s\le t$.
For $x\in Y\subset \CalD(L(s))$ for all $s$, 
the sequence $U^n(s,t)x$ may be shown to be Cauchy
by a Duhamel argument similar to the one used in the
proof of Proposition \ref{semigpcriteria}.
See Theorem 3.1 of [Pa] for details.
\end{proof}

\begin{rem}
\textup{
There exist more general versions applying to systems 
for which the generators are not contractive, but stability
in this case is difficult to verify.  See [Pa], Section 3, for further
discussion.
}
\end{rem}


\bigbreak
\clearpage

\section{Appendix B.  Proof of Proposition \ref{conststab}}

\begin{ex} [{[Kaw]}]
\label{ccbd}
(Optional)
1. Using \eqref{multidiss}, show that 
\begin{equation}
\label{symbd}
|e^{P_-(i\xi)t}|
\le Ce^{-\theta |\xi|^2t/(1+|\xi|^2)},
\end{equation}
where $P_-(i\xi):=
-i\sum_j i\xi_j A^j -\sum_{j,k}\xi_i \xi_j B^{jk}\big)_-$.
\smallbreak
2. Denoting by $S(t)$ the solution operator of the constant-coefficient
problem
\begin{equation}
U_t + \sum_j A^j_- U_{x_j}= \sum_{j,k}(B^{jk}_- U_{x_k})_{x_j},
\end{equation}
we have $\widehat{S(t)U_0}= e^{P_-(i\xi)t}\hat U_0$, where
$\hat{}$ denotes Fourier transform in $x$.
Decomposing 
\begin{equation}
\label{ccdecomp}
S=S_1+ S_2, 
\end{equation}
where
$\widehat{S_j(t)U_0}= \chi_j(\xi) e^{P_-(i\xi)t}\hat U_0$, where
$\chi_1(\xi)$ is one for $|\xi|\le 0$ and zero otherwise, and
$\chi_2=1-\chi_1$, show, using Parseval's identity
$|f|_{L^2(x)}=|\hat f|_{L^2(\xi}$ and the Hausdorff--Young inequality
$|\hat f|_{L^\infty(\xi)} \le |f|_{L^1(x)}$, that
\begin{equation}
\label{s1bd}
|S_1(t)\partial_x^\alpha f|_{L^2}\le C(1+t)^{-d/4- |\alpha|/2}|f|_{L^1},
\quad |\alpha|\ge 0,
\end{equation}
and
\begin{equation}
\label{s2bd}
|S_2(t)f|_{L^2}\le Ce^{-\theta t}|f|_{L^2}.
\end{equation}
\end{ex}

\begin{proof}[Proof of Proposition \ref{conststab}]
Taylor expanding about $U=U_-$, we may rewrite \eqref{viscous} as
\begin{equation}
U_t + \sum_j A^j_- U_{x_j}= \sum_{j,k}(B^{jk}_- U_{x_k})_{x_j}
=
\partial_x Q(U,\partial_x U),
\end{equation}
where
$|Q(U,\partial_x U|\le C|U||\partial_x U|$ and
$|\partial_x Q(U,\partial_x U|\le C(|U||\partial_x^2 U|+ |\partial_x U|^2)$
so long as $|U|$ remains uniformly bounded.
Using Duhamel's principle/variation of constants, we may thus express
the solution of \eqref{viscous} as
\begin{equation}
\label{constduhamel}
U(t)= S(t)U_0 + \int_0^t S(t-s) \partial_x Q(U, \partial U_x)(s)\, ds,
\end{equation}
where $S(\cdot)=S_1(\cdot)+ S_2(\cdot)$ 
is the solution operator discussed in Exercise \ref{ccbd}.

Defining
\begin{equation}
\label{constzeta}
\zeta(t):= \sup_{0\le s\le t} |U(s)|_{H^s}(1+s)^{d/4},
\end{equation}
and using \eqref{s1bd} and \eqref{s2bd},
we may thus bound
\begin{equation}
\begin{aligned}
|U(t)|_{L^2} &\le  |S_1(t)U_0|_{L^2} + 
|S_2(t)U_0|_{L^2} + 
\int_0^t |S_1(t-s) \partial_x Q(U, \partial U_x)(s)|_{L^2}\, ds\\
&\quad + \int_0^t |S_2(t-s) \partial_x Q(U, \partial U_x)(s)|_{L^2}\, ds\\
&\le
C(1+t)^{-d/4}|U_0|_{L^1} + 
Ce^{-\theta t}|U_0|_{L^2} \\
&\quad + 
C\int_0^t (1+ t-s)^{-d/4-1/2} |Q(U, \partial U_x)(s)|_{L^1}\, ds\\
&\quad + C\int_0^t e^{-\theta (t-s)} |\partial_x Q(U, \partial U_x)(s)|_{L^2}\, ds\\
&\le 
C((1+t)^{-d/4}|U(0)|_{L^1\cap L^2} \\
&\quad +
C\zeta(t)^2 \int_0^t (1+t-s)^{-d/4-1/2}(1+s)^{-d/2}\, ds \\
&\le C(1+t)^{-d/4}(|U(0)|_{L^1\cap L^2}+ \zeta(t)^2)
\end{aligned}
\end{equation}
so long as $|U|_{H^s}$ remains uniformly bounded.
Applying now \eqref{ccexp}, we obtain (exercise)
\begin{equation}
\begin{aligned}
|U(t)|_{H^s}^2
&\le Ce^{-\theta t}|U(0)|_{H^s}^2
+ C\int_0^t e^{-\theta(t-s)}|U(s)|_{L^2}^2\, ds \\
&\le Ce^{-\theta t}|U(0)|_{H^s}^2
+ C(1+t)^{-d/2}(|U(0)|_{L^1\cap L^2}+ \zeta(t)^2)^2 \\
&\le
C(1+t)^{-d/2}(|U(0)|_{L^1\cap H^s}+ \zeta(t)^2)^2, \\
\end{aligned}
\end{equation}
and thus
\begin{equation}
\label{ccinduct}
\zeta(t)\le C( |U(0)|_{L^1\cap H^s}+ \zeta(t)^2).
\end{equation}
From \eqref{ccinduct}, it follows by continuous induction (exercise)
that 
\begin{equation}
\zeta(t) \le 2C |U(0)|_{L^1\cap H^s} 
\end{equation}
for $|U(0)|_{L^1\cap H^s}$ sufficiently small, and thus
\begin{equation}
|U(t)|_{H^s}\le 2C(1+t)^{-d/4} |U(0)|_{L^1\cap H^s}
\end{equation}
as claimed.
Applying \eqref{sobolev}, we obtain the same bound for $|U(t)|_{L^\infty}$,
and thus for $|U(t)|_{L^p}$, all $2\le p\le \infty$, by
the $L^p$ interpolation formula
\begin{equation}
\label{Lpinterp}
|f|_{L^{p_*}} \le |f|_{L^{p_1}}^{\beta}|f|_{L^{p_2}}^{1-\beta}
\end{equation}
for all $1\le p_1<p_*<p_2$,
where $\beta:=p_1(p_2-p_*)/p_*(p_2-p_1)$ is
determined by $1/p_*= \beta/p_1 + (1-\beta)/p_2$
(here applied between $L^2$ and $L^\infty$, i.e., with $p_1=2$,
$p_*=p$, $p_2=\infty$).
\end{proof}

\bigbreak
\clearpage

\section{Appendix C.  Proof of Proposition \ref{branch}} 
In this appendix, we complete the proof of Proposition \ref{branch}
in the general {\it symmetrizable, constant-multiplicity} case;
here, we make essential use of  recent results
of M\'etivier [M\'e.4] concerning the
spectral structure of matrix $(A^1)^{-1}(i\tau + iA^\Txi)$.
Without loss of generality, take $A^\xi$ {\it symmetric};
this may be achieved by the change of coordinates 
$A^{\xi}\to \tilde A_0^{1/2}A^\xi \tilde A_0^{-1/2}$.

With these assumptions, the kernel and co-kernel of
$(A^{\xi_0}+ \tau_0)$ are of fixed dimension
$m$, not necessarily equal to one, and are spanned
by a common set of zero-eigenvectors
$r_1,\dots,r_m$.
Vectors $r_1,\dots,r_m$ are necessarily right zero-eigenvectors
of $(A^1)^{-1} (i\tau_0 + iA^{\xi_0})$ as well.  Branch singularities
correspond to the existence of one or more Jordan chains of
generalized zero-eigenvectors extending up from genuine eigenvectors
in their span, which by the argument of Lemma 3.4 is equivalent to
\begin{equation}
\label{vanish2}
\det (r_j^tA^1 r_k)=0.
\end{equation}
In fact, as pointed out by M\'etivier [M\'e.4], the
assumption of constant multiplicity implies
considerable additional structure.

\medbreak
\begin{obs} [{[M\'e.2]}] \label{block} 
\textup{
Let $(\Txi_0, \tau_0)$ lie at a branch singularity
involving root $\alpha_0=i\xi_{0_1}$ in \eqref{slow},
with $\tau_0$ an $m$-fold eigenvalue of $A^{\xi_0}$.
Then, for $(\Txi, \tau)$ in the vicinity of
$(\Txi_0, \tau_0)$,  the roots $\alpha$ bifurcating
from $\alpha_0$ in \eqref{slow} consist of
$m$ copies of $s$ roots $\alpha_1, \dots,\alpha_s$,
where $s$ is some fixed positive integer.
}
\end{obs}

\medbreak
\begin{proof}[Proof of Observation]
Let $a(\Txi, \alpha)$ denote the unique eigenvalue
of $A^\xi$ lying near $-\tau_0$, where, as usual, $-i\xi_1:= \alpha$;
by the constant multiplicity assumption, $a(\cdot,\cdot)$ is
an analytic function of its arguments.
Observing that
\begin{equation}
\begin{aligned}
\det[ (A^1)^{-1} (i\tau + iA^\Txi)-\alpha]
&= \det i(A^1)^{-1}
\det (\tau + A^\xi)\\
&=
e(\Txi,\tau,\alpha)(\tau + a(\Txi,\alpha))^m, 
\end{aligned}
\end{equation}
where $e(\cdot,\cdot,\cdot)$ does not vanish
for $(\Txi,\tau,\alpha)$ sufficiently close to
$(\Txi_0,\tau_0,\alpha_0)$,
we see that the roots in question occur as $m$-fold
copies of the roots of 
\begin{equation}
\label{key}
\tau+ a(\Txi,\alpha)=0.
\end{equation}
But, the lefthand side of \eqref{key} is a family
of analytic functions in $\alpha$, continuously varying
in the parameters $(\Txi,\tau)$, whence the number of
zeroes is constant near $(\Txi_0,\tau_0)$.
\end{proof}
\medbreak

\begin{obs} [{[Z.3]}] \label{ident}
\textup{
The matrix $(r_j^tA^1 r_k)$, $j$, $k=1,\dots,m$ is a 
real multiple of the identity,
\begin{equation}
\label{ident1}
(r_j^tA^1 r_k)= (\partial a/\partial \xi_1)I_m,
\end{equation}
where $a(\xi)$ denotes the (unique, analytic) $m$-fold eigenvalue of $A^\xi$
perturbing from $-\tau_0$.
}

\textup{
More generally, if 
\begin{equation}
\label{szeroes}
(\partial a/\partial \xi_1)= \cdots = 
(\partial^{s-1} a/\partial \xi_1^{s-1})= 0,
\quad
(\partial^{s} a/\partial \xi_1^{s})\ne 0
\end{equation}
at $\xi_0$, then, letting
$r_1(\Txi),\dots,r_m(\Txi)$
denote an analytic choice of basis for the eigenspace 
corresponding to $a(\Txi)$, orthonormal at $(\Txi_0,\tau_0)$, 
we have the relations
\begin{equation}
\label{chain}
(A^1)^{-1}(\tau_0+ A^{\xi_0})r_{j,p}=
 r_{j,{p-1}},
\end{equation}
for $1\le p\le s-1$, 
and
\begin{equation}
\label{idents}
(r_{j,0}^t
A^1 r_{k,p-1})=
p!\,
(\partial^{p} a/\partial \xi_1^{p}) I_m,
\end{equation}
for $1\le p\le s$, 
where
\begin{equation}
r_{j,p}:=
(-1)^{p}p(\partial^{p} r_j/\partial \xi_1^{p}).
\end{equation}
In particular, 
\begin{equation}
\{r_{j,0},\dots, r_{j,s-1}\}, \quad j=1,\dots,m
\end{equation}
is a right Jordan basis for the total zero eigenspace
of $(A^1)^{-1}(\tau_0 + A^{\xi_0})$,
for which the genuine zero-eigenvectors $\tilde l_j$ of the dual, 
left basis are given by
\begin{equation}
\label{corners}
\tilde l_j=(1/s! (\partial^{s} a/\partial \xi_1^{s}) )A^1 r_j.
\end{equation}
}
\end{obs}

\medbreak
\begin{proof}[Proof of Observation]
Considering $A^\xi$ as a matrix perturbation in $\xi_1$,
we find by standard spectral perturbation theory that
the bifurcation of the $m$-fold eigenvalue 
$\tau_0$ as $\xi_1$ is varied is governed to first order
by the spectrum of $(r_j^tA^1 r_k)$.  Since these eigenvalues
in fact do not split, it follows that 
$(r_j^tA^1 r_k)$ has a single eigenvalue.  
But, also, $(r_j^tA^1 r_k)$ is symmetric, hence diagonalizable,
whence we obtain result \eqref{ident1}.

Result \eqref{idents} may be obtained by a more systematic version of
the same argument. Let $R(\xi_1)$ denote 
the matrix of right eigenvectors
\begin{equation}
R(\xi_1):=(r_1,\dots,r_m)(\xi_1).
\end{equation}
Denoting by
\begin{equation}
\label{aexp}
a(\xi_1+h)=: a^0 + a^1 h + \dots + a^p h^p + \dots
\end{equation}
and
\begin{equation}
\label{Rexp}
R(\xi_1+h)=: R^0 + R^1 h + \dots + R^p h^p + \dots
\end{equation}
the Taylor expansions of functions $a(\cdot)$ and $R(\cdot)$
around $\xi_0$ as $\xi_1$ is varied, and recalling that
\begin{equation}
A^\xi= A^{\xi_0}+ hA^1,
\end{equation}
we obtain in the usual way, matching terms of common order in the 
expansion of the defining relation $(A-a)R=0$, the heirarchy
of relations:
\begin{equation}
\label{heirarchy}
\begin{aligned}
(A^{\xi_0}-a^0)R^0&=0, \\
(A^{\xi_0}-a^0)R^1&=-(A^1-a^1)R^0, \\
(A^{\xi_0}-a^0)R^2&=-(A^1-a^1)R^1 + a^2R^0, \\
&\vdots\\
(A^{\xi_0}-a^0)R^p&=-(A^1-a^1)R^{p-1} + a^2R^{p-2}+\dots + a^pR^0. \\
\end{aligned}
\end{equation}
Using $a^0=\cdots=a^{s-1}=0$, we obtain \eqref{chain} immediately,
from equations $p=1,\dots,s-1,$
and $R^p=(1/p!)(r_{1,p}, \dots, r_{m,p})$.
Likewise, \eqref{idents}, follows from equations $p=1,\dots,s$,
upon left multiplication by $L^0:= (R^0)^{-1}= (R^0)^t$,
using relations
$L^0(A^{\xi_0}-a^0)=0$ and $a^p= (\partial^p a/\partial \xi_1^p)/p!$.

From \eqref{chain}, we have the claimed right Jordan basis.
But, defining $\tilde l_j$ as in \eqref{corners}, we can rewrite
\eqref{idents} as
\begin{equation}
\label{newidents}
(\tilde l_j^t r_{k,p-1})=
\begin{cases}
0 & 1\le p \le s-1, \\
I_m, & p=s; \\
\end{cases}
\end{equation}
these $ms$ criteria uniquely define
$\tilde l_j$ (within the $ms$-dimensional total left eigenspace)
as the genuine left eigenvectors dual to the right basis formed 
by vectors $r_{j,p}$ (see also exercise just below). 
\end{proof}
\medbreak

Observation \ref{ident} implies in particular 
that Jordan chains extend from {\it all} or {\it none} 
of the genuine eigenvectors
$r_1,\dots,r_m$, with common height $s$.
As suggested by Observation \ref{block}
(but not directly shown here),
this uniform structure in fact persists under variations
in $\Txi$, $\tau$, see [M\'e.4].
Observation \ref{ident} is a slightly more
concrete version of Lemma 2.5 in [M\'e.4];
note the close similarity between the argument used here,
based on successive variations in basis $r_j$,
and the argument of [M\'e.4], based on variations 
in the associated total projection.

\medbreak

With these preparations, the result goes through essentially
as in the strictly hyperbolic case.
Set
\begin{equation}
\label{p}
p:=
1/(s! (\partial^{s} a/\partial \xi_1^{s}))
\end{equation}
and define
\begin{equation}
\label{q}
pR^t B{\Txi_0,\Txi_0} R =: Q
\end{equation}
Note, as claimed, that $p\ne 0$ by assumption
$(\partial^s a/\partial \xi_1^s)\ne 0$
in Observation \ref{ident}, and $sgn (p) Q >0$ by (K1), Proposition \ref{KSh}.

Thus, working in the Jordan basis defined in Observation \ref{ident},
we find similarly as in the strictly hyperbolic
case that the matrix perturbation problem
\eqref{modsuperslow} reduces to an $ms\times ms$ block-version 
\begin{equation}
\label{syscanonical}
\big(iJ+i\sigma M + \rho N-(\tilde \alpha-\alpha)\big)\V_{I_j}= 0
\end{equation}
of \eqref{canonical} in the strictly hyperbolic case, where
\begin{equation}
\label{blockjordan}
J:=
\begin{pmatrix}
0 & I_m & 0 & \cdots &0 \\
0 & 0 & I_m & 0 & \cdots \\
0 & 0 & 0 & I_m & \cdots \\
\vdots &  \vdots & \vdots &  \vdots & \vdots\\
0 & 0 & 0 & \cdots&0 \\
\end{pmatrix},
\end{equation}
denotes the standard block-Jordan block,
and the lower-lefthand block of $i\sigma M+\rho N$
is $\sigma p I_m -i\rho Q \sim |\sigma|+|\rho|$.
To lowest order $O(|\sigma|+|\rho|)^{1/s}$, therefore, the problem
reduces to the computation of eigenvectors and
eigenvalues of the perturbed block-Jordan block
\begin{equation}
\label{syscomputation}
diag\big\{
i\begin{pmatrix}
0 & I_m & 0 & \cdots &0 \\
0 & 0 & I_m & 0 & \cdots \\
0 & 0 & 0 & I_m & \cdots \\
\vdots &  \vdots & \vdots &  \vdots & \vdots\\
\sigma p I_m -i\rho Q & 0 & 0 & \cdots&0 \\
\end{pmatrix}
\big\},
\end{equation}
from which results \eqref{bevalues}--\eqref{bcondition} follow as before
by standard matrix perturbation theory;
see, e.g.,  Section 2.2.4, {\it Splitting of a block-Jordan block} of [Z.4].
(Note that the simple eigenvalue case $s=1$ 
follows as a special case of the block-Jordan block computation,
with $\sigma \equiv 0$.)
This completes the proof of Proposition \ref{branch} in the general case.
\bigbreak

\end{document}